\documentclass[a4paper,11pt]{article}
\usepackage[top=3cm,bottom=3cm,left=2.5cm,right=2.5cm,a4paper]{geometry}
\usepackage{amsfonts,amsmath,amssymb,amsthm}
\usepackage{upgreek}
\usepackage{float,xcolor,fancyhdr}
\usepackage[T1]{fontenc}
\usepackage{algorithmic}
\usepackage[ruled,vlined]{algorithm2e}
\usepackage{setspace}
\usepackage{lineno}
\usepackage{longtable} 
\usepackage[normalem]{ulem} 
\usepackage[colorlinks]{hyperref}
\hypersetup{
    colorlinks=true, 
    linkcolor={blue!50!black},
    citecolor={green!60!black},
    urlcolor= {red!80!black}
}
\usepackage{enumitem,mdframed}
\usepackage{empheq}
\usepackage[bottom]{footmisc} 
\usepackage{tikz,pgfplots}
\usepgfplotslibrary{groupplots}
\usetikzlibrary{positioning,arrows,fit,patterns}
\pgfplotsset{compat=newest}
\usepackage{siunitx}
\usepackage{subfig}
\DeclareCaptionLabelFormat{opening}{\textbf{#2)}} 
\makeatletter
\renewcommand*\env@matrix[1][\arraystretch]{%
  \edef\arraystretch{#1}%
  \hskip -\arraycolsep
  \let\@ifnextchar\new@ifnextchar
  \array{*\c@MaxMatrixCols c}}
\makeatother

\usepackage[]{todonotes}

\numberwithin{equation}{section}

\newcommand*\samethanks[1][\value{footnote}]{\footnotemark[#1]}

\newcommand{\myblue} {blue!80!white}
\newcommand{\mygreen}{green!50!black}

\newcommand{\forward}     {\mathcal{F}}
\newcommand{\R}           {\mathbb{R}}
\newcommand{\Real}        {\mathrm{Re}}
\newcommand{\Imag}        {\mathrm{Im}}
\newcommand{\Cx}          {\mathbb{C}}

\newcommand{\ii}          {\mathrm{i}}

\newcommand{\bx}          {\mathbf{x}}
\newcommand{\misfit}      {\mathcal{J}}

\newcommand{\data}        {\boldsymbol{d}}
\newcommand{\bm}          {\boldsymbol{m}}
\newcommand{\nrcv}        {n_{\mathrm{rcv}}}
\newcommand{\nsrc}        {n_{\mathrm{src}}}
\newcommand{\velocity}    {\boldsymbol{v}}
\newcommand{\pressure}    {p}

\newcommand{\n}           {\boldsymbol{n}} 

\newcommand{\tausigma}    {\tau_{\sigma}}
\newcommand{\taueps}      {\tau_{\epsilon}}

\newcommand{\omegaR}{\omega_{\mathrm{R}}}
\newcommand{\omegaI}{\omega_{\mathrm{I}}}
\newcommand{\cR}{c_{\mathrm{R}}}
\newcommand{\cI}{c_{\mathrm{I}}}
\newlength{\modelwidth} \newlength{\modelheight}
\newlength{\plotwidth} \newlength{\plotheight}

\newcommand{\modelfile} {}
\newcommand{\modelfileCp} {}
\newcommand{\modelfileRho} {}
\newcommand{\datafile}  {}
\newcommand{\dataA}{}   \newcommand{\dataB}{} \newcommand{\dataC}{}
\newcommand{\dataD}{}   \newcommand{\dataE}{} \newcommand{\dataF}{} 
\newcommand{\dataG}{}   \newcommand{\dataH}{} 
\newcommand{\legendA}{} \newcommand{\legendB}{} \newcommand{\legendC}{}
\newcommand{\legendD}{} \newcommand{\legendE}{} \newcommand{\legendF}{} 
\newcommand{\legendG}{} \newcommand{\legendH}{} 

\newlength {\jumpvert} 

\newtheorem{definition} {Definition}

\newtheorem{remark}     {Remark}

\newtheorem{assumption} {Assumption}

\usepackage[capitalize,nameinlink]{cleveref}
\crefname{section}   {Section}   {Sections}
\crefname{subsection}{Subsection}{Subsections}
\Crefname{section}   {Section}   {Sections}
\Crefname{subsection}{Subsection}{Subsections}
\Crefname{figure}    {Figure}    {Figures}

\crefformat{equation}{\textup{#2(#1)#3}}
\crefrangeformat{equation}{\textup{#3(#1)#4--#5(#2)#6}}
\crefmultiformat{equation}{\textup{#2(#1)#3}}{ and \textup{#2(#1)#3}}
{, \textup{#2(#1)#3}}{, and \textup{#2(#1)#3}}
\crefrangemultiformat{equation}{\textup{#3(#1)#4--#5(#2)#6}}%
{ and \textup{#3(#1)#4--#5(#2)#6}}{, \textup{#3(#1)#4--#5(#2)#6}}{, and \textup{#3(#1)#4--#5(#2)#6}}
\Crefformat{equation}{#2Equation~\textup{(#1)}#3}
\Crefrangeformat{equation}{Equations~\textup{#3(#1)#4--#5(#2)#6}}
\Crefmultiformat{equation}{Equations~\textup{#2(#1)#3}}{ and \textup{#2(#1)#3}}
{, \textup{#2(#1)#3}}{, and \textup{#2(#1)#3}}
\Crefrangemultiformat{equation}{Equations~\textup{#3(#1)#4--#5(#2)#6}}%
{ and \textup{#3(#1)#4--#5(#2)#6}}{, \textup{#3(#1)#4--#5(#2)#6}}{, and \textup{#3(#1)#4--#5(#2)#6}}

\crefdefaultlabelformat{#2\textup{#1}#3}

\crefname{proposition}{Proposition}{Propositions}
\Crefname{proposition}{Proposition}{Propositions}
\crefname{definition} {Definition} {Definitions}
\Crefname{definition} {Definition} {Definitions}
\crefname{theorem}    {Theorem}    {Theorems}
\Crefname{theorem}    {Theorem}    {Theorems}
\crefname{remark}     {Remark}     {Remarks}
\Crefname{remark}     {Remark}     {Remarks}
\crefname{assumption} {Assumption} {Assumptions}
\Crefname{assumption} {Assumption} {Assumptions}

\title{Quantitative inverse problem in 
       visco-acoustic media under attenuation 
       model uncertainty}

\author{
Florian Faucher\thanks{Faculty of Mathematics, University of Vienna, Oskar-Morgenstern-Platz 1,
                       A-1090 Vienna, Austria.}
                \thanks{Project-Team Makutu, Inria Bordeaux Sud-Ouest, Universit\'e de Pau et 
                        des Pays de l'Adour, France, 
                       \href{mailto:florian.faucher@inria.fr}
                      {\texttt{florian.faucher@inria.fr}}.}
\and
Otmar Scherzer\samethanks[1]
\thanks{Johann Radon Institute for Computational and Applied Mathematics (RICAM), Altenbergerstraße 69, A-4040 Linz, Austria}
\thanks{Christian Doppler Laboratory for Mathematical Modeling and Simulation of Next Generations of Ultrasound Devices (MaMSi), Oskar-Morgenstern-Platz 1, A-1090 Vienna, Austria, \href{mailto:otmar.scherzer@univie.ac.at}{\texttt{otmar.scherzer@univie.ac.at}}}
}
\date{\today}

\begin{document}
\maketitle 

\begin{abstract}
  We consider the inverse problem of quantitative 
  reconstruction of properties (e.g., bulk modulus, density) 
  of visco-acoustic materials based on measurements 
  of responding waves after stimulation of the medium.
  Numerical reconstruction 
  is performed by an iterative minimization algorithm.
  Firstly, we investigate the robustness of the algorithm 
  with respect to attenuation model uncertainty, that is, 
  when different attenuation models are used to simulate
  synthetic observation data and for the inversion, respectively.
  Secondly, to handle data-sets with multiple reflections 
  generated by wall boundaries around the domain, we perform 
  inversion using complex frequencies, and show that it 
  offers a robust framework that alleviates the difficulties
  of multiple reflections.
  To illustrate the efficiency of the algorithm, we
  perform numerical simulations of ultrasound imaging 
  experiments to reconstruct a synthetic 
  breast sample that contains an inclusion of 
  high-contrast properties.
  We perform experiments in two and three dimensions, 
  where the latter
  also serves to demonstrate the numerical feasibility
  in a large-scale configuration. 
\end{abstract}

\section{Introduction}

Attenuation refers to the loss of energy of waves propagating in a medium. 
The level of energy loss is quantified as a material property of a viscous medium.
In the context of imaging with waves, visco-acoustic media represent additional 
challenges compared to an ideal (that is a non-attenuating) medium:
For instance, the mathematical description of such materials involve 
more parameters, e.g., to account for the attenuation mechanisms.
In this work we investigate the quantitative reconstruction of 
properties of a visco-acoustic medium, which is relevant, for instance, 
in ultrasound imaging, non-destructive testing, or seismic imaging, 
e.g., \cite{Donald1958,Fenster1998,Virieux2009,Kamei2013,Beniwal2015,Perez2017,Bachmann2020}. 
In such a setup, probing waves are sent to the 
area of interest, and the mechanical response of the medium is indirectly measured at the position 
of the receivers.

Attenuation is a phenomenon that intrinsically depends 
on the frequencies of the propagating waves, where each frequency component 
looses a different amount of energy 
\cite{Bland1960,AkiRichards,Carcione2007}.
It is common to categorize attenuation models 
into different families, such as those that dominantly attenuate
low or high frequency components, see \cref{section:forward}.
\emph{Generalized models} combine multiple mechanisms of attenuation 
and enable the consideration of signals with wide frequency ranges, 
\cite{Carcione2007,BrossierVirieux2016}, see \cref{rk:generalized}.
Therefore, the wave equations 
that model the propagation of attenuated waves
can have different forms, as introduced, for instance, 
in \cite{Ursin2002,Carcione2007}.
From a practical perspective, it is difficult to 
identify the appropriate mathematical model of attenuation that adequately  
describes a given material, cf. \cite{Carcione2007,Royston2011} 
for geophysical applications or \cite{Zhou2018} 
for medical ones. 
Consequently, in practical applications of inverse problems, 
the appropriate attenuation model to describe the material 
behaviour represents an additional unknown. 
Note that despite the different nature of the media
and of the acquisition setups, cf.~\cite{Pratt2007},
geophysical and medical imaging share the same 
models of visco-acoustic wave propagation. 

In our work we assume that the excitation 
that generates the waves within the medium is a 
broad-band time-signal.
Then, we apply a Fourier transform to the time 
measurements to work with frequency-domain data, 
as we illustrate in \cref{subsection:time-to-frequency}.
To simulate the waves that will later
be compared with the measurement data used to reconstruct the 
medium's properties, we use the frequency-domain wave
propagation equations. 
The frequency-domain formulation has been recognized 
to be more convenient to handle attenuation, cf.~\cite{Bland1960,Muller1983,Pratt1999b},
as it allows us to work with Partial Differential Equations 
(PDEs) regardless of the attenuation model. On the other hand, 
the time-domain formulation can lead to integro-differential equations instead. 
This is due to the fact that in the frequency-domain the attenuation is encoded via 
complex-valued parameters in the PDEs, e.g., \cite{Muller1983,Carcione1988,Ursin2002}.
However, causality principles must be carefully 
addressed in the frequency-domain as advocated 
by \cite{Golden1988,Carcione2007,Elbau2017}.
The level of attenuation of a medium is 
quantified by the frequency-dependent \emph{quality factor}, see \cref{definition:Q}.
The higher it is, the less attenuating is the medium.
In our experiments, we consider the reconstruction 
of breast tissues at ultrasound frequencies, where the  
quality factor is higher than 100 
(see \cref{table:interval}, \cref{section:fwi-2d,section:fwi-3d}): this corresponds to 
weakly attenuating  medium, \cite{Vavryvcuk2009}.

For the quantitative reconstruction of medium properties, 
our experiments follow the principles of 
ultrasound tomography, which have been applied, 
for instance, in the context of breast imaging, 
e.g.,  \cite{Duric2005,Duric2007,Li2009}, 
where simplifying assumptions regarding the 
models of wave propagation can be used. Reconstruction
methods are further compared in \cite{Ozmen2015,Faucher2022Book}.
We further use the full model of wave propagation 
for the simulations, and the reconstructions are carried 
out with an iterative method that minimizes a misfit function 
defined to evaluate the distance between the observed 
measurements and simulations.
This approach is usually referred to as the \emph{Full Waveform Inversion} 
(FWI) in seismic imaging, \cite{Lailly1983,Tarantola1984,Pratt1998,Virieux2009}. 
As its name indicates, FWI relies on the full 
measurement signals, e.g., contrary to approaches 
that only use the travel-times.

Imaging using the full waveform has been used in 
the context of ultrasound tomography, for instance with special emphasis on 
bone structures, \cite{Komatitsch2017}, for the brain, 
\cite{Guasch2020} and for breast imaging considering the breast as 
a viscous medium, with implementations both in the 
frequency, \cite{Li2014}, and time-domain, 
\cite{Pratt2007,Perez2017,Bachmann2020,Cox2021}. 
We can highlight two major differences between
our work and these references: 
\begin{enumerate}
	\item The aforementioned references rely on a single 
	      attenuation model, which is fixed a-priori, prior to the 
	      reconstruction. 
	\item These papers assume wave 
	propagation in free-space, that is, assuming that no 
	reflection comes from the boundary of the acquisition 
	setup. 
	For practical experiments this requires particular experimental 
        conditions, such as adding padding materials, 
	or a non-trivial data post-processing step (\cite{Yilmaz2001}) 
	to remove the multiple reflections coming from 
	imperfect boundaries while preserving the ones that come from the sample.
\end{enumerate}
More specifically, the contributions of our work are the following:
\vspace*{-0.5em}
\begin{itemize}\setlength{\itemsep}{-2pt}
  \item We review and implement seven different attenuation  
        models describing visco-acoustic wave propagation.
        Each model encompasses a different 
        effect of wave dissipation and dispersion.
  \item We carry out inversion with \emph{attenuation model uncertainty}, 
        that is, we use a different attenuation model to generate the synthetic data (simulating measurement data) and 
        to carry out the numerical reconstruction procedure. 
        Despite the resulting changes in the forward PDE models, we
        show that FWI is a robust approach that does not suffer from 
        inconsistency in the attenuation model.
  \item We start with experiments where we impose 
        \emph{absorbing boundary conditions} to constrain the 
        numerical domain, hence mimicking a free-space 
        wave propagation.
        We then investigate the consideration 
        of wall boundary on the sides of the domain. 
        The resulting multiple reflections are shown to 
        strongly influence the accuracy of the reconstructions.
        To overcome the difficulty, we use complex frequencies.
        This approach is also referred to as the \emph{Laplace-Fourier} 
        domain method, \cite{Shin2008,Shin2009} and it is shown to improve 
        the convexity of the misfit function for inversion 
        in \cite{Faucher2020basins,Faucher2017}.
        In our work, we show that this transformation can alleviate the 
        difficulty that occur from the multiple reflections coming 
        from wall boundaries, by enhancing wave first arrivals, \cite{Kamei2013}.
  \item In the context of multi-parameter inversion, we investigate 
        the choice of parametrization, that is, the choice of model
        parameters with respect to whom the gradient is computed, 
        \cite{Brossier2011,Kohn2012,Faucher2017}.
        In particular, the density and attenuation properties 
        are known to be hard to reconstruct, \cite{Virieux2009,Jeong2012},
        and possibly require a specific misfit function, \cite{Karaouglu2017}.
  \item We perform experiments in three dimensions to explore
        the feasibility of our methodology and detail the 
        computational cost.
\end{itemize}

In \cref{section:forward}, we provide the visco-acoustic
equation for wave propagation in the frequency domain, 
and review seven different attenuation models from the literature, 
\cite{Ursin2002,Carcione2007,Elbau2017}.
We highlight the benefit of working in the frequency domain, 
in which case all the attenuation models are encoded as PDEs
with complex-valued parameters. Therefore the frequency domain allows 
for using PDE software infrastructure, which the time domain does not allow. 
The inverse procedure for the reconstruction of parameters 
using FWI is detailed in \cref{section:inverse-problem}. 
In \cref{section:fwi-2d}, we carry out reconstructions for 
a two-dimensional sample corresponding to a breast model,
while a three-dimensional case is performed in \cref{section:fwi-3d},
where we also provide the computational cost of the method.

\section{Visco-acoustic forward wave problem}
\label{section:forward}

\subsection{From time-domain to complex frequency-domain data-sets}
\label{subsection:time-to-frequency}

We investigate the reconstruction of material parameters 
of a heterogeneous visco-acoustic medium. 
Waves are excited at the boundary of a sample and propagate 
through it; receivers positioned outside of the sample measure 
the response wavefield which is used to reconstruct the 
medium properties. 
In particular, this means that we investigate a non-invasive experiment.

We assume that the exciting source has a compact support.
Mathematically, it can be represented by a delta-Dirac in space
and a time-domain wavelet: $\delta(\bx)\, S(t)$. 
We refer to such a function as a \emph{point-source}, 
see \cref{fig:time-domain-acquisition-cx-freq_a}. 
We can also consider \emph{line-source} excitations, 
see \cref{fig:time-domain-acquisition-cx-freq_t}. 
Mathematically it is written as a composition of 
point-sources that are excited at the same time-instant.
In this work, the time-domain signal/wavelet $S(t)$ of 
the point-source excitation follows a Ricker wavelet function, 
which represents a broad-band signal containing multiple 
frequencies, as illustrated in \cref{fig:time-domain-acquisition-cx-freq_ricker-f}.

Time-domain signals are transformed with a Laplace-Fourier transform 
\cite{Shin2009,Faucher2020basins,Faucher2020adjoint}, introducing the 
complex frequency $\omega$ such that,
\begin{equation} \label{eq:omegaCx} 
	\omega \, := \, \omegaR \,+\, \ii \, \omegaI \,\,\,\,\,\Rightarrow\,\,\,\,\,
	\ii \omega \,=\, \ii \omegaR - \omegaI \, 
	\qquad \text{with } \qquad \omegaR\,>\,0 \quad \text{ and } \quad \omegaI \,\geq\, 0.
\end{equation}
Namely, the discrete transform of a signal $S$ composed 
of $N_t$ time-steps such that $S=\{\,S_0$, $S_1$, $S_2$, $\ldots$, 
$S_{N_t-1} \}$
is computed with
\begin{equation} \label{eq:complex-dft}
  \mathfrak{F}\big(S \big) (\omega) \,=\,
  \sum_{k=0}^{N_t-1} \,\, S_k \, e^{\ii\omega\,k/N_t} \, ,
\end{equation}
using the same convention as \cite[Box~5.2]{AkiRichards}.
Note that when the imaginary part is zero, $\omegaI=0$, 
it corresponds to the ordinary frequency $\nu$ \si{\Hz}, 
with $\omegaR = 2\pi \nu$.
In the experiments, receivers record time signals 
(e.g., pressure fields), which are thus first 
transformed to complex-frequency data, see 
\cref{fig:time-domain-acquisition-cx-freq}.

\setlength{\modelwidth} {5cm}
\setlength{\modelheight}{5cm}
\graphicspath{{figures/tumor_18x18/main/}}
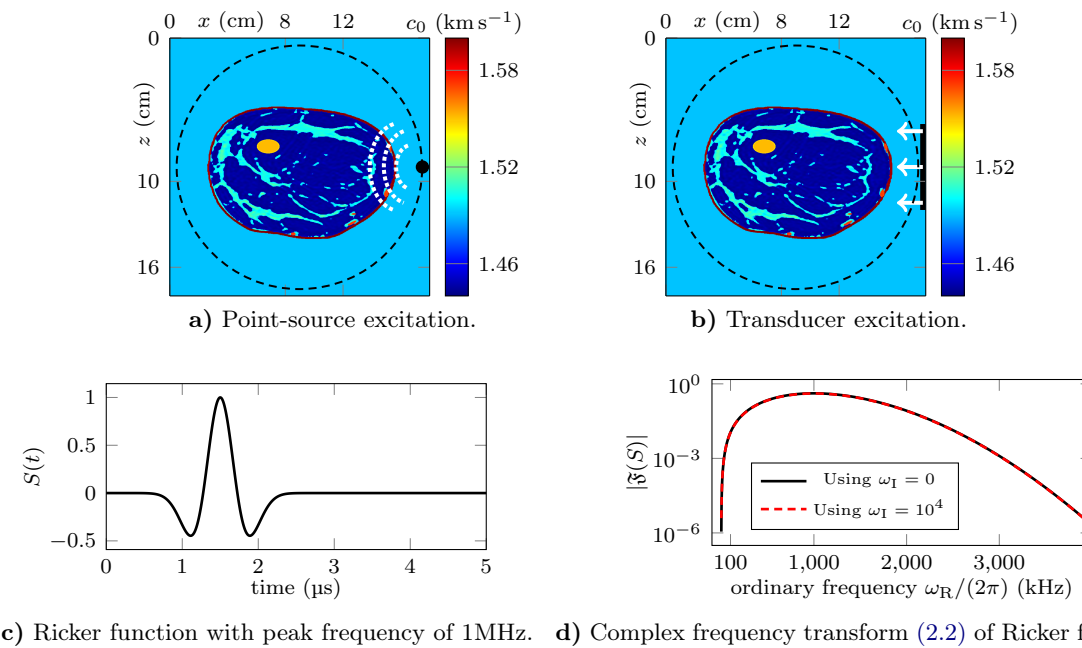
\begin{figure}[ht!] \centering
 \renewcommand{\modelfile}{cp-true_scale1440-1600}
 \subfloat[][Point-source excitation.]
         {\makebox[.40\linewidth]{\begin{tikzpicture}

  \pgfmathsetmacro{\xmin}{0}
  \pgfmathsetmacro{\xmax}{18}
  \pgfmathsetmacro{\zmin}{0}
  \pgfmathsetmacro{\zmax}{18}
  \pgfmathsetmacro{\cmin}{1.44}
  \pgfmathsetmacro{\cmax}{1.60}
  
\begin{axis}[%
width=\modelwidth, height=\modelheight,
axis on top, separate axis lines,
xmin=\xmin, xmax=\xmax, xlabel={$x$ (\si{\cm})},
ymin=\zmin, ymax=\zmax, ylabel={$z$ (\si{\cm})}, y dir=reverse,
xtick={0,8,12},
ytick={0,10,16},
xticklabel pos=right, xlabel near ticks,
x label style={xshift=-0.9cm, yshift=-0.50cm}, 
y label style={xshift= 0.7cm, yshift=-0.50cm},
colormap/jet,colorbar,
colorbar style={title={{\scriptsize{$c_0$ (\si{\km\per\second})}}},
title style={yshift=-3mm, xshift=1mm},
width=3mm,xshift=-1mm,
ytick={1.46,1.52,1.58},
},point meta min=\cmin,point meta max=\cmax,
  label style={font=\scriptsize},
  tick label style={font=\scriptsize},
  legend style={font=\scriptsize\selectfont},
]
\addplot [forget plot] graphics [xmin=\xmin,xmax=\xmax,
                                 ymin=\zmin,ymax=\zmax] {{\modelfile}.png};

   \draw[mark=*,mark size=2pt,mark options={color=black},black,line width=1pt, only marks] plot coordinates{(17.500,   9.0000)};

  \draw[densely dotted, line width=1.5pt, white] (16.50,7.5000) to[out=200,in=160] (16.500,10.50);
  \draw[densely dotted, line width=1.5pt, white] (16.20,6.5000) to[out=200,in=160] (16.200,11.50);  
  \draw[densely dotted, line width=1.5pt, white] (15.50,6.0000) to[out=200,in=160] (15.500,12.00);  
\end{axis}
  
  \draw[black,line width=0.75pt, densely dashed] (1.705,1.705) circle (1.615cm);
\end{tikzpicture}%
}
          \label{fig:time-domain-acquisition-cx-freq_a}}
 \subfloat[][Transducer excitation.]
         {\makebox[.40\linewidth]{\begin{tikzpicture}

  \pgfmathsetmacro{\xmin}{0}
  \pgfmathsetmacro{\xmax}{18}
  \pgfmathsetmacro{\zmin}{0}
  \pgfmathsetmacro{\zmax}{18}
  \pgfmathsetmacro{\cmin}{1.44}
  \pgfmathsetmacro{\cmax}{1.60}
  
\begin{axis}[%
width=\modelwidth, height=\modelheight,
axis on top, separate axis lines,
xmin=\xmin, xmax=\xmax, xlabel={$x$ (\si{\cm})},
ymin=\zmin, ymax=\zmax, ylabel={$z$ (\si{\cm})}, y dir=reverse,
xtick={0,8,12},
ytick={0,10,16},
xticklabel pos=right, xlabel near ticks,
x label style={xshift=-0.9cm, yshift=-0.50cm}, 
y label style={xshift= 0.7cm, yshift=-0.50cm},
colormap/jet,colorbar,
colorbar style={title={{\scriptsize{$c_0$ (\si{\km\per\second})}}},
title style={yshift=-3mm, xshift=1mm},
width=3mm,xshift=-1mm,
ytick={1.46,1.52,1.58},
},point meta min=\cmin,point meta max=\cmax,
  label style={font=\scriptsize},
  tick label style={font=\scriptsize},
  legend style={font=\scriptsize\selectfont},
]
\addplot [forget plot] graphics [xmin=\xmin,xmax=\xmax,
                                 ymin=\zmin,ymax=\zmax] {{\modelfile}.png};

  \draw[line width=2.0pt, black]    (17.80, 6.00) to[] (17.800,12.00);
  \draw[line width=1.5pt, white,->] (17.80, 9.00) to[] (16.000, 9.00);
  \draw[line width=1.5pt, white,->] (17.80,11.50) to[] (16.000,11.50);
  \draw[line width=1.5pt, white,->] (17.80, 6.50) to[] (16.000, 6.50);

\end{axis}
  
  \draw[black,line width=0.75pt, densely dashed] (1.705,1.705) circle (1.615cm);
\end{tikzpicture}%
}
          \label{fig:time-domain-acquisition-cx-freq_t}}

 \renewcommand{\datafile}{figures/tumor_18x18/time-domain/data_ricker_t0-1.50us_f1MHz/ricker-time.txt}
 \setlength{\plotwidth} {5cm} \setlength{\plotheight}{2.20cm}
 \subfloat[][Ricker function with peak frequency of 1\si{\mega\Hz}.]
         {\makebox[.45\linewidth]{
\begin{tikzpicture}
\begin{axis}[
             enlargelimits=false, 
             ylabel={$S(t)$},
             xlabel={time (\si{\micro\second})},
             enlarge y limits=true,enlarge x limits=false,
             yminorticks=true,
             height=\plotheight,width=\plotwidth,
             label style={font=\scriptsize},
             tick label style={font=\scriptsize},
             legend style={font=\scriptsize\selectfont},
             scale only axis,
             ylabel style = {yshift =-2mm, xshift= 0mm},
             xlabel style = {yshift = 2mm, xshift= 0mm},
             legend pos={north east}, 
             clip mode=individual,
             ]  

     \pgfmathsetmacro{\scaletime}{1000000}
     \pgfmathsetmacro{\scalefreq}{1000}

     \addplot[color=black,line width=1]
              table[x expr = \scaletime*\thisrow{time}, y expr=\thisrow{ricker},
                    ]
              {\datafile}; 
%
%
%
              
\end{axis}
\end{tikzpicture}}
                                  \label{fig:time-domain-acquisition-cx-freq_ricker-t}}
 \renewcommand{\datafile}{figures/tumor_18x18/time-domain/data_ricker_t0-1.50us_f1MHz/ricker-fourier.txt}
 \subfloat[][Complex frequency transform~\cref{eq:complex-dft} of Ricker function.]
         {\makebox[.50\linewidth]{
\begin{tikzpicture}
\begin{axis}[
             enlargelimits=false, 
             ylabel={$\vert \mathfrak{F}(S) \vert$},
             xlabel={ordinary frequency $\omegaR/(2\pi)$ (\si{\kilo\Hz})},
             enlarge y limits=true,enlarge x limits=false,
             ymode=log,
             yminorticks=false,
             xmin =-100,
             xtick={100,1000,2000,3000},
             height=\plotheight,width=\plotwidth,
             label style={font=\scriptsize},
             tick label style={font=\scriptsize},
             legend style={font=\tiny\selectfont},
             scale only axis,
             ylabel style = {yshift =-2mm, xshift=0mm},
             xlabel style = {yshift = 2mm, xshift=0mm},
             legend pos={north east}, 
             clip mode=individual,
             legend style={at={(0.65,0.50)}},
             ]  

     \pgfmathsetmacro{\scalefreq}{0.001}
     \pgfmathsetmacro{\scalerik} {1000000}
     
     \addplot[color=black,line width=1]
              table[x expr = \scalefreq*\thisrow{freq}, y expr=\scalerik*\thisrow{ricker},
                    ]  
              {\datafile}; \addlegendentry{Using $\omegaI=0$}

     \addplot[color=red,line width=1,densely dashed]
              table[x expr = \scalefreq*\thisrow{freq}, y expr=\scalerik*\thisrow{rickeromegaI1e4},
                    ]  
              {\datafile}; \addlegendentry{Using $\omegaI=\num{e4}$}

%
%
%
              
\end{axis}
\end{tikzpicture}}
                                  \label{fig:time-domain-acquisition-cx-freq_ricker-f}}
\caption{Non-invasive measurement setup with Ricker source excitation: 
         the data are obtained from a single point-source excitation or 
         a transducer device, which corresponds to an array of fixed length
         composed of multiple point-sources, which are simultaneously excited.
         The source function is a Ricker-wavelet, which is a broad-band signal.
         The receivers are positioned around the sample (black dashed line).}
\label{fig:time-domain-acquisition}
\end{figure}

In the following of the paper, we only study wave propagation in the 
frequency-domain, or time-harmonic waves, 
omitting the preliminary step that consists in transforming the 
measured time-domain signals.
As illustrated in \cref{fig:time-domain-acquisition-cx-freq_ricker-f}, 
the available frequency content depends on the Ricker source frequency peak.

\graphicspath{{figures/tumor_18x18/time-domain/}}
\setlength{\modelwidth} {8.5cm}
\setlength{\modelheight}{7cm}
\begin{figure}[ht!] \centering
\renewcommand{\modelfile}{trace_scale1_total-time2.5ms}
\subfloat[][Data are obtained in the time-domain (left)
            from a point-source excitation (\cref{fig:time-domain-acquisition-cx-freq_a}), 
            and a complex Fourier transform is applied to produce
            the complex frequency data (right). 
            Here, $y_t$ is the signal at time step $t$ with
            $N_t$ time steps in total, 
            and $x_\mathrm{rcv}$ is the position of a receiver.]{
\begin{tikzpicture}[scale=1]
  \setlength{\modelwidth} {5.25cm} \setlength{\modelheight}{5.75cm}
  \node[rectangle,draw=black,line width=1pt] (sismo)  {\begin{tikzpicture}

  \pgfmathsetmacro{\xmin}{-180}
  \pgfmathsetmacro{\xmax}{180}
  \pgfmathsetmacro{\zmin}{0}
  \pgfmathsetmacro{\zmax}{2.5}
  \pgfmathsetmacro{\cmin}{-1}
  \pgfmathsetmacro{\cmax}{1}

  \pgfmathsetmacro{\xminloc}{-180}
  \pgfmathsetmacro{\xmaxloc}{180}
  \pgfmathsetmacro{\zminloc}{0}  
  \pgfmathsetmacro{\zmaxloc}{0.60}    
\begin{axis}[%
width=\modelwidth, height=\modelheight,
axis on top, separate axis lines,
xmin=\xminloc, xmax=\xmaxloc, xlabel={receivers position angle (\si{\degree})},
ymin=\zminloc, ymax=\zmaxloc, ylabel={time (\si{\milli\second})}, y dir=reverse,
x label style={xshift= 0.0cm, yshift= 0.20cm}, 
y label style={xshift= 0.0cm, yshift=-0.40cm},
colorbar,colormap={blackwhite}{gray(0cm)=(0);gray(1cm)=(1)},
colorbar style={title={{\tiny{$\pressure$ (\si{\pascal})}}},
title style={yshift=-2mm, xshift=-2mm},
ytick={},yticklabels={,,},
width=2mm,xshift=-2mm,
},point meta min=\cmin,point meta max=\cmax,
  label style={font=\tiny},
  tick label style={font=\tiny},
  legend style={font=\tiny\selectfont},
]
\addplot [forget plot] graphics [xmin=\xmin,xmax=\xmax,ymin=\zmin,ymax=\zmax] {{\modelfile}.png};
\end{axis}
\end{tikzpicture}%
};
  \node[rectangle,draw=black,line width=1pt,fill=orange!10!white,
        right = 0.2 of sismo, text width=4cm] (f)  {Fourier transform at frequency 
                          $\omega=\omegaR+\ii\omegaI$};
  \node[below =-0.05 of f, text width=5cm,black,xshift=1em] (f2)  
       {\small $\data(\omega,x_\mathrm{rcv})$\\[0.2em]
        $=\sum_{k=0}^{N_t-1} y_k(x_\mathrm{rcv}) \, e^{\ii\omega \, k/N_t}$};
  \setlength{\plotwidth} {4.90cm}
  \setlength{\plotheight}{2.20cm}
  \pgfmathsetmacro{\ymin}{-120}
  \pgfmathsetmacro{\ymax}{ 120}
  \renewcommand{\datafile}{figures/tumor_18x18/modeling_rcv/sol_center.txt}
  \renewcommand{\dataA}{200kHz_0_real}
  \node[right = 0.10 of f, yshift= 2.9em] (fdata1)  {
\begin{tikzpicture}
\begin{axis}[
             enlargelimits=false, 
             ylabel={$\Real(\data)$},
             xlabel={receivers position angle ($^\circ$)},
             enlarge y limits=false,
             enlarge x limits=false,
             yticklabel pos=right,
             xtick={-90,0,90},
             yminorticks=true,
             ymin=\ymin,
             ymax=\ymax,
             height=\plotheight,width=\plotwidth,
             label style={font=\scriptsize},
             tick label style={font=\scriptsize},
             legend style={font=\scriptsize\selectfont},
             scale only axis,
             ylabel style = {yshift = 7mm, xshift=1mm},
             xlabel style = {yshift = 2mm, xshift=0mm},
             legend pos={north east}, 
             clip mode=individual,
             ]  

     \pgfmathsetmacro{\scale}{1}
     \pgfmathsetmacro{\scalefreq}{1000}

     \addplot[color=\myblue,line width=0.2]
              table[x expr = \thisrow{theta}, y expr=\thisrow{\dataA},
                    ]
              {\datafile}; 
%
%
%
              
\end{axis}
\end{tikzpicture}};
  \renewcommand{\dataA}{200kHz_1e4_real}
  \node[right = 0.10 of f, yshift=-4.0em] (fdata2)  {
\begin{tikzpicture}
\begin{axis}[
             enlargelimits=false, 
             ylabel={$\Real(\data)$},
             xticklabel pos=right,
             yticklabel pos=right,
             enlarge y limits=false,
             enlarge x limits=false,
             xtick={-90,0,90},
             yminorticks=true,
             ymin=\ymin,
             ymax=\ymax,
             height=\plotheight,width=\plotwidth,
             label style={font=\scriptsize},
             tick label style={font=\scriptsize},
             legend style={font=\scriptsize\selectfont},
             scale only axis,
             ylabel style = {yshift = 7mm, xshift=1mm},
             xlabel style = {yshift = 2mm, xshift=0mm},
             legend pos={north east}, 
             clip mode=individual,
             ]  

     \pgfmathsetmacro{\scale}{1}
     \pgfmathsetmacro{\scalefreq}{1000}

     \addplot[color=\myblue,line width=0.2]
              table[x expr = \thisrow{theta}, y expr=\thisrow{\dataA},
                    ]
              {\datafile}; 
%
%
%
              
\end{axis}
\end{tikzpicture}};
  \draw[->,line width=1pt] ([yshift=3em] sismo.east) -| ([xshift=-4.5em] f.north);
  \draw[->,line width=1pt] ([xshift=3em]             f.north)    |- ([yshift=1em] fdata1.west);
  \draw[->,line width=1pt] ([xshift=3em,yshift=-2.50em] f.south)    |- ([yshift=-1.5em] fdata2.west);
  \node[left = of fdata1, anchor=east,yshift= 2.0em,xshift=3.3em]
       {{\color{red!75!black}{\small $\Big( \dfrac{\omega_R}{2\pi}=200\si{\kilo\Hz}, \, 
                                                   \omega_I\hspace*{-0.2em}=\hspace*{-0.2em}0 \Big)$}}};
  \node[left = of fdata2, anchor=east,yshift=-2.5em,xshift=3.3em]
       {{\color{red!75!black}{\small $\Big(\dfrac{\omega_R}{2\pi}=200\si{\kilo\Hz}, \,
                                             \omega_I\hspace*{-0.2em}=\hspace*{-0.2em}
                                             \num{e4}\si{\per\second}\Big)$}}};
  \node[above=0.0 of sismo, anchor=south]{time-domain data};
  \node[above=0.0 of sismo, anchor=south,xshift=27em]{frequency-domain data};

\end{tikzpicture}
\label{fig:time-domain-acquisition-cx-freq_b}
} 

\setlength{\modelwidth} {5.5cm}
\setlength{\modelheight}{5.5cm}
\renewcommand{\modelfile}{time-domain_transducers-3cm_tmax-1ms_scale20}
\subfloat[][Time-domain signals using a transducer 
            source (\cref{fig:time-domain-acquisition-cx-freq_t})
            of length 3\si{\centi\meter}, composed of \num{128} 
            Ricker point-sources simultaneously excited]
            {\makebox[.48\linewidth]{\begin{tikzpicture}

  \pgfmathsetmacro{\xmin}{-180}
  \pgfmathsetmacro{\xmax}{180}
  \pgfmathsetmacro{\zmin}{0}
  \pgfmathsetmacro{\zmax}{1}
  \pgfmathsetmacro{\cmin}{-1}
  \pgfmathsetmacro{\cmax}{1}

  \pgfmathsetmacro{\xminloc}{-180}
  \pgfmathsetmacro{\xmaxloc}{180}
  \pgfmathsetmacro{\zminloc}{0}  
  \pgfmathsetmacro{\zmaxloc}{0.60}    
\begin{axis}[%
width=\modelwidth, height=\modelheight,
axis on top, separate axis lines,
xmin=\xminloc, xmax=\xmaxloc, xlabel={receivers position angle (\si{\degree})},
ymin=\zminloc, ymax=\zmaxloc, ylabel={time (\si{\milli\second})}, y dir=reverse,
x label style={xshift= 0.0cm, yshift= 0.20cm}, 
y label style={xshift= 0.0cm, yshift=-0.40cm},
colorbar,colormap={blackwhite}{gray(0cm)=(0);gray(1cm)=(1)},
colorbar style={title={{\tiny{$\pressure$ (\si{\pascal})}}},
title style={yshift=-2mm, xshift=-2mm},
ytick={},yticklabels={,,},
width=2mm,xshift=-2mm,
},point meta min=\cmin,point meta max=\cmax,
  label style={font=\tiny},
  tick label style={font=\tiny},
  legend style={font=\tiny\selectfont},
]
\addplot [forget plot] graphics [xmin=\xmin,xmax=\xmax,ymin=\zmin,ymax=\zmax] {{\modelfile}.png};
\end{axis}
\end{tikzpicture}%
}} \hspace*{2em}
\renewcommand{\modelfile}{time-domain_transducers-9cm_tmax-1ms_scale20}
\subfloat[][Time-domain signals using a transducer 
            source (\cref{fig:time-domain-acquisition-cx-freq_t})
            of length 9\si{\centi\meter}, composed of \num{128} 
            Ricker-point sources simultaneously excited]
            {\makebox[.48\linewidth]{\begin{tikzpicture}

  \pgfmathsetmacro{\xmin}{-180}
  \pgfmathsetmacro{\xmax}{180}
  \pgfmathsetmacro{\zmin}{0}
  \pgfmathsetmacro{\zmax}{1}
  \pgfmathsetmacro{\cmin}{-1}
  \pgfmathsetmacro{\cmax}{1}

  \pgfmathsetmacro{\xminloc}{-180}
  \pgfmathsetmacro{\xmaxloc}{180}
  \pgfmathsetmacro{\zminloc}{0}  
  \pgfmathsetmacro{\zmaxloc}{0.60}    
\begin{axis}[%
width=\modelwidth, height=\modelheight,
axis on top, separate axis lines,
xmin=\xminloc, xmax=\xmaxloc, xlabel={receivers position angle (\si{\degree})},
ymin=\zminloc, ymax=\zmaxloc, ylabel={time (\si{\milli\second})}, y dir=reverse,
x label style={xshift= 0.0cm, yshift= 0.20cm}, 
y label style={xshift= 0.0cm, yshift=-0.40cm},
colorbar,colormap={blackwhite}{gray(0cm)=(0);gray(1cm)=(1)},
colorbar style={title={{\tiny{$\pressure$ (\si{\pascal})}}},
title style={yshift=-2mm, xshift=-2mm},
ytick={},yticklabels={,,},
width=2mm,xshift=-2mm,
},point meta min=\cmin,point meta max=\cmax,
  label style={font=\tiny},
  tick label style={font=\tiny},
  legend style={font=\tiny\selectfont},
]
\addplot [forget plot] graphics [xmin=\xmin,xmax=\xmax,ymin=\zmin,ymax=\zmax] {{\modelfile}.png};
\end{axis}
\end{tikzpicture}%
}}
\caption{Illustration of the steps from time-domain measurements to 
         the complex frequency data in a two-dimensional setup depending
         on the source excitation. 
         The source is positioned at angle 0\si{\degree}, corresponding
         to the Cartesian coordinates $x=\num{17.5}$ \si{\centi\meter}, 
         $z=9$ \si{\centi\meter} in \cref{fig:time-domain-acquisition}. 
         Inversion experiments following this setup are carried out in 
         \cref{section:fwi-2d}, 3D experiments are in \cref{section:fwi-3d}.
         }
\label{fig:time-domain-acquisition-cx-freq}
\end{figure}
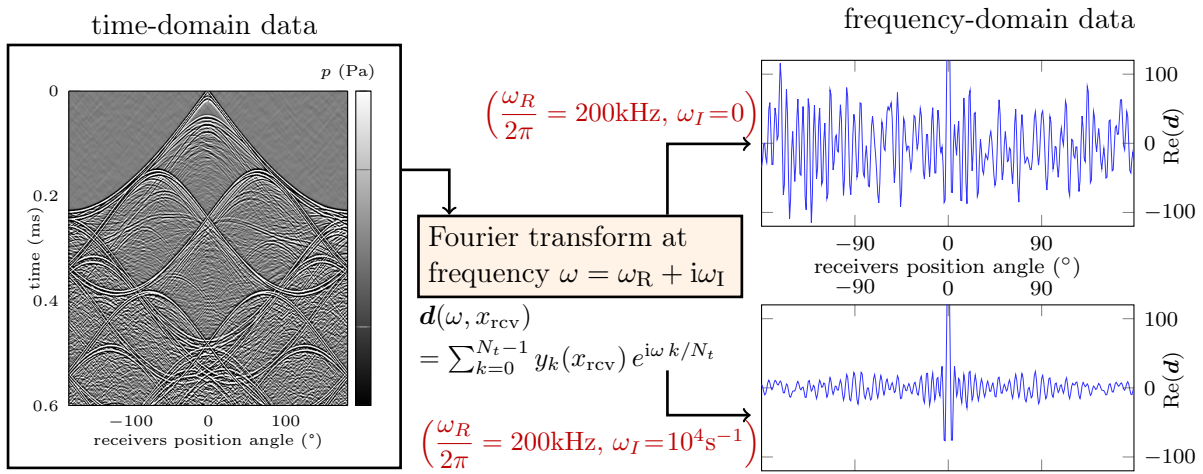
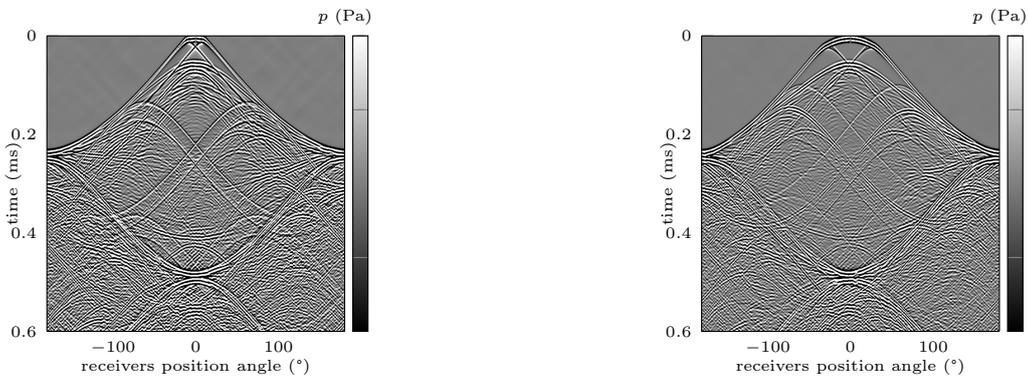


\subsection{Frequency-domain visco-acoustic wave equations}

 
Consequently to the complex Fourier transform of the time-domain data, 
we consider the propagation of waves in the frequency domain.
Let us consider the domain $\Omega$ in dimension $d$ ($2$ 
or $3$ in our experiments), with boundary $\Gamma$. 
The propagation of waves in a visco-acoustic 
medium is described by the particle velocity 
vector field $\velocity:\Omega \rightarrow \Cx^d(\Omega)$ and the scalar 
pressure field $\pressure:\Omega \rightarrow \Cx(\Omega)$, 
that satisfy the first-order system of equations, \cite{Faucher2020adjoint}:
\begin{subequations} \label{eq:euler_main}
\begin{empheq}[left={\empheqlbrace}]{align}
   -\ii \omega \, \rho(\bx) \velocity(\bx, \omega) \,+\, \nabla \pressure(\bx, \omega) & \,=\, 0,  \label{eq:euler_main_a} \\
   -\dfrac{\ii \omega}{\kappa_\dagger(\bx, \omega)} \, \pressure(\bx, \omega) \,+\, \nabla \cdot \velocity(\bx, \omega) 
               & \, = \, g(\bx,\omega), \label{eq:euler_main_b}
\end{empheq} \end{subequations}
where $\bx$ is the space-coordinates and $g$ 
the interior source function. 
The derivation of the complex-frequency time-harmonic 
wave equation is further detailed in \cref{appendix:time-harmonic}.

The material properties describing the medium consist
of the density $\rho:\Omega\rightarrow \R(\Omega)$ and of the bulk 
modulus $\kappa_\dagger:\Omega \rightarrow\Cx(\Omega)$, which is 
complex-valued and frequency dependent to acknowledge 
the effects of attenuation.
The complex wave speed $c_\dagger:\Omega \rightarrow\Cx(\Omega)$ 
is defined from these two coefficients by
\begin{equation} \label{eq:wavespeed:c0}
  c_\dagger(\bx, \omega) \, = \, 
    \sqrt{  \dfrac{\kappa_\dagger(\bx, \omega)}{\rho(\bx)}  } \, ,
\end{equation}
where $\sqrt{\cdot}$ uses the principal 
argument branch $[0, 2\pi )$, 
cf.~\cite[Section~3]{Elbau2017}, that is 
\begin{equation}
  \kappa_\dagger^\beta \,=\, \vert \, \kappa_\dagger\, \vert^\beta 
                       \,\, \big(\beta \,\,\text{arg}(\kappa_\dagger) \big) \,.
\end{equation} 

\paragraph{Free-space propagation}
We refer to free-space propagation when 
waves are assumed to propagate up to infinity, 
without generating reflections. This is typically
assumed in the ultrasound imaging setup, e.g., 
\cite{Pratt2007,Perez2017,Bachmann2020}.
In this case, we impose artificial \emph{absorbing boundary} 
conditions (corresponding to Sommerfeld's radiation condition) 
to constrain the numerical domain in our simulations.
Denoting by $\n$ the normal direction, it corresponds to imposing
the following relation on the numerical boundary $\Gamma$:
\begin{equation}
  \text{absorbing boundary condition:} \qquad  
  -\dfrac{1}{\rho(\bx) \,\, c_\dagger(\bx,\omega)} \, \pressure(\bx,\omega)
  \,+\,  \velocity(\bx,\omega) \cdot \n =0 \,,\,\,\quad \text{on $\Gamma$}. \label{eq:bc:abc}
\end{equation}

\paragraph{Wall boundary conditions}
We also consider a different experiment, where the sample is 
enclosed in a domain with reflecting boundaries. 
Mathematically, this means that Neumann boundary conditions 
for the pressure have to be imposed:
\begin{equation}
  \text{wall boundary condition:}      \qquad\qquad 
      \nabla \pressure(\bx) \cdot \n \,=\,0 \,, 
      \qquad \qquad \qquad \quad \text{ on $\Gamma$}. \label{eq:bc:wall} 
\end{equation} 
This case is more challenging as multiple wave reflections 
occur from the boundary and appear in the data, as we highlight
in the numerical experiments.


\subsection{Attenuation with complex coefficients}

The attenuation of waves is a frequency-dependent
phenomenon, conveniently encoded in the 
frequency-domain equations with complex-valued parameters, 
cf.~\cite{Bland1960,Ursin2002,Carcione2007},
while it can be more challenging to implement
with the time-domain formulation, cf. \cite{Bland1960,Muller1983}.
For propagation with attenuation, the complex wavenumber 
$k$ is decomposed into its real and imaginary parts as follow,
\begin{equation} 
k(\bx,\omega) \,=\, \dfrac{\omega}{c_\dagger(\bx,\omega)} 
              \,=\, \dfrac{\omega_0}{c_p(\bx)} \, + \, \ii \, \alpha(\bx,\omega) \, ,
\end{equation}
where $c_p:\Omega \rightarrow \R_{> 0}(\Omega)$ is the phase velocity, \cite{Ursin2002}, 
the imaginary part $\alpha:\Omega \rightarrow \R_{\geq 0}(\Omega)$ and $\omega_0 \in \R_{>0}$. 

The positiveness of the coefficients leads to conditions 
on the real and imaginary parts of the wave speed and 
frequency. Considering $c_\dagger=\cR + \ii \cI$, we have
\begin{equation}\begin{aligned}
& \dfrac{\omega}{c_\dagger} \,=\, \dfrac{\omegaR \,+\, \ii \omegaI}{\cR \,+\, \ii \cI} 
                            \,=\, \dfrac{\omega_0}{c_p} \, + \, \ii \, \alpha \\
& \Rightarrow \qquad
\dfrac{\omegaR\cR \,-\, \ii \omegaR \cI \,+\, \ii \omegaI \cR \,+\,\omegaI \cI }{\cR^2 \,+\, \cI^2} 
       \,=\, \dfrac{\omega_0}{c_p} \, + \, \ii \, \alpha \\
& \Rightarrow \qquad
\dfrac{\omegaR\cR \,+\,\omegaI \cI }{\cR^2 \,+\, \cI^2} \,=\, \dfrac{\omega_0}{c_p} 
       \qquad \text{and} \qquad 
\dfrac{\omegaI \cR \,-\,\omegaR \cI}{\cR^2 \,+\, \cI^2} 
       \,=\, \alpha \, .
\end{aligned}\end{equation}
Therefore, to ensure that $\omega_0$, $c_p$ and $\alpha$ are positive,
the real and imaginary parts of the wave speed and frequency must verify
\begin{equation} \label{eq:condition-attenuation}
  \omegaR\cR  \,+\,\omegaI \cI \,>\,0 \, , \qquad
  \omegaI \cR \,-\,\omegaR \cI \,\geq\,0 \, .
\end{equation}

Firstly, to verify the first condition in the 
case where $\omegaI=0$ (i.e., for ordinary frequencies), 
we must have that the real part of the wave 
speed is positive: $\cR > 0$.
Now assuming $\cR > 0$, the second condition 
in \cref{eq:condition-attenuation} is always
valid (in particular when $\omegaI=0$) when
\begin{equation}
  \cI \,\leq\,0 \, .
\end{equation}
Therefore we consider that the wave speed has a non positive
imaginary part. This is a minimal condition for attenuation models, see \cite{AkiRichards}.
We summarize the main assumptions regarding attenuation 
below.

\begin{assumption}[Validity of wave propagation with attenuation] \label{assumption:model-validity}
  We assume the following conditions for the 
  formulation of the wave problem \cref{eq:euler_main}:
  \begin{subequations}\label{eq:condition-consistent}
  \begin{align}
    1. \,\, &\text{The complex frequency verifies:}  
    \qquad \qquad \qquad \omegaR > 0 \quad  \text{ and } \quad \omegaI \geq 0 \, ;  \\
    2. \,\, &\text{The real part of the wave speed verifies:}
    \qquad \hspace*{4.4em} \cR > 0 \, ; \\
    3. \,\, &\text{The imaginary part of the wave speed verifies:}
    \qquad \qquad \cI \leq 0 \, ; 
  \end{align}
  \end{subequations}
\end{assumption}
\medskip

In addition, we introduce the \emph{quality factor} $Q$, 
that quantifies the level of attenuation.
\begin{definition}[Quality factor] \label{definition:Q}
The \emph{quality factor} $Q$ 
is defined by the ratio between 
the real and imaginary parts of the complex wave 
speed (\cite{Carcione2007}),
\begin{equation} \label{eq:quality-factor}
  Q(\bx,\omega) \,: = \, \dfrac{\Real\Big( \kappa_\dagger(\bx,\omega)\Big)}{-\Imag \Big(\kappa_\dagger(\bx,\omega)\Big)}
                 \,  \stackrel{\scriptsize \cref{eq:wavespeed:c0}}{=}
                 \, \dfrac{\Real\Big( c_\dagger(\bx,\omega)^2 \Big)}{-\Imag \Big(c_\dagger(\bx,\omega)^2 \Big)} \, > \, 0 \,  .
\end{equation}
The inverse of the quality factor $Q^{-1}$ is referred 
to as the \emph{dissipation factor} and is $0$ in a 
non-attenuating media.
\end{definition}

\subsection{Models of attenuation}

The choice of \emph{attenuation model} indicates the interplay 
between the real and imaginary parts of the wavenumber 
$k$ and its frequency dependency. 
Several models have been
introduced in the literature, cf., e.g.,~\cite{Ursin2002,Carcione2007,Elbau2017}. 
In this work, we are considering seven different models,
which are given in \cref{table:models}. 
In their formulations, we indicate explicitly the 
dependency in frequency, such that all of the variables 
($\kappa_0$, $\eta$, $\upeta$, $\tau$, $\taueps$ and $\tausigma$)
are real and only depend on the space coordinates. 
We further indicate the condition on the parameters to 
ensure the positiveness and negativeness of the signs, 
respectively for the real and imaginary parts of the 
wave speed, see \cref{assumption:model-validity}. 

\begin{table}[H] \begin{center}
\caption{Attenuation models used for the visco-acoustic 
         time-harmonic propagation, extracted from 
         \cite{Ursin2002,Elbau2017} and \cite[Section 2]{Carcione2007}. 
         The parameters $\kappa_0$, $\tau$, $\taueps$, 
         $\tausigma$, $\eta$ and $\upeta$ are real-valued, 
         only depend on the space variable, and are non-negative to 
         validate the conditions given in \cref{eq:condition-consistent}.
         In terms of units, $\kappa_0$ as unit [\si{\pascal}], 
         $\tau$, $\taueps$ and $\tausigma$ are in time units (in [\si{\second}]),
         $\eta$ is in [\si{\pascal\second}] and $\upeta$ is unitless.}
\label{table:models}
\renewcommand{\arraystretch}{1.80}
\begin{tabular}{|>{\arraybackslash}p{.25\linewidth}|
                 >{\arraybackslash}p{.10\linewidth}|
                 >{\arraybackslash}p{.42\linewidth}|
                 >{\arraybackslash}p{.13\linewidth}|}
\hline
\textbf{model} & \textbf{coeff.} & \textbf{complex bulk modulus} & \textbf{condition} \\ \hline
no-attenuation & $\kappa_0$
               & \vspace*{-4em}
                {\begin{flalign} \label{eq:attenuation:no}
               & \kappa_\dagger^{\mathrm{(no)}} \,=\, \kappa_0 &
                \end{flalign}}
                \vspace*{-1.5em}
               & n/a \\ \hline
simplified \newline Kolsky--Futterman & $\kappa_0$, $\upeta$ 
                             & \vspace*{-4em}
                               {\begin{flalign} \label{eq:attenuation:kf}
                               & \kappa_\dagger^{\mathrm{(kf)}} \,=\, 
                                \kappa_0  \, - \, \ii \, \dfrac{\kappa_0}{\upeta} &
                               \end{flalign}}
                               \vspace*{-1em} 
                             & $\upeta \geq 0$ 
                             \\ \hline 
Cole--Cole,
\cite{ColeCole1941,Ursin2002}
       & $\kappa_0$, $\taueps$, $\tausigma$, $\beta$
                 & \vspace*{-4em}
                   {\begin{flalign} \label{eq:attenuation:cc}
                    &\kappa_\dagger^{\mathrm{(cc)}} \,=\, \kappa_0 \, \dfrac{1 \,+\, (- \, \ii \omegaR \, \taueps)^\beta}
                     {1 \,+\, (- \, \ii \omegaR \, \tausigma)^\beta}&
                   \end{flalign}}
                   \vspace*{-1.5em}
                 & $\taueps \geq \tausigma \geq 0$, \newline
                   $0\leq\beta\leq1$ \\ \hline
Zener, \cite{Zener1948}, \newline
       \cite[Section~2.4.3]{Carcione2007}
                 & $\kappa_0$, $\taueps$, $\tausigma$ 
                 & \vspace*{-4em}
                   {\begin{flalign} \label{eq:attenuation:z}
                   & \kappa_\dagger^{\mathrm{(z)}}  \,=\, \kappa_0 \, \dfrac{1 \, - \, \ii \omegaR \,  \taueps}
                  {1 \, - \, \ii \omegaR \, \tausigma} &
                   \end{flalign}}
                   \vspace*{-1.5em}
                 & $\taueps \geq \tausigma \geq 0$ \\ \hline
Kelvin--Voigt, \newline
               \cite[Section~2.4.2]{Carcione2007}
    & $\kappa_0$, $\taueps$ 
                 & \vspace*{-4em}
                   {\begin{flalign} \label{eq:attenuation:kv}
                   & \kappa_\dagger^{\mathrm{(kv)}} \,=\, \kappa_0 \, - \, \ii \omegaR \, \kappa_0 \, \taueps &
                   \end{flalign}}
                   \vspace*{-2em}
                 & $\taueps \geq 0$ \\ \hline
Maxwell, \newline
         \cite[Section~2.4.1]{Carcione2007}
          & $\kappa_0$, $\eta$ 
                 & \vspace*{-4em}
                   {\begin{flalign} \label{eq:attenuation:m}
                   & \kappa_\dagger^{\mathrm{(m)}} \,=\, \dfrac{- \ii \omegaR \, \kappa_0 \, \eta}
                                               {\kappa_0 \, - \, \ii \, \omegaR \, \eta} &
                   \end{flalign}}
                   \vspace*{-1.5em}
                 & $\eta > 0$ \\ \hline
KSB (Kowar--Scherzer--Bonnefond), \cite{Elbau2017,KSB2010}  & $\kappa_0$, $\upeta$, $\tau$, $\beta$ 
                 & \vspace*{-4em}
                   {\begin{flalign} \label{eq:attenuation:ksb}
                   & \kappa_\dagger^{\mathrm{(ksb)}} \hspace*{-0.3em}=\hspace*{-0.3em} \dfrac{\kappa_0}
                                               {\bigg( 1\hspace*{-0.15em}+\hspace*{-0.15em}
                                                        \dfrac{\upeta}{\sqrt{1\hspace*{-0.15em}+\hspace*{-0.15em}
                                                        (-\ii\omegaR\tau)^\beta}}\bigg)^2} &
                   \end{flalign}}
                   \vspace*{-1.5em}
                 & $\upeta\hspace*{-0.1em}>\hspace*{-0.1em}0$, $\tau\hspace*{-0.1em}>\hspace*{-0.1em}0$, \newline 
                   $0<\beta<1$ \\ \hline
modified Szabo, \cite{Szabo1994,Elbau2017}    & $\kappa_0$, $\tau$, $\beta$ 
                 & \vspace*{-4em}
                   {\begin{flalign} \label{eq:attenuation:sz}
                   & \kappa_\dagger^{\mathrm{(sz)}} \,=\, \dfrac{\kappa_0}
                                               {1 \,+\, (-\ii\omegaR\,\tau\,)^{\beta-1}}                                               
                                               &
                   \end{flalign}}
                   \vspace*{-1.5em}
                 & $\tau > 0$, \newline
                   $0<\beta<1$ \\ \hline
\end{tabular} \end{center}
\end{table}

We note that in the case where $\beta=1$, 
the Cole--Cole model is equivalent to 
the Zener model.
The Kolsky--Futterman model \cref{eq:attenuation:kf} 
originates from \cite{Kolsky1956,Futterman1962}
and in our work, we use a simplified version 
which is frequency independent. It is obtained 
under the assumptions of \emph{weak attenuation} 
($4Q^2 \gg 1$), as detailed in \cite[Section~1.7]{Faucher2017}) 
This model (in its frequency-dependent version) 
is also used in \cite{Ribodetti1998,Malinowski2011}.
In the definition of the modified Szabo model {\cref{eq:attenuation:sz}}, 
        $\tau$ is placed together with the frequency and both are at power
        $(\beta-1)$, contrary to {\cite{Elbau2017}} where $\tau$ 
        is separated from the power of $\omega$; This is only motivated 
        for unit consistency, such that $\tau$ is expressed in time 
        unit in {\cref{eq:attenuation:sz}}.

\begin{remark}[Causality]
 The models summarized in \cref{table:models} are causal 
 as they are directly derived from their time-domain 
 counterpart, cf.~\cite{Golden1988,Ursin2002,Carcione2007,Elbau2017}. 
 More generally, causality principles can be verified
 via the Kramers--Kronig relations, or ensuring that
 the parameters are analytic in the lower-half complex 
 plane, and we refer to \cite{Golden1988,Carcione2007,Elbau2017}
 for more details.
\end{remark}

\begin{remark}[Generalized models] \label{rk:generalized}
  The attenuation models given in \cref{table:models} accounts 
  for one frequency mechanism of attenuation, which is sufficient when 
  the data contain a relatively narrow band of frequency, 
  as highlighted in our experiments.
  For broadband signals, one has to consider multiple mechanisms 
  each associated with a different frequency, introducing a 
  \emph{generalized} model, cf. \cite{Carcione1988,Nachman1990,Carcione2007,Mclaughlin2011,BrossierVirieux2016}.
  For instance, considering $L$ attenuation mechanisms each 
  associated with frequency $\omega_{l}$, one can define, 
  e.g., \cite{BrossierVirieux2016},
  \begin{equation}
       \kappa_\dagger^{\mathrm{(G)}}(\bx,\omega)
       \,=\, \kappa_0(\bx) \, \bigg( 1 \, - \, \sum_{l=1}^L \, B_l(\bx,\omega) 
          \,\, \dfrac{\omega_{l}}{\omega_{l} \,+\,  \ii \omega}\,\bigg) \, ,
  \end{equation}
  where $B_l$ represents the attenuation function of the $l^{\text{th}}$ mechanism.
\end{remark}

\subsection{Comparison of attenuation models}
\label{subsection:numerical_model_comparisons}

We illustrate the dependence of the quality factor $Q$ 
defined in \cref{eq:quality-factor} with the frequency
for the different attenuation models in \cref{fig:plot:comparison-attenuation}.
We fix the values of the density to $\rho=$ \num{1000} \si{\kg\per\meter\cubed},
and of the bulk modulus $\kappa_0=$ \num{2.25} \si{\giga\Pa}, 
such that the wave speed is $c_0 = \sqrt{\kappa_0/\rho} = $\num{1500} \si{\meter\per\second}.
To select and compare the parameters of each of the 
attenuation models, we impose the quality factor at 
the reference frequency $\omegaR/(2\pi)=\num{300} \si{\kilo\Hz}$ 
to be $Q=118$, leading to the values represented in \cref{table:models:Q118}.

\begin{table}[H] 
\begin{center}
\caption{Choice of coefficients for the attenuation models
         to compute $Q$ in \cref{fig:plot:comparison-attenuation}. 
         These are chosen such that the 
         quality factor at \num{300} \si{\kHz} is $Q=$\num{118}.
         For models that involve multiple coefficients,
         there exist different combinations that would 
         also give $Q=118$ at \num{300} \si{\kHz} and we 
         only choose one of them.}
\label{table:models:Q118}
\renewcommand{\arraystretch}{1.10}
\begin{tabular}{ >{\arraybackslash}p{.30\linewidth}|
                 >{\arraybackslash}p{.50\linewidth}}
\hline
\textbf{attenuation model} & \textbf{parameters} \\ \hline
simplified Kolsky--Futterman & $\upeta=\num{118}$ \\ \hline 
Kelvin--Voigt    &  $\taueps=\num{4.5}$ \si{\nano\second} \\ \hline
Maxwell          &  $\eta=\num{1.4e5}$ \si{\pascal\second} \\ \hline
Zener            &  $\taueps=\num{90}$ \si{\nano\second}, $\tausigma=\num{85.4}$ \si{\nano\second} \\ \hline
Cole--Cole       &  $\taueps=\num{90.5}$ \si{\nano\second}, $\tausigma=\num{85.5}$ \si{\nano\second}, $\beta=0.8$ \\ \hline
KSB              &  $\tau=\num{2e5}$ \si{\second}, $\upeta=\num{8.75}$, $\beta=\num{0.5}$ \\ \hline
modified Szabo   &  $\tau=\num{13.28}$ \si{\second}, $\beta=\num{0.6}$ \\ \hline
\end{tabular} 
\end{center}
\end{table}

\setlength{\plotheight}{5.25cm}
\setlength{\plotwidth} {11.0cm}
\renewcommand{\datafile}{figures/attenuation-models/data/attenuation_models.txt}
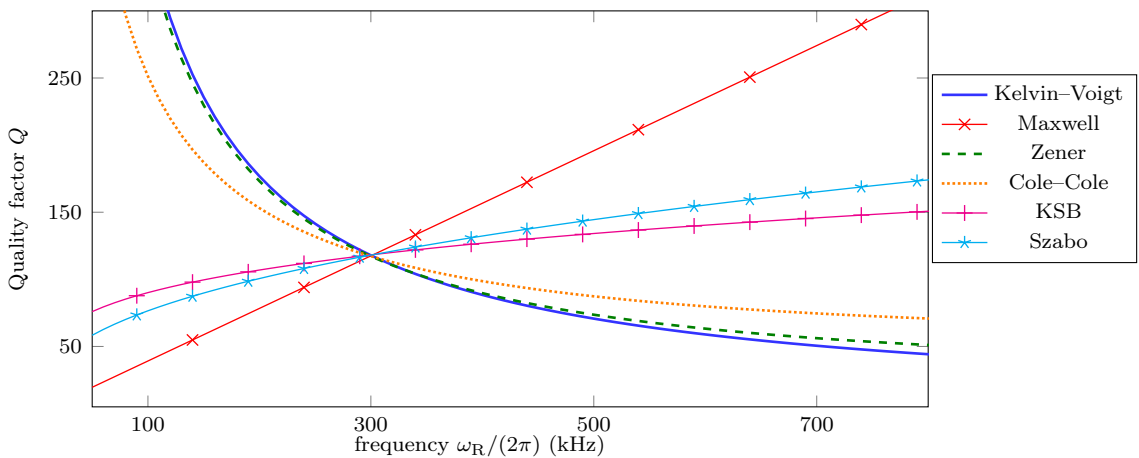
\begin{figure}[ht!] \centering
  \renewcommand{\dataA}{KV}       \renewcommand{\legendA}{Kelvin--Voigt}
  \renewcommand{\dataB}{Maxwell}  \renewcommand{\legendB}{Maxwell}
  \renewcommand{\dataC}{ZenerA}   \renewcommand{\legendC}{Zener}
  \renewcommand{\dataD}{ZenerB}   \renewcommand{\legendD}{Zener ($\tausigma=\num{75}\si{\nano\second}$)}
  \renewcommand{\dataE}{ColeA}    \renewcommand{\legendE}{Cole--Cole}
  \renewcommand{\dataF}{ColeB}    \renewcommand{\legendF}{Cole--Cole ($\beta=\num{0.70}$)}
  \renewcommand{\dataG}{KSB}      \renewcommand{\legendG}{KSB}
  \renewcommand{\dataH}{Szabo}    \renewcommand{\legendH}{Szabo}
\begin{tikzpicture}
\begin{axis}[
             enlargelimits=false, 
             ylabel={Quality factor $Q$},
             xlabel={frequency $\omegaR/(2\pi)$ (\si{\kilo\Hz})},
             enlarge y limits=false,
             enlarge x limits=false,
             yminorticks=true,
             ymin=5,
             ymax=300,
             xmin=50,
             xmax=800,
             height=\plotheight,width=\plotwidth,
             label style={font=\scriptsize},
             tick label style={font=\scriptsize},
             legend style={font=\scriptsize\selectfont},
             scale only axis,
             ylabel style = {yshift =0mm, xshift=0mm},
             xlabel style = {yshift =2mm, xshift=-4mm},
             legend pos={north east}, 
             clip mode=individual,
             legend columns=1,
             legend style={at={(1.25,0.84)}},
             xtick={100,300,500,700},
             ytick={50,150,250,350}
             ]  

     \pgfmathsetmacro{\scale}{1}
     \pgfmathsetmacro{\scalefreq}{1000}


     \addplot[color=\myblue,line width=1]
              table[x expr = \thisrow{freq}/\scalefreq, y expr=\scale*\thisrow{\dataA},
                    ]
              {\datafile}; \addlegendentry{\legendA}    
     \addplot[color=red,line width=0.5,mark=x, mark repeat=10, mark phase=10, mark size=3]
              table[x expr = \thisrow{freq}/\scalefreq, y expr=\thisrow{\dataB}]
              {\datafile}; \addlegendentry{\legendB}    

     \addplot[color=\mygreen,line width=1,dashed]
              table[x expr = \thisrow{freq}/\scalefreq, y expr=\thisrow{\dataC}]
              {\datafile}; \addlegendentry{\legendC} 

     \addplot[color=orange,line width=1,densely dotted]
              table[x expr = \thisrow{freq}/\scalefreq, y expr=\thisrow{\dataE}]
              {\datafile}; \addlegendentry{\legendE} 

     \addplot[color=magenta,line width=0.5,mark=+, mark repeat=5, mark phase=5, mark size=3]
              table[x expr = \thisrow{freq}/\scalefreq, y expr=\thisrow{\dataG}]
              {\datafile}; \addlegendentry{\legendG} 
     \addplot[color=cyan   ,line width=0.5,mark=star, mark repeat=5, mark phase=5, mark size=2.50]
              table[x expr = \thisrow{freq}/\scalefreq, y expr=\thisrow{\dataH}]
              {\datafile}; \addlegendentry{\legendH} 
                 
\end{axis}
\end{tikzpicture}
  \vspace*{-0.5em}
  \caption{Comparisons of attenuation models with frequency
           using $c_0=\num{1500}$ \si{\meter\per\second} and 
           $\rho=\num{1000}$ \si{\kg\per\meter\cubed}. 
           The attenuation parameters for each model are chosen 
           such that their quality factor $Q$ is equal to \num{118} 
           at frequency 300 \si{\kilo\Hz}, cf. \cref{table:models:Q118}.
           }
\label{fig:plot:comparison-attenuation}
\end{figure}

One main difference between the attenuation models 
is that either the low or high frequencies are 
attenuated the most.
For instance, the quality factor decreases with 
frequency for the Kelvin--Voigt, Zener and Cole--Cole 
models, while it increases with
frequency for the Maxwell, KSB and modified
Szabo ones.
For the latter family, it means that low-frequency 
waves are more attenuated while it is the high-frequency 
waves that are attenuated the most with the first
family of models.
In the investigated frequency band, 
we observe that the Kelvin--Voigt, Zener and Cole--Cole 
models have the same pattern with high quality factors at 
low frequencies, a rapid decrease, and then a stabilization.
The KSB and modified Szabo models instead have low
quality factors at low frequencies, before it 
slowly increases with higher frequencies. 
One could also use the \emph{generalized models}, as 
discussed in \cref{rk:generalized}, to allow more 
flexibility on the quality factor's evolution.
Here, we see that, 
while the quality factors of the attenuation models all coincide 
at frequency \num{300} \si{\kilo\Hz}, they become all different
when stepping away from this reference frequency.

\section{Quantitative inverse problem using iterative minimization}
\label{section:inverse-problem}

The inverse wave problem associated to the forward 
equation \cref{eq:euler_main} 
corresponds to the reconstruction of the heterogeneous
properties $\kappa_\dagger(\bx)$ and $\rho(\bx)$ 
given some wave measurements.
In this section, we review the steps for quantitative 
imaging based upon iterative minimization. 
In the context of seismic imaging where one uses the 
phase and amplitude information recorded by the 
seismograms 
this is usually 
referred to as the \emph{Full Waveform Inversion} (FWI), 
\cite{Lailly1983,Tarantola1984,Virieux2009,Faucher2017}.

\subsection{Forward problem}

We assume measurements of the pressure field 
$\pressure$ at the position of receivers. 
As described in \cref{subsection:time-to-frequency}, we 
consider that the time-domain measurements results from 
a source, which is mathematically represented as a Ricker wavelet. 
Then, a complex Fourier 
(or Laplace--Fourier) transform is applied to work with 
frequency-domain data, see \cref{fig:time-domain-acquisition-cx-freq}.
Therefore, we define the \emph{forward problem} $\forward$ 
directly in the frequency domain: For a source $g$ (which corresponds 
to the complex Fourier transform of the original Ricker source), it is 
defined such that,
\begin{equation}
  \forward(\bm, \omega, g) \,=\, \Big\lbrace \, 
                           \pressure(\bx_1, \omega, g), \,
                           \pressure(\bx_2, \omega, g), \, 
                           \ldots , \,
                           \pressure(\bx_{\nrcv}, \omega, g)
                           \Big\rbrace \, ,
\end{equation}
where $\pressure$ solves \cref{eq:euler_main} 
with right-hand side $g$.
Here, $\bx_1, \ldots, \bx_{\nrcv}$ denotes the 
position of the $\nrcv$ receivers.
The model parameters are represented by $\bm$ and correspond 
to the bulk modulus and density such that
$\bm:=\{ \kappa_\dagger, \, \rho \}$.
We have $\nsrc$ independent sources in the 
experiment, and define the data-set associated 
to the frequency $\omega$ by
\begin{equation}
  \forward(\bm, \omega) \,=\, \Big\lbrace \, 
                        \forward(\bm, \omega, g_1), \,    
                        \forward(\bm, \omega, g_2), \,
                        \ldots, \,
                        \forward(\bm, \omega, g_{\nsrc})
                        \Big\rbrace \, .
\end{equation}
Here, we keep the notation $\forward$ when there is no ambiguity.
Note that we have assumed that the receivers remains in 
the same position for each source (that is, they do not depend 
on $g$), but modifying the receiver positions together with the 
source would not modify the methodology.

\subsection{Reconstruction with iterative minimization, FWI}

The reconstruction is carried out following the minimization
of a misfit functional $\misfit$, that evaluates a distance 
between measurements $\data$ and simulations:
\begin{equation}
  \misfit(\bm, \omega) \,  = \, \mathrm{dist}\Big( \forward(\bm,\omega)  , \, \data(\omega) \Big) \, ,
\end{equation}
where the choice of misfit function, depending on the data-sets and
inverted parameters, is the subject of numerous studies, e.g., 
\cite{Shin2008,Fichtner2008,Brossier2010,VanLeeuwen2010,Karaouglu2017,Faucher2019FRgWIGeo,Faucher2020EV,Faucher2020DAS}.
Per simplicity, we rely on the 
$l^2$-distance of the difference, such that
\begin{equation}
  \misfit(\bm, \omega) \,=\, 
    \dfrac{1}{2} \,\sum_{k=1}^{\nsrc}\, \big\Vert \, \forward(\bm,\omega,g_k) \,-\, \data(\omega,g_k) 
                   \big\Vert^2_2 \, .
\end{equation}

The reconstruction is performed by minimizing $\misfit$ 
with respect to the model parameters $\bm$, with a truncated-Newton
method.
We further follow a progressive increase in the frequency-content 
to improve the convergence of the algorithm and mitigate ill-posedness, 
as advocated in \cite{Bunks1995,Faucher2019basins,Faucher2020basins}.

The pressure wave propagating through the medium is broadband 
(see \cref{fig:time-domain-acquisition-cx-freq_ricker-f}). It 
is recorded at a finite number of receivers.  
Then the complex Fourier transform of the time signals is performed, and we 
get information for all frequencies
$\omega=\omegaR + \ii\omegaI$ 
such that
$\omegaR^{(\min)} \leq \omegaR \leq \omegaR^{(\max)}$ and 
$\omegaI^{(\min)} \leq \omegaI \leq \omegaI^{(\max)}$.
Following the guidelines of \cite{Faucher2019basins,Faucher2020basins}, 
we further use a \emph{sequential} frequency progression (i.e., one $\omega$ 
at a time) in the reconstruction algorithm (instead of a band of 
frequencies inverted at once). 
The progression in the selection of $\omegaR$ and $\omegaI$ is 
detailed and motivated below.
The iterative minimization algorithm is summarized in \cref{algo:FWI}.

\paragraph{Ordinary frequency progression}
When the imaginary part of the frequency is zero, $\omegaI=0$, 
(see \cref{fig:time-domain-acquisition-cx-freq}), one works with 
ordinary frequencies only. In this case the progression in frequency 
content is chosen from low to high values, 
\cite{Bunks1995,Faucher2019basins,Faucher2020basins,Faucher2017,Virieux2009}.
That is, we start with $\omega^{(1)} \,=\, \omegaR^{(1)} \geq \omegaR^{(min)}$,
carry out minimization iterations, and then update to 
$\omega^{(2)} \,=\, \omegaR^{(2)} > \omegaR^{(1)}$, and 
repeat with the next frequency, see \cref{algo:FWI}.

\paragraph{Complex frequency progression}
In the case where the imaginary part of the frequency is 
non-zero (illustrated in \cref{fig:time-domain-acquisition-cx-freq}), 
the progression in the frequency content follows the strategy of 
\cite{Faucher2020basins}:
$\omegaI$ varies first, from high to low values with $\omegaR$ fixed.
Then $\omegaR$ is increased and we repeat the sequence of $\omegaI$
(from high to low). This choice is motivated from the estimates of 
attraction basins size provided in \cite{Faucher2020basins}. 
Namely, we start with $\omega^{(1)} = \omegaR^{(1)} + \ii \omegaI^{(1)}$,
such that $\omegaR^{(1)}\geq\omegaR^{(\min)}$ and $\omegaI^{(1)}\leq\omegaI^{(\max)}$.
The next frequency is $\omega^{(2)} = \omegaR^{(1)} + \ii \omegaI^{(2)}$,
with $\omegaI^{(2)}\leq\omegaI^{(1)}$, see \cref{algo:FWI}.

\begin{algorithm}[ht!]
 \caption{Iterative minimization procedure for quantitative reconstruction 
          using complex frequencies. The real and imaginary parts of the 
          frequency are ordered such that
          $\omegaI^{(1)} > \omegaI^{(2)} > \ldots > \omegaI^{(N_{\omegaI})}$
          and
          $\omegaR^{(1)} < \omegaR^{(2)} < \ldots < \omegaR^{(N_{\omegaR})}$.
          }
 \label{algo:FWI}
 \SetAlgoLined
 Initialization: starting model parameters $\bm_1 = (\kappa_{\dagger,1}, \, \rho_1)$,
                 and measurements $\data$. \\

 \For{$i = 1, \, \ldots, \, {(N_{\omegaI})}$}{
   \For{$j = 1, \, \ldots, \, {(N_{\omegaR})}$}{
     Set $\omega := \omegaR^{(j)} \,+\, \ii \omegaI^{(i)}$ \;
     Compute the complex Fourier transform of 
          the time-domain measurement at $\omega$, 
          see \cref{eq:complex-dft,fig:time-domain-acquisition-cx-freq} \;

     \For{$k = 1, \, \ldots, \, n_\text{iter}$}{
     Set $l:= (i-1)\,N_{\omegaR}\,n_\text{iter} \,+\,  (j-1)n_\text{iter} + k$ \;
     Solve Problem~\cref{eq:euler_main} using current models $\bm_l$ and frequency $\omega$
           for all sources in the acquisition\;
     Compute the misfit function $\misfit(\bm_l, \omega)$ \;
     Compute the gradient of the misfit function with adjoint-state 
             method, \cite{Faucher2020adjoint}\;
     Compute the search direction (e.g., nonlinear conjugate gradient, \cite{Nocedal2006}) \;
     Compute the step length, $\alpha_l$ using line-search method, \cite{Nocedal2006} \;
     Update the model: $\bm_{l+1}  = \bm_l - \alpha_l s_l$.
   }
  }
 }
\end{algorithm}

\subsection{Numerical implementation}
\label{subsection:hdg}

For the discretization of the forward wave 
problem \cref{eq:euler_main}, we use the 
Hybridizable Discontinuous Galerkin method (HDG), 
\cite{Cockburn2009,Faucher2020adjoint}.
The HDG method used the mixed formulation
and 
static condensation
to solve Problem~\cref{eq:euler_main}.
It results in a linear system to be solved, 
which size is shown to be less than, e.g., for finite 
elements method, depending on the degree of the 
approximation polynomials,  see.~\cite{Kirby2012,Kirby2016,Faucher2020adjoint}.
Our implementation for the forward and inverse
problems follows the steps described in 
\cite{Faucher2020adjoint}, where the only change
is to consider a complex-valued bulk modulus $\kappa_\dagger$.
The computation of the gradient is performed with
the adjoint-state method (\cite{Pratt1998,Plessix2006,Faucher2019basins})
specifically derived for the HDG discretization in \cite{Faucher2020adjoint}.
For the two and three dimensional computational 
experiments that are carried out in the following 
sections, we use the open-source software 
\texttt{hawen}\footnote{\url{https://ffaucher.gitlab.io/hawen-website/} \label{footnote:haven}}, 
\cite{Hawen2021}.
The precise computational cost is further discussed in \cref{subsection:computational-cost}

%

\section{Numerical experiments in two dimensions}
\label{section:fwi-2d}
\pgfmathsetmacro{\xminloc}{2}
\pgfmathsetmacro{\xmaxloc}{16}
\pgfmathsetmacro{\zminloc}{4}
\pgfmathsetmacro{\zmaxloc}{14.5}

In this section, we investigate a two-dimensional experiment, 
following the setup described in \cref{fig:time-domain-acquisition-cx-freq}.
The reconstruction is carried out following \cref{algo:FWI},
and we study the following:
\begin{enumerate}\setlength{\itemsep}{-1pt}
  \item We evaluate the robustness of the reconstruction 
        procedure with attenuation model uncertainty, that is, 
        we select an attenuation model to generate 
        synthetic measurement data and use a different model to carry out the 
        reconstruction. 
  \item While the primary focus of this paper is 
        on the reconstruction of the wave speed variations, 
        we also investigate the choice of parametrization 
        to recover the medium's density, which is known to 
        be harder to recover, see \cite{Virieux2009,Jeong2012,Faucher2017} 
        and the references therein. 
  \item We investigate the effect of the domain boundary conditions, that is, 
        whether we allow waves to freely escape the domain using 
        absorbing boundary conditions (free-space) or assuming wall 
        boundaries around the sample.
\end{enumerate}
The following experiments have all been realized using 
software \texttt{hawen}, \cite{Hawen2021}, see \cref{footnote:haven}. 
We further detail the computational cost of the 2D and 3D 
reconstructions in \cref{subsection:computational-cost}.

\subsection{Experimental setup}

We consider a synthetic experiment of ultrasound imaging 
with a two-dimensional breast model, described by its wave 
speed $c_0$, density $\rho$, and quality factor $Q$, pictured 
in \cref{fig:models:true}.
The medium is of size \num{18}$\times$\num{18} \si{\cm\squared}
and is composed of tissues, blood and fat. These layers are
obtained from a cross-section of the OA-Breast Phantom 
data-set\footnote{\url{https://anastasiolab.wustl.edu/downloadable-content/oa-breast-database/}.
\label{footnote:OA-breast}},
\cite{Lou2017}, where the values of the parameters depending
on the layers (e.g., blood, skin) are taken following 
the IT'IS Database\footnote{\url{https://itis.swiss/virtual-population/tissue-properties/}.\label{footnote:itis}},
see \cref{table:interval}. 
In addition, we incorporate a high-contrast ellipsoid inclusion 
in the model parameters, which can be seen as a defect that has 
to be identified with ultrasound imaging. 
From the values given in \cref{table:interval}, we see that
the quality factor is the lowest in the skin layer, that is, 
waves are attenuated the most by skin.  In the other layers, 
the quality factor is relatively high (at least 280),
such that we have a weakly attenuating medium, \cite{Vavryvcuk2009}.

\begin{table}[H] \begin{center}
\caption{Interval of the model values depending on the type of tissues.
         These follow the IT'IS database, see \cref{footnote:itis}.}
\label{table:interval}
\renewcommand{\arraystretch}{1.20}
\begin{tabular}{|>{\arraybackslash}p{.17\linewidth}|
                 >{\arraybackslash}p{.26\linewidth}|
                 >{\arraybackslash}p{.22\linewidth}|
                 >{\arraybackslash}p{.20\linewidth}|}
\hline
\textbf{medium} & \textbf{wave speed} $c_0$  (\si{\meter\per\second}) 
                & \textbf{density}    $\rho$ (\si{\kg\per\meter\cubed})
                & \textbf{$Q$} at \num{300} \si{\kilo\Hz} \\ \hline
background      & \num{1490} & 1000 & 800 \\ \hline
skin            & $(\num{1590}, \, \num{1610})$ & $(\num{1100}, \, \num{1120})$ & $(\num{100}, \, \num{120})$ \\ \hline
blood           & $(\num{1565}, \, \num{1575})$ & $(\num{1090}, \, \num{1110})$ & $(\num{290}, \, \num{310})$ \\ \hline
fat             & $(\num{1440}, \, \num{1460})$ & $(\num{ 920}, \, \num{ 940})$ & $(\num{410}, \, \num{430})$ \\ \hline
grandular tissue& $(\num{1490}, \, \num{1520})$ & $(\num{1030}, \, \num{1050})$ & $(\num{280}, \, \num{300})$ \\ \hline
inclusion & \num{1550}                    & \num{1050}                    & \num{350}                   \\ \hline
\end{tabular} \end{center}
\end{table}

\setlength{\modelwidth} {4.70cm}
\setlength{\modelheight}{4.70cm}
\graphicspath{{figures/tumor_18x18/main/}}
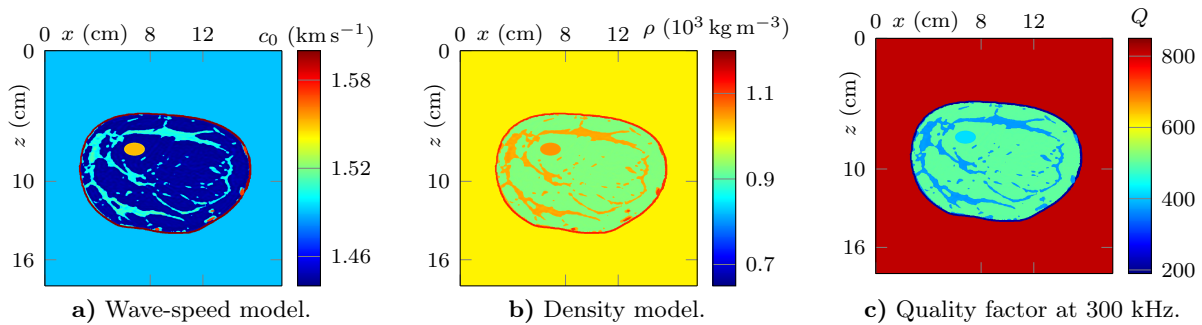
\begin{figure}[ht!] \centering
\renewcommand{\modelfile}{cp-true_scale1440-1600}
\subfloat[][Wave-speed model.]{\begin{tikzpicture}

  \pgfmathsetmacro{\xmin}{0}
  \pgfmathsetmacro{\xmax}{18}
  \pgfmathsetmacro{\zmin}{0}
  \pgfmathsetmacro{\zmax}{18}
  \pgfmathsetmacro{\cmin}{1.44}
  \pgfmathsetmacro{\cmax}{1.60}
  
\begin{axis}[%
width=\modelwidth, height=\modelheight,
axis on top, separate axis lines,
xmin=\xmin, xmax=\xmax, xlabel={$x$ (\si{\cm})},
ymin=\zmin, ymax=\zmax, ylabel={$z$ (\si{\cm})}, y dir=reverse,
xtick={0,8,12},
ytick={0,10,16},
xticklabel pos=right, xlabel near ticks,
x label style={xshift=-0.9cm, yshift=-0.50cm}, 
y label style={xshift= 0.7cm, yshift=-0.50cm},
colormap/jet,colorbar,
colorbar style={title={{\scriptsize{$c_0$ (\si{\km\per\second})}}},
title style={yshift=-3mm, xshift=1mm},
width=3mm,xshift=-1mm,
ytick={1.46,1.52,1.58},
},point meta min=\cmin,point meta max=\cmax,
  label style={font=\scriptsize},
  tick label style={font=\scriptsize},
  legend style={font=\scriptsize\selectfont},
]
\addplot [forget plot] graphics [xmin=\xmin,xmax=\xmax,ymin=\zmin,ymax=\zmax] {{\modelfile}.png};
\end{axis}
\end{tikzpicture}%
} \hfill
\renewcommand{\modelfile}{rho-true_scale650-1200}
\subfloat[][Density model.]{\begin{tikzpicture}

  \pgfmathsetmacro{\xmin}{0}
  \pgfmathsetmacro{\xmax}{18}
  \pgfmathsetmacro{\zmin}{0}
  \pgfmathsetmacro{\zmax}{18}
  \pgfmathsetmacro{\cmin}{0.65}
  \pgfmathsetmacro{\cmax}{1.20}
  
\begin{axis}[%
width=\modelwidth, height=\modelheight,
axis on top, separate axis lines,
xmin=\xmin, xmax=\xmax, xlabel={$x$ (\si{\cm})},
ymin=\zmin, ymax=\zmax, ylabel={$z$ (\si{\cm})}, y dir=reverse,
xtick={0,8,12},
ytick={0,10,16},
xticklabel pos=right, xlabel near ticks,
x label style={xshift=-0.9cm, yshift=-0.50cm}, 
y label style={xshift= 0.7cm, yshift=-0.50cm},
colormap/jet,colorbar,
colorbar style={title={{\scriptsize{$\rho$ (\num{e3}$\,$\si{\kg\per\meter\cubed})}}},
title style={yshift=-2mm, xshift=-1mm},
width=3mm,xshift=-1mm,
ytick={0.7,0.9,1.10},
},point meta min=\cmin,point meta max=\cmax,
  label style={font=\scriptsize},
  tick label style={font=\scriptsize},
  legend style={font=\scriptsize\selectfont},
]
\addplot [forget plot] graphics [xmin=\xmin,xmax=\xmax,
                                 ymin=\zmin,ymax=\zmax] {{\modelfile}.png};
\end{axis}
\end{tikzpicture}%
}   \hfill 
\renewcommand{\modelfile}{Q-true_scale190-850}
\subfloat[][Quality factor at \num{300} \si{\kilo\Hz}.]
           {\begin{tikzpicture}

  \pgfmathsetmacro{\xmin}{0}
  \pgfmathsetmacro{\xmax}{18}
  \pgfmathsetmacro{\zmin}{0}
  \pgfmathsetmacro{\zmax}{18}
  \pgfmathsetmacro{\cmin}{190}
  \pgfmathsetmacro{\cmax}{850}
  
\begin{axis}[%
width=\modelwidth, height=\modelheight,
axis on top, separate axis lines,
xmin=\xmin, xmax=\xmax, xlabel={$x$ (\si{\cm})},
ymin=\zmin, ymax=\zmax, ylabel={$z$ (\si{\cm})}, y dir=reverse,
xtick={0,8,12},
ytick={0,10,16},
xticklabel pos=right, xlabel near ticks,
x label style={xshift=-0.9cm, yshift=-0.50cm}, 
y label style={xshift= 0.7cm, yshift=-0.50cm},
colormap/jet,colorbar,
colorbar style={title={{\scriptsize{$Q$}}},
title style={yshift=-2mm, xshift=-0mm},
width=3mm,xshift=-1mm,
},point meta min=\cmin,point meta max=\cmax,
  label style={font=\scriptsize},
  tick label style={font=\scriptsize},
  legend style={font=\scriptsize\selectfont},
]
\addplot [forget plot] graphics [xmin=\xmin,xmax=\xmax,
                                 ymin=\zmin,ymax=\zmax] {{\modelfile}.png};
\end{axis}
\end{tikzpicture}%
            \label{fig:models:true_Q}}
\caption{Two-dimensional breast models of wave-speed, density 
         and quality factor at frequency $\omegaR/(2\pi)=\num{300}$ \si{\kilo\Hz} 
         including a high-contrast inclusion. 
         The values of the parameters depending
         on the type of tissues are detailed 
         in \cref{table:interval}. The domain is 
         of size \num{18}$\times$\num{18} \si{\cm\squared}.}
\label{fig:models:true}
\end{figure}

To carry out the reconstruction, we assume measurements 
of the pressure field, following the setup described in 
\cref{subsection:time-to-frequency}:
A Ricker source wavelet excites the medium and time-domain 
data are measured. The Ricker wavelet conveys multiple 
frequency contents (\cref{fig:time-domain-acquisition-cx-freq_ricker-f}) 
and one transforms the recorded time-domain data into 
complex frequency-domain ones, see \cref{fig:time-domain-acquisition-cx-freq} 
Nonetheless, as we consider a synthetic experiment, 
for simplicity we directly generate the data-sets in
the frequency domain and incorporate white noise a-posteriori.
The acquisition is composed of \num{36} independent 
point-sources positioned in a circular pattern, 
as illustrated in \cref{fig:time-domain-acquisition-cx-freq_a}.
The receivers measuring the pressure field for each of the source 
are positioned onto the same disk, with a total of \num{360} receivers.

For the reconstruction procedure, the initial models 
consist in the constant values of the 
parameters in the background, that is,  
the inversion starts with
$c_0 \,=\,\num{1500}$ \si{\meter\per\second}, 
$\rho\,=\,\num{1000}$ \si{\kilo\gram\per\meter\cubed} 
and attenuation parameters chosen (depending on the 
attenuation model) such that 
$Q=800$ at $\omegaR/(2\pi) = \num{300}$ \si{\kilo\Hz}.

\subsection{Comparison of wave propagation with attenuation models and boundary conditions}
\label{subsection:2d-data}

We generate synthetic data for each of the 
seven attenuation models given in \cref{table:models},
leading to seven different data-sets. 
The choice of the parameters that describe the attenuation
is such that the quality factors $Q$ of each model 
are equal at \si{300} \si{\kilo\Hz}, and can be seen in 
\cref{fig:models:true_Q}.
This amounts to the intervals of values given in \cref{table:models:2d}.

\begin{table}[H] 
\begin{center}
\caption{Interval of values for the coefficients of 
         the attenuation models to represent the 
         parameters given \cref{fig:models:true,table:interval}.}
\label{table:models:2d}
\renewcommand{\arraystretch}{1.10}
\begin{tabular}{ >{\arraybackslash}p{.30\linewidth}|
                 >{\arraybackslash}p{.50\linewidth}}
\hline
\textbf{model} & \textbf{parameters} \\ \hline
simplified Kolsky--Futterman & $\upeta\in(\num{100},\num{800})$ \\ \hline 
Kelvin--Voigt    &  $\taueps\in(\num{0.663},\num{5.305})$ \si{\nano\second} \\ \hline
Maxwell          &  $\eta\in(\num{0.148},\num{0.942})$ \si{\mega\pascal\second} \\ \hline
Zener            &  $\taueps\in(\num{4.729},\num{8.663})$ \si{\nano\second}, 
                    $\tausigma\in(\num{1},\num{8})$ \si{\nano\second} \\ \hline
Cole--Cole       &  $\taueps\in(\num{3.150},\num{8.387})$ \si{\nano\second}, 
                    $\tausigma\in(\num{1},\num{8})$ \si{\nano\second}, $\beta=0.8$. \\ \hline
KSB              &  $\tau\in(\num{1},\num{8})$ \si{\second}, 
                    $\upeta\in(\num{17.332},\num{3.036e8})$, $\beta=\num{0.5}$.\\ \hline
modified Szabo   &  $\tau\in(\num{5.905},\num{384.552})$ \si{\second}, $\beta=\num{0.5}$\\ \hline
\end{tabular} 
\end{center}
\end{table}

\begin{figure}[ht!]\centering
\subfloat[Signals at \num{100} \si{\kilo\Hz} using wall boundary conditions on the sides,
          with one data point per \num{0.1}\si{\degree}.]
           {\makebox[.95\linewidth][c]{\includegraphics[]{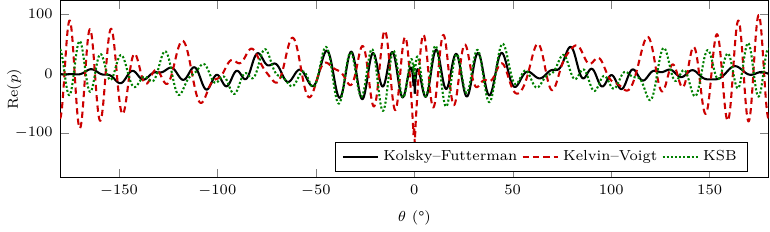}}
            \label{fig:data_a}} \\

  \subfloat[Signals at \num{500} \si{\kilo\Hz} using wall boundary conditions on the sides, 
            with one data point per \num{0.1}\si{\degree}.]
           {\makebox[.48\linewidth][c]{\includegraphics[]{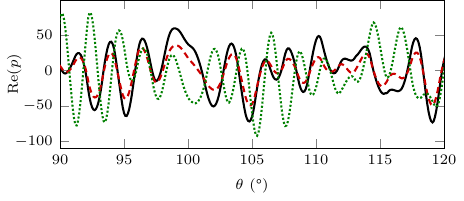}}
            \label{fig:data_b}} \hfill
  \subfloat[Signals at \num{500} \si{\kilo\Hz} using absorbing boundary conditions on the sides, 
            with one data point per \num{0.1}\si{\degree}.]
           {\makebox[.48\linewidth][c]{\includegraphics[]{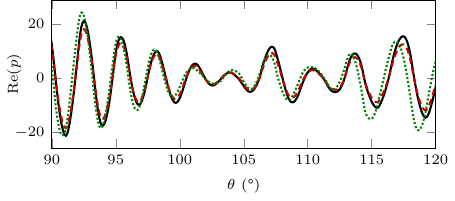}}
            \label{fig:data_c}} \\

  \subfloat[Signals at \num{300} \si{\kilo\Hz} using wall boundary conditions on the sides, 
            with one data point per \num{0.1}\si{\degree}.]
           {\makebox[.48\linewidth][c]{\includegraphics[]{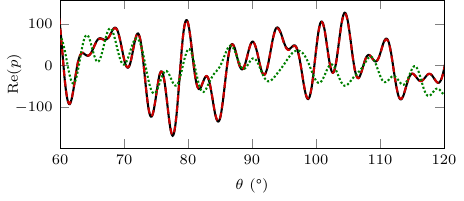}}
            \label{fig:data_d}} \hfill
  \subfloat[Signals at \num{300} \si{\kilo\Hz} using wall boundary conditions on the sides, 
            with one data point per \num{1}\si{\degree}.]
           {\makebox[.48\linewidth][c]{\includegraphics[]{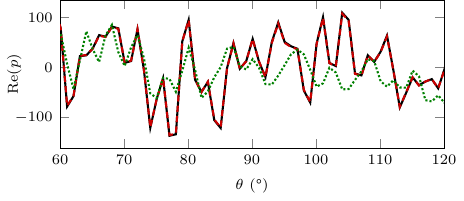}}
            \label{fig:data_e}}

\caption{Comparison of the pressure field at the position
         of the receivers depending on 
         the attenuation model of the medium,
         with medium parameters pictured in
         \cref{fig:models:true}. 
         The receiver devices are positioned around  
         the sample, in the configuration of 
         \cref{fig:time-domain-acquisition-cx-freq},
         and the source is at 
         angle $\theta=$\num{0}\si{\degree}. 
         For visualization, we first consider 
         one receiver every \num{0.1}\si{\degree}.
         Nevertheless, the data used for inversion
         consider one receiver every \num{1}\si{\degree}
         instead, i.e., \num{360} receivers.
         }
\label{fig:data}
\end{figure}

For the models that have several degrees of freedom, 
we consider different magnitudes 
in the parameters. 
We illustrate the data depending on the medium attenuation 
model in \cref{fig:data}, for a source positioned at angle
0\si{\degree} (see \cref{fig:time-domain-acquisition-cx-freq})
and the medium properties pictured in \cref{fig:models:true}.
We further compare the simulations with absorbing or wall 
boundary conditions. For the sake of clarity,  we only 
compare three ordinary (real) frequencies, that is, we 
keep $\omegaI=0$,
	and picture the results for 
the Kolsky--Futterman, Kelvin--Voigt, and KSB models. 
The data used for inversion are obtained every 1\si{\degree} 
(i.e., \num{360} receivers, see \cref{fig:data_e}), but we also 
plot in \cref{fig:data} the signals with data-point 
every \num{0.1}\si{\degree} for better visualization 
and comparisons.

\graphicspath{{figures/tumor_18x18/modeling_compare-models/}}

We observe that the three selected attenuation models 
all give different signals. In \cref{fig:data_a}, we 
see that the difference exists for all angles $\theta$, 
and particularly at large offset, with drastic changes in 
terms of amplitude.
In \cref{fig:data_b,fig:data_c}, we picture the 
solution at \num{500}\si{\kHz} for $\theta\in(90\si{\degree},120\si{\degree})$
depending on the choice of the boundary condition:
We observe that when using absorbing conditions 
(i.e., assuming free-space propagation, \cref{fig:data_c} 
), 
the difference between the signals is small both in terms 
of amplitude and phase. On the other hand, when assuming 
the wall boundaries on the sides of the domain, we see 
drastic changes in the wavefields, cf. \cref{fig:data_c}.
In \cref{fig:data_d}, we plot the results for 
frequency 300 \si{\kilo\Hz}, at which the quality factors 
of each attenuation model is the same, here we see that 
the signals corresponding to the Kolsky--Futterman and 
Kelvin--Voigt models are indeed very close. Nonetheless,
the signal corresponding to the KSB model is different 
in amplitude and phase. 
Furthermore, to carry out the inversion, we see that we have a 
sub-sampled signal with one data-point every 1\si{\degree}, 
as plotted in \cref{fig:data_e}.

Note also that there is some flexibility with the position of the 
absorbing boundary conditions to constrain the computational 
domain, as waves are anyway supposed to freely propagate 
up until infinity.
On the other hand, when considering wall boundaries, it is
essential that the computational domain boundaries follow
the actual setup, so that the reflections from the wall are 
adequately simulated.

\subsection{Reconstruction with absorbing boundary conditions (free-space)}

We consider the medium in free-space, such that
acoustic waves propagate according to \cref{eq:euler_main} 
in the domain $\Omega$ with absorbing boundary conditions
\cref{eq:bc:abc} on $\Gamma$.
In this configuration, one assumes that either the 
propagation domain is sufficiently large so that 
reflections from the boundaries are negligible, or
that these reflections can somehow be removed 
in a pre-processing stage of the data. 
Note that the later option is \emph{not} a trivial 
task in general (cf. the suppression of the so-called 
multiples in seismic acquisition) and can lead to 
inconsistency in the data.
We further investigate the case of reflecting boundaries, 
which is more challenging, in \cref{subsection:2d:wall-bc}.

\subsubsection{Wave speed reconstruction with attenuation model uncertainty}

In this test-case, we investigate the robustness of the
iterative reconstruction procedure with respect to the 
attenuation model.
We generate a synthetic data-set for each of the 
attenuation model (see \cref{subsection:2d-data}) 
and carry out the reconstruction using a different 
one.
All of the data-sets have the same acquisition with
\num{36} independent point-sources and \num{360}
receivers measuring the pressure field, as illustrated
in \cref{fig:time-domain-acquisition-cx-freq,fig:data}.
We use synthetic data, but the numerical setup 
is changed between the forward and inverse problems.
Namely, the data are generated using a discretization 
mesh of more than \num{200000} cells to ensure the 
accuracy of the different layers, and the solution is 
approximated with polynomials of order 5. 
On the other hand, the inversion is carried out on a mesh
of less than \num{10000} cells, with variable orders 
using $p$-adaptivity for efficiency, cf.~\cite{Faucher2020adjoint}.
In addition, white noise is added to the synthetic data, 
using a signal-to-noise ratio of \num{20}\si{\deci\bel}.
During inversion, the model parameter is represented as a 
piecewise-polynomial function per cell, using Lagrange 
basis function of linear order. 
Therefore, the model parameters
are represented by \num{3} unknowns on each cell of the discretized mesh.

In this experiment we only use ordinary (real) frequencies 
(i.e., regular Fourier transform of the time-domain data) such that 
$\omegaI=0$. 
We consider the increasing frequency progression as follows:
$\omegaR/(2\pi) = \{\num{200}, \, \num{300}, \, \num{400}, \, \num{600} \}$ \si{\kilo\Hz}.
The evolution in the frequency content of the data
is chosen sequentially, \cite{Faucher2020basins}, 
and \num{30} iterations are performed per frequency for 
the minimization, following \cref{algo:FWI}.
Therefore, there is a total of \num{120} iterations.
The reconstruction is carried out with respect to the 
bulk modulus $\kappa_0$, that is, the density $\rho$ 
and quality factor $Q$ remain at their constant 
initial values. 
For visualization, we instead picture the resulting wave
speed $c_0=\sqrt{\kappa_0/\rho}$, 
assembled from the reconstructed real part of the 
bulk modulus $\kappa_0$ and the (homogeneous) initial 
density $\rho$. 
In \cref{fig:fwi2d:abosrbing:all} are shown the reconstructed 
wave speed models,
comparing for all data-sets the use of a different model for
inversion where, for visualization, the raw results are slightly 
smoothed a-posteriori, see~\cref{rk:visualization}.

\setlength{\modelwidth} {3.80cm}
\setlength{\modelheight}{3.40cm}
\graphicspath{{figures/tumor_18x18/fwi_absorbing_noise20dB/}}
\newlength{\myhspace}
\setlength{\myhspace}{-0.25em}
\newlength{\myvspace}
\setlength{\myvspace}{-0.10em}
\begin{figure}[ht!]\centering

\setlength{\myhspace}{-0.60em}
\subfloat[][Reconstruction using the data-set 
            generated with the simplified Kolsky--Futterman model.]{
\renewcommand{\modelfile}{cp-600kHz-data-kf-1}\begin{tikzpicture}

  \pgfmathsetmacro{\xmin}{0}
  \pgfmathsetmacro{\xmax}{18}
  \pgfmathsetmacro{\zmin}{0}
  \pgfmathsetmacro{\zmax}{18}
  \pgfmathsetmacro{\cmin}{1.47}
  \pgfmathsetmacro{\cmax}{1.55}
  
\begin{axis}[%
width=\modelwidth, height=\modelheight,
axis on top, separate axis lines,
xmin=\xminloc, xmax=\xmaxloc, xlabel={$x$ (\si{\cm})},
ymin=\zminloc, ymax=\zmaxloc, ylabel={$z$ (\si{\cm})}, y dir=reverse,
xtick={4,10,12},
ytick={4,10,14},
xticklabel pos=right, xlabel near ticks,
x label style={xshift=-0.5cm, yshift=-0.50cm}, 
y label style={xshift= 0.5cm, yshift=-0.50cm},
hide axis,
  label style={font=\scriptsize},
  tick label style={font=\scriptsize},
  legend style={font=\scriptsize\selectfont},
]
\addplot [forget plot] graphics [xmin=\xmin,xmax=\xmax,ymin=\zmin,ymax=\zmax] {{\modelfile}.png};
\end{axis}
\end{tikzpicture}%
 \hspace*{\myhspace}
\renewcommand{\modelfile}{cp-600kHz-data-kf-5}\begin{tikzpicture}

  \pgfmathsetmacro{\xmin}{0}
  \pgfmathsetmacro{\xmax}{18}
  \pgfmathsetmacro{\zmin}{0}
  \pgfmathsetmacro{\zmax}{18}
  \pgfmathsetmacro{\cmin}{1.47}
  \pgfmathsetmacro{\cmax}{1.55}
  
\begin{axis}[%
width=\modelwidth, height=\modelheight,
axis on top, separate axis lines,
xmin=\xminloc, xmax=\xmaxloc, xlabel={$x$ (\si{\cm})},
ymin=\zminloc, ymax=\zmaxloc, ylabel={$z$ (\si{\cm})}, y dir=reverse,
xtick={4,10,12},
ytick={4,10,14},
xticklabel pos=right, xlabel near ticks,
x label style={xshift=-0.5cm, yshift=-0.50cm}, 
y label style={xshift= 0.5cm, yshift=-0.50cm},
hide axis,
  label style={font=\scriptsize},
  tick label style={font=\scriptsize},
  legend style={font=\scriptsize\selectfont},
]
\addplot [forget plot] graphics [xmin=\xmin,xmax=\xmax,ymin=\zmin,ymax=\zmax] {{\modelfile}.png};
\end{axis}
\end{tikzpicture}%
 \hspace*{\myhspace}
\renewcommand{\modelfile}{cp-600kHz-data-kf-7}\begin{tikzpicture}

  \pgfmathsetmacro{\xmin}{0}
  \pgfmathsetmacro{\xmax}{18}
  \pgfmathsetmacro{\zmin}{0}
  \pgfmathsetmacro{\zmax}{18}
  \pgfmathsetmacro{\cmin}{1.47}
  \pgfmathsetmacro{\cmax}{1.55}
  
\begin{axis}[%
width=\modelwidth, height=\modelheight,
axis on top, separate axis lines,
xmin=\xminloc, xmax=\xmaxloc, xlabel={$x$ (\si{\cm})},
ymin=\zminloc, ymax=\zmaxloc, ylabel={$z$ (\si{\cm})}, y dir=reverse,
xtick={4,10,12},
ytick={4,10,14},
xticklabel pos=right, xlabel near ticks,
x label style={xshift=-0.5cm, yshift=-0.50cm}, 
y label style={xshift= 0.5cm, yshift=-0.50cm},
hide axis,
  label style={font=\scriptsize},
  tick label style={font=\scriptsize},
  legend style={font=\scriptsize\selectfont},
]
\addplot [forget plot] graphics [xmin=\xmin,xmax=\xmax,ymin=\zmin,ymax=\zmax] {{\modelfile}.png};
\end{axis}
\end{tikzpicture}%
 \hspace*{\myhspace}
\renewcommand{\modelfile}{cp-600kHz-data-kf-2}\begin{tikzpicture}

  \pgfmathsetmacro{\xmin}{0}
  \pgfmathsetmacro{\xmax}{18}
  \pgfmathsetmacro{\zmin}{0}
  \pgfmathsetmacro{\zmax}{18}
  \pgfmathsetmacro{\cmin}{1.47}
  \pgfmathsetmacro{\cmax}{1.55}
  
\begin{axis}[%
width=\modelwidth, height=\modelheight,
axis on top, separate axis lines,
xmin=\xminloc, xmax=\xmaxloc, xlabel={$x$ (\si{\cm})},
ymin=\zminloc, ymax=\zmaxloc, ylabel={$z$ (\si{\cm})}, y dir=reverse,
xtick={4,10,12},
ytick={4,10,14},
xticklabel pos=right, xlabel near ticks,
x label style={xshift=-0.5cm, yshift=-0.50cm}, 
y label style={xshift= 0.5cm, yshift=-0.50cm},
hide axis,
  label style={font=\scriptsize},
  tick label style={font=\scriptsize},
  legend style={font=\scriptsize\selectfont},
]
\addplot [forget plot] graphics [xmin=\xmin,xmax=\xmax,ymin=\zmin,ymax=\zmax] {{\modelfile}.png};
\end{axis}
\end{tikzpicture}%
 \hspace*{\myhspace}
\renewcommand{\modelfile}{cp-600kHz-data-kf-6}\begin{tikzpicture}

  \pgfmathsetmacro{\xmin}{0}
  \pgfmathsetmacro{\xmax}{18}
  \pgfmathsetmacro{\zmin}{0}
  \pgfmathsetmacro{\zmax}{18}
  \pgfmathsetmacro{\cmin}{1.47}
  \pgfmathsetmacro{\cmax}{1.55}
  
\begin{axis}[%
width=\modelwidth, height=\modelheight,
axis on top, separate axis lines,
xmin=\xminloc, xmax=\xmaxloc, xlabel={$x$ (\si{\cm})},
ymin=\zminloc, ymax=\zmaxloc, ylabel={$z$ (\si{\cm})}, y dir=reverse,
xtick={4,10,12},
ytick={4,10,14},
xticklabel pos=right, xlabel near ticks,
x label style={xshift=-0.5cm, yshift=-0.50cm}, 
y label style={xshift= 0.5cm, yshift=-0.50cm},
hide axis,
  label style={font=\scriptsize},
  tick label style={font=\scriptsize},
  legend style={font=\scriptsize\selectfont},
]
\addplot [forget plot] graphics [xmin=\xmin,xmax=\xmax,ymin=\zmin,ymax=\zmax] {{\modelfile}.png};
\end{axis}
\end{tikzpicture}%
 \hspace*{\myhspace}
\renewcommand{\modelfile}{cp-600kHz-data-kf-3}\begin{tikzpicture}

  \pgfmathsetmacro{\xmin}{0}
  \pgfmathsetmacro{\xmax}{18}
  \pgfmathsetmacro{\zmin}{0}
  \pgfmathsetmacro{\zmax}{18}
  \pgfmathsetmacro{\cmin}{1.47}
  \pgfmathsetmacro{\cmax}{1.55}
  
\begin{axis}[%
width=\modelwidth, height=\modelheight,
axis on top, separate axis lines,
xmin=\xminloc, xmax=\xmaxloc, xlabel={$x$ (\si{\cm})},
ymin=\zminloc, ymax=\zmaxloc, ylabel={$z$ (\si{\cm})}, y dir=reverse,
xtick={4,10,12},
ytick={4,10,14},
xticklabel pos=right, xlabel near ticks,
x label style={xshift=-0.5cm, yshift=-0.50cm}, 
y label style={xshift= 0.5cm, yshift=-0.50cm},
hide axis,
  label style={font=\scriptsize},
  tick label style={font=\scriptsize},
  legend style={font=\scriptsize\selectfont},
]
\addplot [forget plot] graphics [xmin=\xmin,xmax=\xmax,ymin=\zmin,ymax=\zmax] {{\modelfile}.png};
\end{axis}
\end{tikzpicture}%
 \hspace*{\myhspace}
\renewcommand{\modelfile}{cp-600kHz-data-kf-4}\begin{tikzpicture}

  \pgfmathsetmacro{\xmin}{0}
  \pgfmathsetmacro{\xmax}{18}
  \pgfmathsetmacro{\zmin}{0}
  \pgfmathsetmacro{\zmax}{18}
  \pgfmathsetmacro{\cmin}{1.47}
  \pgfmathsetmacro{\cmax}{1.55}
  
\begin{axis}[%
width=\modelwidth, height=\modelheight,
axis on top, separate axis lines,
xmin=\xminloc, xmax=\xmaxloc, xlabel={$x$ (\si{\cm})},
ymin=\zminloc, ymax=\zmaxloc, ylabel={$z$ (\si{\cm})}, y dir=reverse,
xtick={4,10,12},
ytick={4,10,14},
xticklabel pos=right, xlabel near ticks,
x label style={xshift=-0.5cm, yshift=-0.50cm}, 
y label style={xshift= 0.5cm, yshift=-0.50cm},
hide axis,
  label style={font=\scriptsize},
  tick label style={font=\scriptsize},
  legend style={font=\scriptsize\selectfont},
]
\addplot [forget plot] graphics [xmin=\xmin,xmax=\xmax,ymin=\zmin,ymax=\zmax] {{\modelfile}.png};
\end{axis}
\end{tikzpicture}%
} \\[\myvspace]

\subfloat[][Reconstruction using the data-set 
            generated with the simplified Kelvin--Voigt model.]{
\renewcommand{\modelfile}{cp-600kHz-data-kv-1}\begin{tikzpicture}

  \pgfmathsetmacro{\xmin}{0}
  \pgfmathsetmacro{\xmax}{18}
  \pgfmathsetmacro{\zmin}{0}
  \pgfmathsetmacro{\zmax}{18}
  \pgfmathsetmacro{\cmin}{1.47}
  \pgfmathsetmacro{\cmax}{1.55}
  
\begin{axis}[%
width=\modelwidth, height=\modelheight,
axis on top, separate axis lines,
xmin=\xminloc, xmax=\xmaxloc, xlabel={$x$ (\si{\cm})},
ymin=\zminloc, ymax=\zmaxloc, ylabel={$z$ (\si{\cm})}, y dir=reverse,
xtick={4,10,12},
ytick={4,10,14},
xticklabel pos=right, xlabel near ticks,
x label style={xshift=-0.5cm, yshift=-0.50cm}, 
y label style={xshift= 0.5cm, yshift=-0.50cm},
hide axis,
  label style={font=\scriptsize},
  tick label style={font=\scriptsize},
  legend style={font=\scriptsize\selectfont},
]
\addplot [forget plot] graphics [xmin=\xmin,xmax=\xmax,ymin=\zmin,ymax=\zmax] {{\modelfile}.png};
\end{axis}
\end{tikzpicture}%
 \hspace*{\myhspace}
\renewcommand{\modelfile}{cp-600kHz-data-kv-5}\begin{tikzpicture}

  \pgfmathsetmacro{\xmin}{0}
  \pgfmathsetmacro{\xmax}{18}
  \pgfmathsetmacro{\zmin}{0}
  \pgfmathsetmacro{\zmax}{18}
  \pgfmathsetmacro{\cmin}{1.47}
  \pgfmathsetmacro{\cmax}{1.55}
  
\begin{axis}[%
width=\modelwidth, height=\modelheight,
axis on top, separate axis lines,
xmin=\xminloc, xmax=\xmaxloc, xlabel={$x$ (\si{\cm})},
ymin=\zminloc, ymax=\zmaxloc, ylabel={$z$ (\si{\cm})}, y dir=reverse,
xtick={4,10,12},
ytick={4,10,14},
xticklabel pos=right, xlabel near ticks,
x label style={xshift=-0.5cm, yshift=-0.50cm}, 
y label style={xshift= 0.5cm, yshift=-0.50cm},
hide axis,
  label style={font=\scriptsize},
  tick label style={font=\scriptsize},
  legend style={font=\scriptsize\selectfont},
]
\addplot [forget plot] graphics [xmin=\xmin,xmax=\xmax,ymin=\zmin,ymax=\zmax] {{\modelfile}.png};
\end{axis}
\end{tikzpicture}%
 \hspace*{\myhspace}
\renewcommand{\modelfile}{cp-600kHz-data-kv-7}\begin{tikzpicture}

  \pgfmathsetmacro{\xmin}{0}
  \pgfmathsetmacro{\xmax}{18}
  \pgfmathsetmacro{\zmin}{0}
  \pgfmathsetmacro{\zmax}{18}
  \pgfmathsetmacro{\cmin}{1.47}
  \pgfmathsetmacro{\cmax}{1.55}
  
\begin{axis}[%
width=\modelwidth, height=\modelheight,
axis on top, separate axis lines,
xmin=\xminloc, xmax=\xmaxloc, xlabel={$x$ (\si{\cm})},
ymin=\zminloc, ymax=\zmaxloc, ylabel={$z$ (\si{\cm})}, y dir=reverse,
xtick={4,10,12},
ytick={4,10,14},
xticklabel pos=right, xlabel near ticks,
x label style={xshift=-0.5cm, yshift=-0.50cm}, 
y label style={xshift= 0.5cm, yshift=-0.50cm},
hide axis,
  label style={font=\scriptsize},
  tick label style={font=\scriptsize},
  legend style={font=\scriptsize\selectfont},
]
\addplot [forget plot] graphics [xmin=\xmin,xmax=\xmax,ymin=\zmin,ymax=\zmax] {{\modelfile}.png};
\end{axis}
\end{tikzpicture}%
 \hspace*{\myhspace}
\renewcommand{\modelfile}{cp-600kHz-data-kv-2}\begin{tikzpicture}

  \pgfmathsetmacro{\xmin}{0}
  \pgfmathsetmacro{\xmax}{18}
  \pgfmathsetmacro{\zmin}{0}
  \pgfmathsetmacro{\zmax}{18}
  \pgfmathsetmacro{\cmin}{1.47}
  \pgfmathsetmacro{\cmax}{1.55}
  
\begin{axis}[%
width=\modelwidth, height=\modelheight,
axis on top, separate axis lines,
xmin=\xminloc, xmax=\xmaxloc, xlabel={$x$ (\si{\cm})},
ymin=\zminloc, ymax=\zmaxloc, ylabel={$z$ (\si{\cm})}, y dir=reverse,
xtick={4,10,12},
ytick={4,10,14},
xticklabel pos=right, xlabel near ticks,
x label style={xshift=-0.5cm, yshift=-0.50cm}, 
y label style={xshift= 0.5cm, yshift=-0.50cm},
hide axis,
  label style={font=\scriptsize},
  tick label style={font=\scriptsize},
  legend style={font=\scriptsize\selectfont},
]
\addplot [forget plot] graphics [xmin=\xmin,xmax=\xmax,ymin=\zmin,ymax=\zmax] {{\modelfile}.png};
\end{axis}
\end{tikzpicture}%
 \hspace*{\myhspace}
\renewcommand{\modelfile}{cp-600kHz-data-kv-6}\begin{tikzpicture}

  \pgfmathsetmacro{\xmin}{0}
  \pgfmathsetmacro{\xmax}{18}
  \pgfmathsetmacro{\zmin}{0}
  \pgfmathsetmacro{\zmax}{18}
  \pgfmathsetmacro{\cmin}{1.47}
  \pgfmathsetmacro{\cmax}{1.55}
  
\begin{axis}[%
width=\modelwidth, height=\modelheight,
axis on top, separate axis lines,
xmin=\xminloc, xmax=\xmaxloc, xlabel={$x$ (\si{\cm})},
ymin=\zminloc, ymax=\zmaxloc, ylabel={$z$ (\si{\cm})}, y dir=reverse,
xtick={4,10,12},
ytick={4,10,14},
xticklabel pos=right, xlabel near ticks,
x label style={xshift=-0.5cm, yshift=-0.50cm}, 
y label style={xshift= 0.5cm, yshift=-0.50cm},
hide axis,
  label style={font=\scriptsize},
  tick label style={font=\scriptsize},
  legend style={font=\scriptsize\selectfont},
]
\addplot [forget plot] graphics [xmin=\xmin,xmax=\xmax,ymin=\zmin,ymax=\zmax] {{\modelfile}.png};
\end{axis}
\end{tikzpicture}%
 \hspace*{\myhspace}
\renewcommand{\modelfile}{cp-600kHz-data-kv-3}\begin{tikzpicture}

  \pgfmathsetmacro{\xmin}{0}
  \pgfmathsetmacro{\xmax}{18}
  \pgfmathsetmacro{\zmin}{0}
  \pgfmathsetmacro{\zmax}{18}
  \pgfmathsetmacro{\cmin}{1.47}
  \pgfmathsetmacro{\cmax}{1.55}
  
\begin{axis}[%
width=\modelwidth, height=\modelheight,
axis on top, separate axis lines,
xmin=\xminloc, xmax=\xmaxloc, xlabel={$x$ (\si{\cm})},
ymin=\zminloc, ymax=\zmaxloc, ylabel={$z$ (\si{\cm})}, y dir=reverse,
xtick={4,10,12},
ytick={4,10,14},
xticklabel pos=right, xlabel near ticks,
x label style={xshift=-0.5cm, yshift=-0.50cm}, 
y label style={xshift= 0.5cm, yshift=-0.50cm},
hide axis,
  label style={font=\scriptsize},
  tick label style={font=\scriptsize},
  legend style={font=\scriptsize\selectfont},
]
\addplot [forget plot] graphics [xmin=\xmin,xmax=\xmax,ymin=\zmin,ymax=\zmax] {{\modelfile}.png};
\end{axis}
\end{tikzpicture}%
 \hspace*{\myhspace}
\renewcommand{\modelfile}{cp-600kHz-data-kv-4}\begin{tikzpicture}

  \pgfmathsetmacro{\xmin}{0}
  \pgfmathsetmacro{\xmax}{18}
  \pgfmathsetmacro{\zmin}{0}
  \pgfmathsetmacro{\zmax}{18}
  \pgfmathsetmacro{\cmin}{1.47}
  \pgfmathsetmacro{\cmax}{1.55}
  
\begin{axis}[%
width=\modelwidth, height=\modelheight,
axis on top, separate axis lines,
xmin=\xminloc, xmax=\xmaxloc, xlabel={$x$ (\si{\cm})},
ymin=\zminloc, ymax=\zmaxloc, ylabel={$z$ (\si{\cm})}, y dir=reverse,
xtick={4,10,12},
ytick={4,10,14},
xticklabel pos=right, xlabel near ticks,
x label style={xshift=-0.5cm, yshift=-0.50cm}, 
y label style={xshift= 0.5cm, yshift=-0.50cm},
hide axis,
  label style={font=\scriptsize},
  tick label style={font=\scriptsize},
  legend style={font=\scriptsize\selectfont},
]
\addplot [forget plot] graphics [xmin=\xmin,xmax=\xmax,ymin=\zmin,ymax=\zmax] {{\modelfile}.png};
\end{axis}
\end{tikzpicture}%
}\\[\myvspace]

\subfloat[][Reconstruction using the data-set 
            generated with the simplified Maxwell model.]{
\renewcommand{\modelfile}{cp-600kHz-data-ma-1}\begin{tikzpicture}

  \pgfmathsetmacro{\xmin}{0}
  \pgfmathsetmacro{\xmax}{18}
  \pgfmathsetmacro{\zmin}{0}
  \pgfmathsetmacro{\zmax}{18}
  \pgfmathsetmacro{\cmin}{1.47}
  \pgfmathsetmacro{\cmax}{1.55}
  
\begin{axis}[%
width=\modelwidth, height=\modelheight,
axis on top, separate axis lines,
xmin=\xminloc, xmax=\xmaxloc, xlabel={$x$ (\si{\cm})},
ymin=\zminloc, ymax=\zmaxloc, ylabel={$z$ (\si{\cm})}, y dir=reverse,
xtick={4,10,12},
ytick={4,10,14},
xticklabel pos=right, xlabel near ticks,
x label style={xshift=-0.5cm, yshift=-0.50cm}, 
y label style={xshift= 0.5cm, yshift=-0.50cm},
hide axis,
  label style={font=\scriptsize},
  tick label style={font=\scriptsize},
  legend style={font=\scriptsize\selectfont},
]
\addplot [forget plot] graphics [xmin=\xmin,xmax=\xmax,ymin=\zmin,ymax=\zmax] {{\modelfile}.png};
\end{axis}
\end{tikzpicture}%
 \hspace*{\myhspace}
\renewcommand{\modelfile}{cp-600kHz-data-ma-5}\begin{tikzpicture}

  \pgfmathsetmacro{\xmin}{0}
  \pgfmathsetmacro{\xmax}{18}
  \pgfmathsetmacro{\zmin}{0}
  \pgfmathsetmacro{\zmax}{18}
  \pgfmathsetmacro{\cmin}{1.47}
  \pgfmathsetmacro{\cmax}{1.55}
  
\begin{axis}[%
width=\modelwidth, height=\modelheight,
axis on top, separate axis lines,
xmin=\xminloc, xmax=\xmaxloc, xlabel={$x$ (\si{\cm})},
ymin=\zminloc, ymax=\zmaxloc, ylabel={$z$ (\si{\cm})}, y dir=reverse,
xtick={4,10,12},
ytick={4,10,14},
xticklabel pos=right, xlabel near ticks,
x label style={xshift=-0.5cm, yshift=-0.50cm}, 
y label style={xshift= 0.5cm, yshift=-0.50cm},
hide axis,
  label style={font=\scriptsize},
  tick label style={font=\scriptsize},
  legend style={font=\scriptsize\selectfont},
]
\addplot [forget plot] graphics [xmin=\xmin,xmax=\xmax,ymin=\zmin,ymax=\zmax] {{\modelfile}.png};
\end{axis}
\end{tikzpicture}%
 \hspace*{\myhspace}
\renewcommand{\modelfile}{cp-600kHz-data-ma-7}\begin{tikzpicture}

  \pgfmathsetmacro{\xmin}{0}
  \pgfmathsetmacro{\xmax}{18}
  \pgfmathsetmacro{\zmin}{0}
  \pgfmathsetmacro{\zmax}{18}
  \pgfmathsetmacro{\cmin}{1.47}
  \pgfmathsetmacro{\cmax}{1.55}
  
\begin{axis}[%
width=\modelwidth, height=\modelheight,
axis on top, separate axis lines,
xmin=\xminloc, xmax=\xmaxloc, xlabel={$x$ (\si{\cm})},
ymin=\zminloc, ymax=\zmaxloc, ylabel={$z$ (\si{\cm})}, y dir=reverse,
xtick={4,10,12},
ytick={4,10,14},
xticklabel pos=right, xlabel near ticks,
x label style={xshift=-0.5cm, yshift=-0.50cm}, 
y label style={xshift= 0.5cm, yshift=-0.50cm},
hide axis,
  label style={font=\scriptsize},
  tick label style={font=\scriptsize},
  legend style={font=\scriptsize\selectfont},
]
\addplot [forget plot] graphics [xmin=\xmin,xmax=\xmax,ymin=\zmin,ymax=\zmax] {{\modelfile}.png};
\end{axis}
\end{tikzpicture}%
 \hspace*{\myhspace}
\renewcommand{\modelfile}{cp-600kHz-data-ma-2}\begin{tikzpicture}

  \pgfmathsetmacro{\xmin}{0}
  \pgfmathsetmacro{\xmax}{18}
  \pgfmathsetmacro{\zmin}{0}
  \pgfmathsetmacro{\zmax}{18}
  \pgfmathsetmacro{\cmin}{1.47}
  \pgfmathsetmacro{\cmax}{1.55}
  
\begin{axis}[%
width=\modelwidth, height=\modelheight,
axis on top, separate axis lines,
xmin=\xminloc, xmax=\xmaxloc, xlabel={$x$ (\si{\cm})},
ymin=\zminloc, ymax=\zmaxloc, ylabel={$z$ (\si{\cm})}, y dir=reverse,
xtick={4,10,12},
ytick={4,10,14},
xticklabel pos=right, xlabel near ticks,
x label style={xshift=-0.5cm, yshift=-0.50cm}, 
y label style={xshift= 0.5cm, yshift=-0.50cm},
hide axis,
  label style={font=\scriptsize},
  tick label style={font=\scriptsize},
  legend style={font=\scriptsize\selectfont},
]
\addplot [forget plot] graphics [xmin=\xmin,xmax=\xmax,ymin=\zmin,ymax=\zmax] {{\modelfile}.png};
\end{axis}
\end{tikzpicture}%
 \hspace*{\myhspace}
\renewcommand{\modelfile}{cp-600kHz-data-ma-6}\begin{tikzpicture}

  \pgfmathsetmacro{\xmin}{0}
  \pgfmathsetmacro{\xmax}{18}
  \pgfmathsetmacro{\zmin}{0}
  \pgfmathsetmacro{\zmax}{18}
  \pgfmathsetmacro{\cmin}{1.47}
  \pgfmathsetmacro{\cmax}{1.55}
  
\begin{axis}[%
width=\modelwidth, height=\modelheight,
axis on top, separate axis lines,
xmin=\xminloc, xmax=\xmaxloc, xlabel={$x$ (\si{\cm})},
ymin=\zminloc, ymax=\zmaxloc, ylabel={$z$ (\si{\cm})}, y dir=reverse,
xtick={4,10,12},
ytick={4,10,14},
xticklabel pos=right, xlabel near ticks,
x label style={xshift=-0.5cm, yshift=-0.50cm}, 
y label style={xshift= 0.5cm, yshift=-0.50cm},
hide axis,
  label style={font=\scriptsize},
  tick label style={font=\scriptsize},
  legend style={font=\scriptsize\selectfont},
]
\addplot [forget plot] graphics [xmin=\xmin,xmax=\xmax,ymin=\zmin,ymax=\zmax] {{\modelfile}.png};
\end{axis}
\end{tikzpicture}%
 \hspace*{\myhspace}
\renewcommand{\modelfile}{cp-600kHz-data-ma-3}\begin{tikzpicture}

  \pgfmathsetmacro{\xmin}{0}
  \pgfmathsetmacro{\xmax}{18}
  \pgfmathsetmacro{\zmin}{0}
  \pgfmathsetmacro{\zmax}{18}
  \pgfmathsetmacro{\cmin}{1.47}
  \pgfmathsetmacro{\cmax}{1.55}
  
\begin{axis}[%
width=\modelwidth, height=\modelheight,
axis on top, separate axis lines,
xmin=\xminloc, xmax=\xmaxloc, xlabel={$x$ (\si{\cm})},
ymin=\zminloc, ymax=\zmaxloc, ylabel={$z$ (\si{\cm})}, y dir=reverse,
xtick={4,10,12},
ytick={4,10,14},
xticklabel pos=right, xlabel near ticks,
x label style={xshift=-0.5cm, yshift=-0.50cm}, 
y label style={xshift= 0.5cm, yshift=-0.50cm},
hide axis,
  label style={font=\scriptsize},
  tick label style={font=\scriptsize},
  legend style={font=\scriptsize\selectfont},
]
\addplot [forget plot] graphics [xmin=\xmin,xmax=\xmax,ymin=\zmin,ymax=\zmax] {{\modelfile}.png};
\end{axis}
\end{tikzpicture}%
 \hspace*{\myhspace}
\renewcommand{\modelfile}{cp-600kHz-data-ma-4}\begin{tikzpicture}

  \pgfmathsetmacro{\xmin}{0}
  \pgfmathsetmacro{\xmax}{18}
  \pgfmathsetmacro{\zmin}{0}
  \pgfmathsetmacro{\zmax}{18}
  \pgfmathsetmacro{\cmin}{1.47}
  \pgfmathsetmacro{\cmax}{1.55}
  
\begin{axis}[%
width=\modelwidth, height=\modelheight,
axis on top, separate axis lines,
xmin=\xminloc, xmax=\xmaxloc, xlabel={$x$ (\si{\cm})},
ymin=\zminloc, ymax=\zmaxloc, ylabel={$z$ (\si{\cm})}, y dir=reverse,
xtick={4,10,12},
ytick={4,10,14},
xticklabel pos=right, xlabel near ticks,
x label style={xshift=-0.5cm, yshift=-0.50cm}, 
y label style={xshift= 0.5cm, yshift=-0.50cm},
hide axis,
  label style={font=\scriptsize},
  tick label style={font=\scriptsize},
  legend style={font=\scriptsize\selectfont},
]
\addplot [forget plot] graphics [xmin=\xmin,xmax=\xmax,ymin=\zmin,ymax=\zmax] {{\modelfile}.png};
\end{axis}
\end{tikzpicture}%
} \\[\myvspace]

\subfloat[][Reconstruction using the data-set 
            generated with the simplified Zener model.]{
\renewcommand{\modelfile}{cp-600kHz-data-ze-1}\begin{tikzpicture}

  \pgfmathsetmacro{\xmin}{0}
  \pgfmathsetmacro{\xmax}{18}
  \pgfmathsetmacro{\zmin}{0}
  \pgfmathsetmacro{\zmax}{18}
  \pgfmathsetmacro{\cmin}{1.47}
  \pgfmathsetmacro{\cmax}{1.55}
  
\begin{axis}[%
width=\modelwidth, height=\modelheight,
axis on top, separate axis lines,
xmin=\xminloc, xmax=\xmaxloc, xlabel={$x$ (\si{\cm})},
ymin=\zminloc, ymax=\zmaxloc, ylabel={$z$ (\si{\cm})}, y dir=reverse,
xtick={4,10,12},
ytick={4,10,14},
xticklabel pos=right, xlabel near ticks,
x label style={xshift=-0.5cm, yshift=-0.50cm}, 
y label style={xshift= 0.5cm, yshift=-0.50cm},
hide axis,
  label style={font=\scriptsize},
  tick label style={font=\scriptsize},
  legend style={font=\scriptsize\selectfont},
]
\addplot [forget plot] graphics [xmin=\xmin,xmax=\xmax,ymin=\zmin,ymax=\zmax] {{\modelfile}.png};
\end{axis}
\end{tikzpicture}%
 \hspace*{\myhspace}
\renewcommand{\modelfile}{cp-600kHz-data-ze-5}\begin{tikzpicture}

  \pgfmathsetmacro{\xmin}{0}
  \pgfmathsetmacro{\xmax}{18}
  \pgfmathsetmacro{\zmin}{0}
  \pgfmathsetmacro{\zmax}{18}
  \pgfmathsetmacro{\cmin}{1.47}
  \pgfmathsetmacro{\cmax}{1.55}
  
\begin{axis}[%
width=\modelwidth, height=\modelheight,
axis on top, separate axis lines,
xmin=\xminloc, xmax=\xmaxloc, xlabel={$x$ (\si{\cm})},
ymin=\zminloc, ymax=\zmaxloc, ylabel={$z$ (\si{\cm})}, y dir=reverse,
xtick={4,10,12},
ytick={4,10,14},
xticklabel pos=right, xlabel near ticks,
x label style={xshift=-0.5cm, yshift=-0.50cm}, 
y label style={xshift= 0.5cm, yshift=-0.50cm},
hide axis,
  label style={font=\scriptsize},
  tick label style={font=\scriptsize},
  legend style={font=\scriptsize\selectfont},
]
\addplot [forget plot] graphics [xmin=\xmin,xmax=\xmax,ymin=\zmin,ymax=\zmax] {{\modelfile}.png};
\end{axis}
\end{tikzpicture}%
 \hspace*{\myhspace}
\renewcommand{\modelfile}{cp-600kHz-data-ze-7}\begin{tikzpicture}

  \pgfmathsetmacro{\xmin}{0}
  \pgfmathsetmacro{\xmax}{18}
  \pgfmathsetmacro{\zmin}{0}
  \pgfmathsetmacro{\zmax}{18}
  \pgfmathsetmacro{\cmin}{1.47}
  \pgfmathsetmacro{\cmax}{1.55}
  
\begin{axis}[%
width=\modelwidth, height=\modelheight,
axis on top, separate axis lines,
xmin=\xminloc, xmax=\xmaxloc, xlabel={$x$ (\si{\cm})},
ymin=\zminloc, ymax=\zmaxloc, ylabel={$z$ (\si{\cm})}, y dir=reverse,
xtick={4,10,12},
ytick={4,10,14},
xticklabel pos=right, xlabel near ticks,
x label style={xshift=-0.5cm, yshift=-0.50cm}, 
y label style={xshift= 0.5cm, yshift=-0.50cm},
hide axis,
  label style={font=\scriptsize},
  tick label style={font=\scriptsize},
  legend style={font=\scriptsize\selectfont},
]
\addplot [forget plot] graphics [xmin=\xmin,xmax=\xmax,ymin=\zmin,ymax=\zmax] {{\modelfile}.png};
\end{axis}
\end{tikzpicture}%
 \hspace*{\myhspace}
\renewcommand{\modelfile}{cp-600kHz-data-ze-2}\begin{tikzpicture}

  \pgfmathsetmacro{\xmin}{0}
  \pgfmathsetmacro{\xmax}{18}
  \pgfmathsetmacro{\zmin}{0}
  \pgfmathsetmacro{\zmax}{18}
  \pgfmathsetmacro{\cmin}{1.47}
  \pgfmathsetmacro{\cmax}{1.55}
  
\begin{axis}[%
width=\modelwidth, height=\modelheight,
axis on top, separate axis lines,
xmin=\xminloc, xmax=\xmaxloc, xlabel={$x$ (\si{\cm})},
ymin=\zminloc, ymax=\zmaxloc, ylabel={$z$ (\si{\cm})}, y dir=reverse,
xtick={4,10,12},
ytick={4,10,14},
xticklabel pos=right, xlabel near ticks,
x label style={xshift=-0.5cm, yshift=-0.50cm}, 
y label style={xshift= 0.5cm, yshift=-0.50cm},
hide axis,
  label style={font=\scriptsize},
  tick label style={font=\scriptsize},
  legend style={font=\scriptsize\selectfont},
]
\addplot [forget plot] graphics [xmin=\xmin,xmax=\xmax,ymin=\zmin,ymax=\zmax] {{\modelfile}.png};
\end{axis}
\end{tikzpicture}%
 \hspace*{\myhspace}
\renewcommand{\modelfile}{cp-600kHz-data-ze-6}\begin{tikzpicture}

  \pgfmathsetmacro{\xmin}{0}
  \pgfmathsetmacro{\xmax}{18}
  \pgfmathsetmacro{\zmin}{0}
  \pgfmathsetmacro{\zmax}{18}
  \pgfmathsetmacro{\cmin}{1.47}
  \pgfmathsetmacro{\cmax}{1.55}
  
\begin{axis}[%
width=\modelwidth, height=\modelheight,
axis on top, separate axis lines,
xmin=\xminloc, xmax=\xmaxloc, xlabel={$x$ (\si{\cm})},
ymin=\zminloc, ymax=\zmaxloc, ylabel={$z$ (\si{\cm})}, y dir=reverse,
xtick={4,10,12},
ytick={4,10,14},
xticklabel pos=right, xlabel near ticks,
x label style={xshift=-0.5cm, yshift=-0.50cm}, 
y label style={xshift= 0.5cm, yshift=-0.50cm},
hide axis,
  label style={font=\scriptsize},
  tick label style={font=\scriptsize},
  legend style={font=\scriptsize\selectfont},
]
\addplot [forget plot] graphics [xmin=\xmin,xmax=\xmax,ymin=\zmin,ymax=\zmax] {{\modelfile}.png};
\end{axis}
\end{tikzpicture}%
 \hspace*{\myhspace}
\renewcommand{\modelfile}{cp-600kHz-data-ze-3}\begin{tikzpicture}

  \pgfmathsetmacro{\xmin}{0}
  \pgfmathsetmacro{\xmax}{18}
  \pgfmathsetmacro{\zmin}{0}
  \pgfmathsetmacro{\zmax}{18}
  \pgfmathsetmacro{\cmin}{1.47}
  \pgfmathsetmacro{\cmax}{1.55}
  
\begin{axis}[%
width=\modelwidth, height=\modelheight,
axis on top, separate axis lines,
xmin=\xminloc, xmax=\xmaxloc, xlabel={$x$ (\si{\cm})},
ymin=\zminloc, ymax=\zmaxloc, ylabel={$z$ (\si{\cm})}, y dir=reverse,
xtick={4,10,12},
ytick={4,10,14},
xticklabel pos=right, xlabel near ticks,
x label style={xshift=-0.5cm, yshift=-0.50cm}, 
y label style={xshift= 0.5cm, yshift=-0.50cm},
hide axis,
  label style={font=\scriptsize},
  tick label style={font=\scriptsize},
  legend style={font=\scriptsize\selectfont},
]
\addplot [forget plot] graphics [xmin=\xmin,xmax=\xmax,ymin=\zmin,ymax=\zmax] {{\modelfile}.png};
\end{axis}
\end{tikzpicture}%
 \hspace*{\myhspace}
\renewcommand{\modelfile}{cp-600kHz-data-ze-4}\begin{tikzpicture}

  \pgfmathsetmacro{\xmin}{0}
  \pgfmathsetmacro{\xmax}{18}
  \pgfmathsetmacro{\zmin}{0}
  \pgfmathsetmacro{\zmax}{18}
  \pgfmathsetmacro{\cmin}{1.47}
  \pgfmathsetmacro{\cmax}{1.55}
  
\begin{axis}[%
width=\modelwidth, height=\modelheight,
axis on top, separate axis lines,
xmin=\xminloc, xmax=\xmaxloc, xlabel={$x$ (\si{\cm})},
ymin=\zminloc, ymax=\zmaxloc, ylabel={$z$ (\si{\cm})}, y dir=reverse,
xtick={4,10,12},
ytick={4,10,14},
xticklabel pos=right, xlabel near ticks,
x label style={xshift=-0.5cm, yshift=-0.50cm}, 
y label style={xshift= 0.5cm, yshift=-0.50cm},
hide axis,
  label style={font=\scriptsize},
  tick label style={font=\scriptsize},
  legend style={font=\scriptsize\selectfont},
]
\addplot [forget plot] graphics [xmin=\xmin,xmax=\xmax,ymin=\zmin,ymax=\zmax] {{\modelfile}.png};
\end{axis}
\end{tikzpicture}%
}  \\[\myvspace]

\subfloat[][Reconstruction using the data-set 
            generated with the simplified Cole--Cole model.]{
\renewcommand{\modelfile}{cp-600kHz-data-cc-1}\begin{tikzpicture}

  \pgfmathsetmacro{\xmin}{0}
  \pgfmathsetmacro{\xmax}{18}
  \pgfmathsetmacro{\zmin}{0}
  \pgfmathsetmacro{\zmax}{18}
  \pgfmathsetmacro{\cmin}{1.47}
  \pgfmathsetmacro{\cmax}{1.55}
  
\begin{axis}[%
width=\modelwidth, height=\modelheight,
axis on top, separate axis lines,
xmin=\xminloc, xmax=\xmaxloc, xlabel={$x$ (\si{\cm})},
ymin=\zminloc, ymax=\zmaxloc, ylabel={$z$ (\si{\cm})}, y dir=reverse,
xtick={4,10,12},
ytick={4,10,14},
xticklabel pos=right, xlabel near ticks,
x label style={xshift=-0.5cm, yshift=-0.50cm}, 
y label style={xshift= 0.5cm, yshift=-0.50cm},
hide axis,
  label style={font=\scriptsize},
  tick label style={font=\scriptsize},
  legend style={font=\scriptsize\selectfont},
]
\addplot [forget plot] graphics [xmin=\xmin,xmax=\xmax,ymin=\zmin,ymax=\zmax] {{\modelfile}.png};
\end{axis}
\end{tikzpicture}%
 \hspace*{\myhspace}
\renewcommand{\modelfile}{cp-600kHz-data-cc-5}\begin{tikzpicture}

  \pgfmathsetmacro{\xmin}{0}
  \pgfmathsetmacro{\xmax}{18}
  \pgfmathsetmacro{\zmin}{0}
  \pgfmathsetmacro{\zmax}{18}
  \pgfmathsetmacro{\cmin}{1.47}
  \pgfmathsetmacro{\cmax}{1.55}
  
\begin{axis}[%
width=\modelwidth, height=\modelheight,
axis on top, separate axis lines,
xmin=\xminloc, xmax=\xmaxloc, xlabel={$x$ (\si{\cm})},
ymin=\zminloc, ymax=\zmaxloc, ylabel={$z$ (\si{\cm})}, y dir=reverse,
xtick={4,10,12},
ytick={4,10,14},
xticklabel pos=right, xlabel near ticks,
x label style={xshift=-0.5cm, yshift=-0.50cm}, 
y label style={xshift= 0.5cm, yshift=-0.50cm},
hide axis,
  label style={font=\scriptsize},
  tick label style={font=\scriptsize},
  legend style={font=\scriptsize\selectfont},
]
\addplot [forget plot] graphics [xmin=\xmin,xmax=\xmax,ymin=\zmin,ymax=\zmax] {{\modelfile}.png};
\end{axis}
\end{tikzpicture}%
 \hspace*{\myhspace}
\renewcommand{\modelfile}{cp-600kHz-data-cc-7}\begin{tikzpicture}

  \pgfmathsetmacro{\xmin}{0}
  \pgfmathsetmacro{\xmax}{18}
  \pgfmathsetmacro{\zmin}{0}
  \pgfmathsetmacro{\zmax}{18}
  \pgfmathsetmacro{\cmin}{1.47}
  \pgfmathsetmacro{\cmax}{1.55}
  
\begin{axis}[%
width=\modelwidth, height=\modelheight,
axis on top, separate axis lines,
xmin=\xminloc, xmax=\xmaxloc, xlabel={$x$ (\si{\cm})},
ymin=\zminloc, ymax=\zmaxloc, ylabel={$z$ (\si{\cm})}, y dir=reverse,
xtick={4,10,12},
ytick={4,10,14},
xticklabel pos=right, xlabel near ticks,
x label style={xshift=-0.5cm, yshift=-0.50cm}, 
y label style={xshift= 0.5cm, yshift=-0.50cm},
hide axis,
  label style={font=\scriptsize},
  tick label style={font=\scriptsize},
  legend style={font=\scriptsize\selectfont},
]
\addplot [forget plot] graphics [xmin=\xmin,xmax=\xmax,ymin=\zmin,ymax=\zmax] {{\modelfile}.png};
\end{axis}
\end{tikzpicture}%
 \hspace*{\myhspace}
\renewcommand{\modelfile}{cp-600kHz-data-cc-2}\begin{tikzpicture}

  \pgfmathsetmacro{\xmin}{0}
  \pgfmathsetmacro{\xmax}{18}
  \pgfmathsetmacro{\zmin}{0}
  \pgfmathsetmacro{\zmax}{18}
  \pgfmathsetmacro{\cmin}{1.47}
  \pgfmathsetmacro{\cmax}{1.55}
  
\begin{axis}[%
width=\modelwidth, height=\modelheight,
axis on top, separate axis lines,
xmin=\xminloc, xmax=\xmaxloc, xlabel={$x$ (\si{\cm})},
ymin=\zminloc, ymax=\zmaxloc, ylabel={$z$ (\si{\cm})}, y dir=reverse,
xtick={4,10,12},
ytick={4,10,14},
xticklabel pos=right, xlabel near ticks,
x label style={xshift=-0.5cm, yshift=-0.50cm}, 
y label style={xshift= 0.5cm, yshift=-0.50cm},
hide axis,
  label style={font=\scriptsize},
  tick label style={font=\scriptsize},
  legend style={font=\scriptsize\selectfont},
]
\addplot [forget plot] graphics [xmin=\xmin,xmax=\xmax,ymin=\zmin,ymax=\zmax] {{\modelfile}.png};
\end{axis}
\end{tikzpicture}%
 \hspace*{\myhspace}
\renewcommand{\modelfile}{cp-600kHz-data-cc-6}\begin{tikzpicture}

  \pgfmathsetmacro{\xmin}{0}
  \pgfmathsetmacro{\xmax}{18}
  \pgfmathsetmacro{\zmin}{0}
  \pgfmathsetmacro{\zmax}{18}
  \pgfmathsetmacro{\cmin}{1.47}
  \pgfmathsetmacro{\cmax}{1.55}
  
\begin{axis}[%
width=\modelwidth, height=\modelheight,
axis on top, separate axis lines,
xmin=\xminloc, xmax=\xmaxloc, xlabel={$x$ (\si{\cm})},
ymin=\zminloc, ymax=\zmaxloc, ylabel={$z$ (\si{\cm})}, y dir=reverse,
xtick={4,10,12},
ytick={4,10,14},
xticklabel pos=right, xlabel near ticks,
x label style={xshift=-0.5cm, yshift=-0.50cm}, 
y label style={xshift= 0.5cm, yshift=-0.50cm},
hide axis,
  label style={font=\scriptsize},
  tick label style={font=\scriptsize},
  legend style={font=\scriptsize\selectfont},
]
\addplot [forget plot] graphics [xmin=\xmin,xmax=\xmax,ymin=\zmin,ymax=\zmax] {{\modelfile}.png};
\end{axis}
\end{tikzpicture}%
 \hspace*{\myhspace}
\renewcommand{\modelfile}{cp-600kHz-data-cc-3}\begin{tikzpicture}

  \pgfmathsetmacro{\xmin}{0}
  \pgfmathsetmacro{\xmax}{18}
  \pgfmathsetmacro{\zmin}{0}
  \pgfmathsetmacro{\zmax}{18}
  \pgfmathsetmacro{\cmin}{1.47}
  \pgfmathsetmacro{\cmax}{1.55}
  
\begin{axis}[%
width=\modelwidth, height=\modelheight,
axis on top, separate axis lines,
xmin=\xminloc, xmax=\xmaxloc, xlabel={$x$ (\si{\cm})},
ymin=\zminloc, ymax=\zmaxloc, ylabel={$z$ (\si{\cm})}, y dir=reverse,
xtick={4,10,12},
ytick={4,10,14},
xticklabel pos=right, xlabel near ticks,
x label style={xshift=-0.5cm, yshift=-0.50cm}, 
y label style={xshift= 0.5cm, yshift=-0.50cm},
hide axis,
  label style={font=\scriptsize},
  tick label style={font=\scriptsize},
  legend style={font=\scriptsize\selectfont},
]
\addplot [forget plot] graphics [xmin=\xmin,xmax=\xmax,ymin=\zmin,ymax=\zmax] {{\modelfile}.png};
\end{axis}
\end{tikzpicture}%
 \hspace*{\myhspace}
\renewcommand{\modelfile}{cp-600kHz-data-cc-4}\begin{tikzpicture}

  \pgfmathsetmacro{\xmin}{0}
  \pgfmathsetmacro{\xmax}{18}
  \pgfmathsetmacro{\zmin}{0}
  \pgfmathsetmacro{\zmax}{18}
  \pgfmathsetmacro{\cmin}{1.47}
  \pgfmathsetmacro{\cmax}{1.55}
  
\begin{axis}[%
width=\modelwidth, height=\modelheight,
axis on top, separate axis lines,
xmin=\xminloc, xmax=\xmaxloc, xlabel={$x$ (\si{\cm})},
ymin=\zminloc, ymax=\zmaxloc, ylabel={$z$ (\si{\cm})}, y dir=reverse,
xtick={4,10,12},
ytick={4,10,14},
xticklabel pos=right, xlabel near ticks,
x label style={xshift=-0.5cm, yshift=-0.50cm}, 
y label style={xshift= 0.5cm, yshift=-0.50cm},
hide axis,
  label style={font=\scriptsize},
  tick label style={font=\scriptsize},
  legend style={font=\scriptsize\selectfont},
]
\addplot [forget plot] graphics [xmin=\xmin,xmax=\xmax,ymin=\zmin,ymax=\zmax] {{\modelfile}.png};
\end{axis}
\end{tikzpicture}%
}  \\[\myvspace]

\subfloat[][Reconstruction using the data-set 
            generated with the simplified KSB model.]{
\renewcommand{\modelfile}{cp-600kHz-data-ks-1}\begin{tikzpicture}

  \pgfmathsetmacro{\xmin}{0}
  \pgfmathsetmacro{\xmax}{18}
  \pgfmathsetmacro{\zmin}{0}
  \pgfmathsetmacro{\zmax}{18}
  \pgfmathsetmacro{\cmin}{1.47}
  \pgfmathsetmacro{\cmax}{1.55}
  
\begin{axis}[%
width=\modelwidth, height=\modelheight,
axis on top, separate axis lines,
xmin=\xminloc, xmax=\xmaxloc, xlabel={$x$ (\si{\cm})},
ymin=\zminloc, ymax=\zmaxloc, ylabel={$z$ (\si{\cm})}, y dir=reverse,
xtick={4,10,12},
ytick={4,10,14},
xticklabel pos=right, xlabel near ticks,
x label style={xshift=-0.5cm, yshift=-0.50cm}, 
y label style={xshift= 0.5cm, yshift=-0.50cm},
hide axis,
  label style={font=\scriptsize},
  tick label style={font=\scriptsize},
  legend style={font=\scriptsize\selectfont},
]
\addplot [forget plot] graphics [xmin=\xmin,xmax=\xmax,ymin=\zmin,ymax=\zmax] {{\modelfile}.png};
\end{axis}
\end{tikzpicture}%
 \hspace*{\myhspace}
\renewcommand{\modelfile}{cp-600kHz-data-ks-5}\begin{tikzpicture}

  \pgfmathsetmacro{\xmin}{0}
  \pgfmathsetmacro{\xmax}{18}
  \pgfmathsetmacro{\zmin}{0}
  \pgfmathsetmacro{\zmax}{18}
  \pgfmathsetmacro{\cmin}{1.47}
  \pgfmathsetmacro{\cmax}{1.55}
  
\begin{axis}[%
width=\modelwidth, height=\modelheight,
axis on top, separate axis lines,
xmin=\xminloc, xmax=\xmaxloc, xlabel={$x$ (\si{\cm})},
ymin=\zminloc, ymax=\zmaxloc, ylabel={$z$ (\si{\cm})}, y dir=reverse,
xtick={4,10,12},
ytick={4,10,14},
xticklabel pos=right, xlabel near ticks,
x label style={xshift=-0.5cm, yshift=-0.50cm}, 
y label style={xshift= 0.5cm, yshift=-0.50cm},
hide axis,
  label style={font=\scriptsize},
  tick label style={font=\scriptsize},
  legend style={font=\scriptsize\selectfont},
]
\addplot [forget plot] graphics [xmin=\xmin,xmax=\xmax,ymin=\zmin,ymax=\zmax] {{\modelfile}.png};
\end{axis}
\end{tikzpicture}%
 \hspace*{\myhspace}
\renewcommand{\modelfile}{cp-600kHz-data-ks-7}\begin{tikzpicture}

  \pgfmathsetmacro{\xmin}{0}
  \pgfmathsetmacro{\xmax}{18}
  \pgfmathsetmacro{\zmin}{0}
  \pgfmathsetmacro{\zmax}{18}
  \pgfmathsetmacro{\cmin}{1.47}
  \pgfmathsetmacro{\cmax}{1.55}
  
\begin{axis}[%
width=\modelwidth, height=\modelheight,
axis on top, separate axis lines,
xmin=\xminloc, xmax=\xmaxloc, xlabel={$x$ (\si{\cm})},
ymin=\zminloc, ymax=\zmaxloc, ylabel={$z$ (\si{\cm})}, y dir=reverse,
xtick={4,10,12},
ytick={4,10,14},
xticklabel pos=right, xlabel near ticks,
x label style={xshift=-0.5cm, yshift=-0.50cm}, 
y label style={xshift= 0.5cm, yshift=-0.50cm},
hide axis,
  label style={font=\scriptsize},
  tick label style={font=\scriptsize},
  legend style={font=\scriptsize\selectfont},
]
\addplot [forget plot] graphics [xmin=\xmin,xmax=\xmax,ymin=\zmin,ymax=\zmax] {{\modelfile}.png};
\end{axis}
\end{tikzpicture}%
 \hspace*{\myhspace}
\renewcommand{\modelfile}{cp-600kHz-data-ks-2}\begin{tikzpicture}

  \pgfmathsetmacro{\xmin}{0}
  \pgfmathsetmacro{\xmax}{18}
  \pgfmathsetmacro{\zmin}{0}
  \pgfmathsetmacro{\zmax}{18}
  \pgfmathsetmacro{\cmin}{1.47}
  \pgfmathsetmacro{\cmax}{1.55}
  
\begin{axis}[%
width=\modelwidth, height=\modelheight,
axis on top, separate axis lines,
xmin=\xminloc, xmax=\xmaxloc, xlabel={$x$ (\si{\cm})},
ymin=\zminloc, ymax=\zmaxloc, ylabel={$z$ (\si{\cm})}, y dir=reverse,
xtick={4,10,12},
ytick={4,10,14},
xticklabel pos=right, xlabel near ticks,
x label style={xshift=-0.5cm, yshift=-0.50cm}, 
y label style={xshift= 0.5cm, yshift=-0.50cm},
hide axis,
  label style={font=\scriptsize},
  tick label style={font=\scriptsize},
  legend style={font=\scriptsize\selectfont},
]
\addplot [forget plot] graphics [xmin=\xmin,xmax=\xmax,ymin=\zmin,ymax=\zmax] {{\modelfile}.png};
\end{axis}
\end{tikzpicture}%
 \hspace*{\myhspace}
\renewcommand{\modelfile}{cp-600kHz-data-ks-6}\begin{tikzpicture}

  \pgfmathsetmacro{\xmin}{0}
  \pgfmathsetmacro{\xmax}{18}
  \pgfmathsetmacro{\zmin}{0}
  \pgfmathsetmacro{\zmax}{18}
  \pgfmathsetmacro{\cmin}{1.47}
  \pgfmathsetmacro{\cmax}{1.55}
  
\begin{axis}[%
width=\modelwidth, height=\modelheight,
axis on top, separate axis lines,
xmin=\xminloc, xmax=\xmaxloc, xlabel={$x$ (\si{\cm})},
ymin=\zminloc, ymax=\zmaxloc, ylabel={$z$ (\si{\cm})}, y dir=reverse,
xtick={4,10,12},
ytick={4,10,14},
xticklabel pos=right, xlabel near ticks,
x label style={xshift=-0.5cm, yshift=-0.50cm}, 
y label style={xshift= 0.5cm, yshift=-0.50cm},
hide axis,
  label style={font=\scriptsize},
  tick label style={font=\scriptsize},
  legend style={font=\scriptsize\selectfont},
]
\addplot [forget plot] graphics [xmin=\xmin,xmax=\xmax,ymin=\zmin,ymax=\zmax] {{\modelfile}.png};
\end{axis}
\end{tikzpicture}%
 \hspace*{\myhspace}
\renewcommand{\modelfile}{cp-600kHz-data-ks-3}\begin{tikzpicture}

  \pgfmathsetmacro{\xmin}{0}
  \pgfmathsetmacro{\xmax}{18}
  \pgfmathsetmacro{\zmin}{0}
  \pgfmathsetmacro{\zmax}{18}
  \pgfmathsetmacro{\cmin}{1.47}
  \pgfmathsetmacro{\cmax}{1.55}
  
\begin{axis}[%
width=\modelwidth, height=\modelheight,
axis on top, separate axis lines,
xmin=\xminloc, xmax=\xmaxloc, xlabel={$x$ (\si{\cm})},
ymin=\zminloc, ymax=\zmaxloc, ylabel={$z$ (\si{\cm})}, y dir=reverse,
xtick={4,10,12},
ytick={4,10,14},
xticklabel pos=right, xlabel near ticks,
x label style={xshift=-0.5cm, yshift=-0.50cm}, 
y label style={xshift= 0.5cm, yshift=-0.50cm},
hide axis,
  label style={font=\scriptsize},
  tick label style={font=\scriptsize},
  legend style={font=\scriptsize\selectfont},
]
\addplot [forget plot] graphics [xmin=\xmin,xmax=\xmax,ymin=\zmin,ymax=\zmax] {{\modelfile}.png};
\end{axis}
\end{tikzpicture}%
 \hspace*{\myhspace}
\renewcommand{\modelfile}{cp-600kHz-data-ks-4}\begin{tikzpicture}

  \pgfmathsetmacro{\xmin}{0}
  \pgfmathsetmacro{\xmax}{18}
  \pgfmathsetmacro{\zmin}{0}
  \pgfmathsetmacro{\zmax}{18}
  \pgfmathsetmacro{\cmin}{1.47}
  \pgfmathsetmacro{\cmax}{1.55}
  
\begin{axis}[%
width=\modelwidth, height=\modelheight,
axis on top, separate axis lines,
xmin=\xminloc, xmax=\xmaxloc, xlabel={$x$ (\si{\cm})},
ymin=\zminloc, ymax=\zmaxloc, ylabel={$z$ (\si{\cm})}, y dir=reverse,
xtick={4,10,12},
ytick={4,10,14},
xticklabel pos=right, xlabel near ticks,
x label style={xshift=-0.5cm, yshift=-0.50cm}, 
y label style={xshift= 0.5cm, yshift=-0.50cm},
hide axis,
  label style={font=\scriptsize},
  tick label style={font=\scriptsize},
  legend style={font=\scriptsize\selectfont},
]
\addplot [forget plot] graphics [xmin=\xmin,xmax=\xmax,ymin=\zmin,ymax=\zmax] {{\modelfile}.png};
\end{axis}
\end{tikzpicture}%
} \\[\myvspace]

\subfloat[][Reconstruction using the data-set 
            generated with the simplified modified Szabo model.]{
\renewcommand{\modelfile}{cp-600kHz-data-sz-1}\begin{tikzpicture}

  \pgfmathsetmacro{\xmin}{0}
  \pgfmathsetmacro{\xmax}{18}
  \pgfmathsetmacro{\zmin}{0}
  \pgfmathsetmacro{\zmax}{18}
  \pgfmathsetmacro{\cmin}{1.47}
  \pgfmathsetmacro{\cmax}{1.55}
  
\begin{axis}[%
width=\modelwidth, height=\modelheight,
axis on top, separate axis lines,
xmin=\xminloc, xmax=\xmaxloc, xlabel={$x$ (\si{\cm})},
ymin=\zminloc, ymax=\zmaxloc, ylabel={$z$ (\si{\cm})}, y dir=reverse,
xtick={4,10,12},
ytick={4,10,14},
xticklabel pos=right, xlabel near ticks,
x label style={xshift=-0.5cm, yshift=-0.50cm}, 
y label style={xshift= 0.5cm, yshift=-0.50cm},
hide axis,
  label style={font=\scriptsize},
  tick label style={font=\scriptsize},
  legend style={font=\scriptsize\selectfont},
]
\addplot [forget plot] graphics [xmin=\xmin,xmax=\xmax,ymin=\zmin,ymax=\zmax] {{\modelfile}.png};
\end{axis}
\end{tikzpicture}%
 \hspace*{\myhspace}
\renewcommand{\modelfile}{cp-600kHz-data-sz-5}\begin{tikzpicture}

  \pgfmathsetmacro{\xmin}{0}
  \pgfmathsetmacro{\xmax}{18}
  \pgfmathsetmacro{\zmin}{0}
  \pgfmathsetmacro{\zmax}{18}
  \pgfmathsetmacro{\cmin}{1.47}
  \pgfmathsetmacro{\cmax}{1.55}
  
\begin{axis}[%
width=\modelwidth, height=\modelheight,
axis on top, separate axis lines,
xmin=\xminloc, xmax=\xmaxloc, xlabel={$x$ (\si{\cm})},
ymin=\zminloc, ymax=\zmaxloc, ylabel={$z$ (\si{\cm})}, y dir=reverse,
xtick={4,10,12},
ytick={4,10,14},
xticklabel pos=right, xlabel near ticks,
x label style={xshift=-0.5cm, yshift=-0.50cm}, 
y label style={xshift= 0.5cm, yshift=-0.50cm},
hide axis,
  label style={font=\scriptsize},
  tick label style={font=\scriptsize},
  legend style={font=\scriptsize\selectfont},
]
\addplot [forget plot] graphics [xmin=\xmin,xmax=\xmax,ymin=\zmin,ymax=\zmax] {{\modelfile}.png};
\end{axis}
\end{tikzpicture}%
 \hspace*{\myhspace}
\renewcommand{\modelfile}{cp-600kHz-data-sz-7}\begin{tikzpicture}

  \pgfmathsetmacro{\xmin}{0}
  \pgfmathsetmacro{\xmax}{18}
  \pgfmathsetmacro{\zmin}{0}
  \pgfmathsetmacro{\zmax}{18}
  \pgfmathsetmacro{\cmin}{1.47}
  \pgfmathsetmacro{\cmax}{1.55}
  
\begin{axis}[%
width=\modelwidth, height=\modelheight,
axis on top, separate axis lines,
xmin=\xminloc, xmax=\xmaxloc, xlabel={$x$ (\si{\cm})},
ymin=\zminloc, ymax=\zmaxloc, ylabel={$z$ (\si{\cm})}, y dir=reverse,
xtick={4,10,12},
ytick={4,10,14},
xticklabel pos=right, xlabel near ticks,
x label style={xshift=-0.5cm, yshift=-0.50cm}, 
y label style={xshift= 0.5cm, yshift=-0.50cm},
hide axis,
  label style={font=\scriptsize},
  tick label style={font=\scriptsize},
  legend style={font=\scriptsize\selectfont},
]
\addplot [forget plot] graphics [xmin=\xmin,xmax=\xmax,ymin=\zmin,ymax=\zmax] {{\modelfile}.png};
\end{axis}
\end{tikzpicture}%
 \hspace*{\myhspace}
\renewcommand{\modelfile}{cp-600kHz-data-sz-2}\begin{tikzpicture}

  \pgfmathsetmacro{\xmin}{0}
  \pgfmathsetmacro{\xmax}{18}
  \pgfmathsetmacro{\zmin}{0}
  \pgfmathsetmacro{\zmax}{18}
  \pgfmathsetmacro{\cmin}{1.47}
  \pgfmathsetmacro{\cmax}{1.55}
  
\begin{axis}[%
width=\modelwidth, height=\modelheight,
axis on top, separate axis lines,
xmin=\xminloc, xmax=\xmaxloc, xlabel={$x$ (\si{\cm})},
ymin=\zminloc, ymax=\zmaxloc, ylabel={$z$ (\si{\cm})}, y dir=reverse,
xtick={4,10,12},
ytick={4,10,14},
xticklabel pos=right, xlabel near ticks,
x label style={xshift=-0.5cm, yshift=-0.50cm}, 
y label style={xshift= 0.5cm, yshift=-0.50cm},
hide axis,
  label style={font=\scriptsize},
  tick label style={font=\scriptsize},
  legend style={font=\scriptsize\selectfont},
]
\addplot [forget plot] graphics [xmin=\xmin,xmax=\xmax,ymin=\zmin,ymax=\zmax] {{\modelfile}.png};
\end{axis}
\end{tikzpicture}%
 \hspace*{\myhspace}
\renewcommand{\modelfile}{cp-600kHz-data-sz-6}\begin{tikzpicture}

  \pgfmathsetmacro{\xmin}{0}
  \pgfmathsetmacro{\xmax}{18}
  \pgfmathsetmacro{\zmin}{0}
  \pgfmathsetmacro{\zmax}{18}
  \pgfmathsetmacro{\cmin}{1.47}
  \pgfmathsetmacro{\cmax}{1.55}
  
\begin{axis}[%
width=\modelwidth, height=\modelheight,
axis on top, separate axis lines,
xmin=\xminloc, xmax=\xmaxloc, xlabel={$x$ (\si{\cm})},
ymin=\zminloc, ymax=\zmaxloc, ylabel={$z$ (\si{\cm})}, y dir=reverse,
xtick={4,10,12},
ytick={4,10,14},
xticklabel pos=right, xlabel near ticks,
x label style={xshift=-0.5cm, yshift=-0.50cm}, 
y label style={xshift= 0.5cm, yshift=-0.50cm},
hide axis,
  label style={font=\scriptsize},
  tick label style={font=\scriptsize},
  legend style={font=\scriptsize\selectfont},
]
\addplot [forget plot] graphics [xmin=\xmin,xmax=\xmax,ymin=\zmin,ymax=\zmax] {{\modelfile}.png};
\end{axis}
\end{tikzpicture}%
 \hspace*{\myhspace}
\renewcommand{\modelfile}{cp-600kHz-data-sz-3}\begin{tikzpicture}

  \pgfmathsetmacro{\xmin}{0}
  \pgfmathsetmacro{\xmax}{18}
  \pgfmathsetmacro{\zmin}{0}
  \pgfmathsetmacro{\zmax}{18}
  \pgfmathsetmacro{\cmin}{1.47}
  \pgfmathsetmacro{\cmax}{1.55}
  
\begin{axis}[%
width=\modelwidth, height=\modelheight,
axis on top, separate axis lines,
xmin=\xminloc, xmax=\xmaxloc, xlabel={$x$ (\si{\cm})},
ymin=\zminloc, ymax=\zmaxloc, ylabel={$z$ (\si{\cm})}, y dir=reverse,
xtick={4,10,12},
ytick={4,10,14},
xticklabel pos=right, xlabel near ticks,
x label style={xshift=-0.5cm, yshift=-0.50cm}, 
y label style={xshift= 0.5cm, yshift=-0.50cm},
hide axis,
  label style={font=\scriptsize},
  tick label style={font=\scriptsize},
  legend style={font=\scriptsize\selectfont},
]
\addplot [forget plot] graphics [xmin=\xmin,xmax=\xmax,ymin=\zmin,ymax=\zmax] {{\modelfile}.png};
\end{axis}
\end{tikzpicture}%
 \hspace*{\myhspace}
\renewcommand{\modelfile}{cp-600kHz-data-sz-4}\begin{tikzpicture}

  \pgfmathsetmacro{\xmin}{0}
  \pgfmathsetmacro{\xmax}{18}
  \pgfmathsetmacro{\zmin}{0}
  \pgfmathsetmacro{\zmax}{18}
  \pgfmathsetmacro{\cmin}{1.47}
  \pgfmathsetmacro{\cmax}{1.55}
  
\begin{axis}[%
width=\modelwidth, height=\modelheight,
axis on top, separate axis lines,
xmin=\xminloc, xmax=\xmaxloc, xlabel={$x$ (\si{\cm})},
ymin=\zminloc, ymax=\zmaxloc, ylabel={$z$ (\si{\cm})}, y dir=reverse,
xtick={4,10,12},
ytick={4,10,14},
xticklabel pos=right, xlabel near ticks,
x label style={xshift=-0.5cm, yshift=-0.50cm}, 
y label style={xshift= 0.5cm, yshift=-0.50cm},
hide axis,
  label style={font=\scriptsize},
  tick label style={font=\scriptsize},
  legend style={font=\scriptsize\selectfont},
]
\addplot [forget plot] graphics [xmin=\xmin,xmax=\xmax,ymin=\zmin,ymax=\zmax] {{\modelfile}.png};
\end{axis}
\end{tikzpicture}%
}  \\[-0.25em]
\caption{Reconstruction of the wave speed $c_0=\sqrt{\kappa_0/\rho}$
          using the iterative minimization 
         \cref{algo:FWI} with sequential frequency progression 
         $\omegaR/(2\pi)=\{200, \,300, \, 400, \, \num{600}\}$ \si{\kilo\Hz} and
         fixed $\omegaI=0$, starting from homogeneous 
         backgrounds (\cref{table:interval}).
         The visualization corresponds to $x$ from 
         $2$ to $16$ \si{\cm} and $z$ from $4$ to 
         \num{14.5} \si{\cm}, with
         color scale similar to \ref{fig:models:true}.
         Each line corresponds to a different data-set, 
         each column corresponds to a different 
         attenuation model used for FWI:
         (column 1) Kolsky--Futterman model;
         (column 2) Kelvin--Voigt model;
         (column 3) Maxwell model;
         (column 4) Zener model;
         (column 5) Cole--Cole model;
         (column 6) KSB model;
         (column 7) Szabo model. 
         Diagonal elements are reconstructions 
         with the ``right'' model.}
\label{fig:fwi2d:abosrbing:all}
\end{figure}
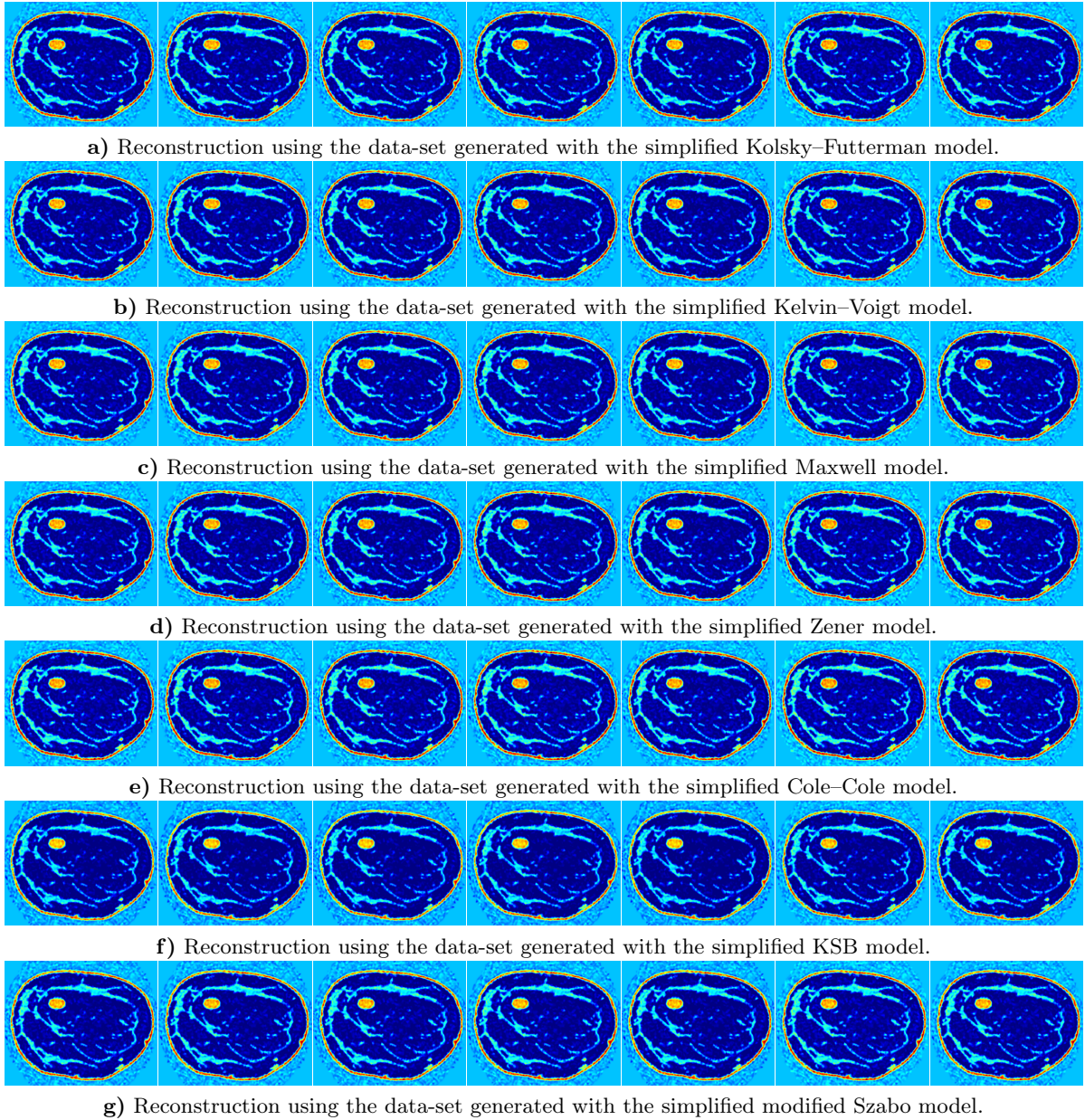

\bigskip

While each attenuation model gives a quite different 
data-set, as highlighted in \cref{fig:data}, we see 
that all of the configurations 
provide the appropriate layers in the sample, and have similar accuracy. 
The main tissue features
are reconstructed with the appropriate values
and the skin of the sample is correctly obtained.
Namely, using an attenuation model for the reconstruction that 
is different from the one of the original medium \emph{does not} 
alter the accuracy of the features reconstructed, 
demonstrating the robustness of the iterative minimization 
procedure for this two-dimensional test-case.
In addition, we see that while the density $\rho$ 
and attenuation $Q$ are not inverted and remain to their 
initial values, the reconstructed profiles of wave speed
$c_0$ is accurate.

In all of the cases, the ellipse-shaped included defect
is found at the right position and with its precise shape. 
We note that the speed in the skin layer in slightly less than
expected, while the one in the contrasting object is slightly
higher than expected. 
The robustness of the reconstruction algorithm with respect
to changes in attenuation models could be explained by two 
reasons: 
\vspace*{-0.5em}
\begin{enumerate} \setlength{\itemsep}{-2pt}
  \item Data of relatively limited frequency bandwidth are 
        sufficient for the reconstruction (from \num{200} 
        to \num{600} \si{\kilo\Hz}) and the difference 
        between attenuation models is less significant is
        such a narrow band, see \cref{fig:plot:comparison-attenuation}.
  \item The reconstruction with FWI is known to be more 
        sensitive to phase shifts rather than amplitude,
        hence robust with respect to incorrect attenuation
        which mostly affects the signal amplitude.
        Indeed, the dispersion resulting from the change of
        attenuation models remains small compared to the one
        that comes from changes in the bulk modulus $\kappa_0$. 
\end{enumerate}
In addition, we are in a case of tissues which 
have weak attenuation properties, which may also
justifies the lack of sensitivity to change in attenuation
model, even though the wave fields of \cref{fig:data}
show strong differences.
In the following experiments, for the sake of clarity, 
we only consider one of the cases, with the data-set 
generated with the Kolsky--Futterman model, while FWI 
is conducted using the Kelvin--Voigt model.

\begin{remark}[Post-processing visualization] \label{rk:visualization}
  The visualization of the reconstruction may suffer 
  from the coarse discretized mesh (about \num{10000} cells)
  employed for the numerical discretization, which is 
  chosen to reduce the numerical cost. 
  To enhance the visualization, we use the function 
  \texttt{imgaussfilt} of MATLAB, which corresponds 
  to a Gaussian smoothing filter, that we illustrate 
  in \cref{fig:fwi2d:visualization}. 
  This procedure is done a posteriori, independently of the reconstruction 
  algorithm, and is therefore costless. 
  One could instead rely on an extra regularization term in 
  the optimization, such that Tikhonov or Total Variation 
  regularization, but it needs an additional threshold in 
  the formulation of the minimization problem, which can
  be difficult to select, e.g., \cite{Faucher2020EV}. 
  This post-processing is only performed to enhance 
  the visual aspect of the reconstructions by smoothing 
  some numerical artifacts.
  \setlength{\modelwidth} {5.10cm}
  \setlength{\modelheight}{4.20cm} 
  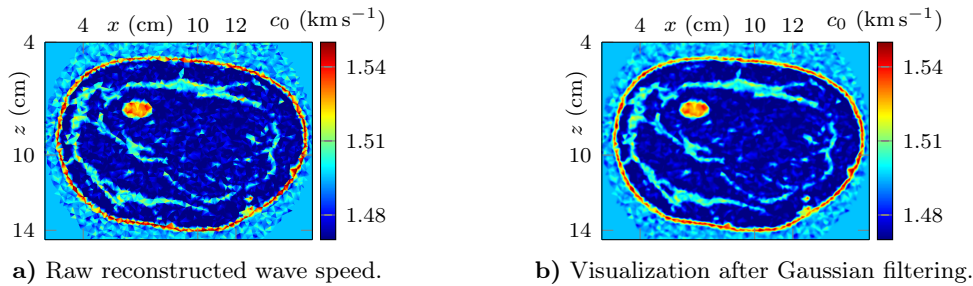
\begin{figure}[ht!]\centering
  \renewcommand{\modelfile} {cp-400kHz_data-kf_fwi-kf_raw}
  \subfloat[][Raw reconstructed wave speed.]
             {\makebox[.45\linewidth]{\begin{tikzpicture}

  \pgfmathsetmacro{\xmin}{0}
  \pgfmathsetmacro{\xmax}{18}
  \pgfmathsetmacro{\zmin}{0}
  \pgfmathsetmacro{\zmax}{18}
  \pgfmathsetmacro{\cmin}{1.47}
  \pgfmathsetmacro{\cmax}{1.55}
  
\begin{axis}[%
width=\modelwidth, height=\modelheight,
axis on top, separate axis lines,
xmin=\xminloc, xmax=\xmaxloc, xlabel={$x$ (\si{\cm})},
ymin=\zminloc, ymax=\zmaxloc, ylabel={$z$ (\si{\cm})}, y dir=reverse,
xtick={4,10,12},
ytick={4,10,14},
xticklabel pos=right, xlabel near ticks,
x label style={xshift=-0.5cm, yshift=-0.50cm}, 
y label style={xshift= 0.5cm, yshift=-0.50cm},
colormap/jet,colorbar,
colorbar style={title={{\scriptsize{$c_0$ (\si{\km\per\second})}}},
title style={yshift=-2mm, xshift=0mm},
width=2mm,xshift=-2mm,
ytick={1.48,1.51,1.54},
},point meta min=\cmin,point meta max=\cmax,
  label style={font=\scriptsize},
  tick label style={font=\scriptsize},
  legend style={font=\scriptsize\selectfont},
]
\addplot [forget plot] graphics [xmin=\xmin,xmax=\xmax,ymin=\zmin,ymax=\zmax] {{\modelfile}.png};
\end{axis}
\end{tikzpicture}%
}}
  \renewcommand{\modelfile} {cp-400kHz_data-kf_fwi-kf_gauss2}
  \subfloat[][Visualization after Gaussian filtering.]
             {\makebox[.45\linewidth]{\begin{tikzpicture}

  \pgfmathsetmacro{\xmin}{0}
  \pgfmathsetmacro{\xmax}{18}
  \pgfmathsetmacro{\zmin}{0}
  \pgfmathsetmacro{\zmax}{18}
  \pgfmathsetmacro{\cmin}{1.47}
  \pgfmathsetmacro{\cmax}{1.55}
  
\begin{axis}[%
width=\modelwidth, height=\modelheight,
axis on top, separate axis lines,
xmin=\xminloc, xmax=\xmaxloc, xlabel={$x$ (\si{\cm})},
ymin=\zminloc, ymax=\zmaxloc, ylabel={$z$ (\si{\cm})}, y dir=reverse,
xtick={4,10,12},
ytick={4,10,14},
xticklabel pos=right, xlabel near ticks,
x label style={xshift=-0.5cm, yshift=-0.50cm}, 
y label style={xshift= 0.5cm, yshift=-0.50cm},
colormap/jet,colorbar,
colorbar style={title={{\scriptsize{$c_0$ (\si{\km\per\second})}}},
title style={yshift=-2mm, xshift=0mm},
width=2mm,xshift=-2mm,
ytick={1.48,1.51,1.54},
},point meta min=\cmin,point meta max=\cmax,
  label style={font=\scriptsize},
  tick label style={font=\scriptsize},
  legend style={font=\scriptsize\selectfont},
]
\addplot [forget plot] graphics [xmin=\xmin,xmax=\xmax,ymin=\zmin,ymax=\zmax] {{\modelfile}.png};
\end{axis}
\end{tikzpicture}%
}}
  \caption{Illustration of the Gaussian filtering of the reconstruction 
           to a-posteriori improve the visualization. 
           We use a standard deviation for the Gaussian $\sigma=2$. 
           This reconstruction corresponds to the upper left image 
           of \cref{fig:fwi2d:abosrbing:all}.}
  \label{fig:fwi2d:visualization}
  \end{figure}
\end{remark}

\subsubsection{Multi-parameter reconstruction}

We have inverted with respect to model parameter 
$\kappa_0$ in the previous experiments, and now
investigate the reconstruction of two model parameters simultaneously.
We refer to the choices of the parameters to invert 
as the \emph{parametrization}, which is shown to 
have a strong influence the reconstruction, cf. \cite{Brossier2011,Kohn2012,Faucher2017}.
Namely, the physical properties characterizing 
the medium can be expressed with different physical
parameters, for instance, omitting the attenuation 
property, we have 
the bulk modulus $\kappa_0$,
the density $\rho$,
the wave speed $c_0 = \sqrt{\kappa_0/\rho}$,
or the impedance $I_0=\sqrt{\kappa_0 \, \rho}$.
For inversion, here we select two parameters among 
the aforementioned four, that will be reconstructed 
simultaneously, following \cref{algo:FWI}.
It amounts to twelve possibilities but, for the sake of 
conciseness, we consider four cases: \newline
\makebox[20em]{-- Inversion with respect to $(\kappa_0,\,\rho)$,}
\makebox[20em]{-- Inversion with respect to $(I_0,\,\rho)$,} \newline
\noindent \makebox[20.15em]{-- Inversion with respect to $(c_0,\,\rho)$,}
\makebox[20em]{-- Inversion with respect to $(I_0,\,c_0)$.} \newline
In particular, the last choice is motivated by \cite{Virieux2009}.
To evaluate the gradient of the misfit function with respect
to the selected model parameter, we start by computing the derivatives
with respect to $\kappa_0$ and $\rho$ 
(which are the main unknowns appearing in the wave equation~\cref{eq:euler_main}),
and then use the chain rule formula, cf.~\cite[Sections~5.5 and 5.6]{Faucher2017}.
In \cref{fig:fwi2d:absorbing:rho}, we compare the reconstructions: 
For all cases, we picture the wave speed $c_0$ and density $\rho$ that are 
reassembled from the reconstructed model parameters.
In addition, we maintain the attenuation model uncertainty, with data 
generated with the Kolsky--Futterman attenuation model, while inversion 
is carried out with the Kelvin--Voigt model.

\setlength{\modelwidth} {4.25cm}
\setlength{\modelheight}{3.60cm}
\graphicspath{{figures/tumor_18x18/fwi_absorbing_no-noise_all/}}
\begin{figure}[ht!]\centering
  \renewcommand{\modelfileCp} {paramK_cp_data-kf_fwi-kv_600kHz}
  \renewcommand{\modelfileRho}{paramK_rho_data-kf_fwi-kv_600kHz}
  \subfloat[][Wave speed and density obtained after
              inverting with respect to $(\kappa_0, \rho)$.]
              {\makebox[0.48\linewidth][c]{
\begin{tikzpicture}

  \pgfmathsetmacro{\xmin}{0}
  \pgfmathsetmacro{\xmax}{18}
  \pgfmathsetmacro{\zmin}{0}
  \pgfmathsetmacro{\zmax}{18}
  \pgfmathsetmacro{\cmin}{1.47}
  \pgfmathsetmacro{\cmax}{1.55}
  \pgfmathsetmacro{\cminr}{0.65}
  \pgfmathsetmacro{\cmaxr}{1.20}

\begin{groupplot}[group style={group size=2 by 1,
                  horizontal sep=11mm,
                  group name=mygroupplot},
                  enlargelimits=false
                  enlarge y limits=true,
                  enlarge x limits=false,
                  legend columns=1,
                  label style     ={font=\scriptsize},
                  tick label style={font=\scriptsize},
                  legend style    ={font=\scriptsize\selectfont},
                  enlargelimits=false,
                  height=\modelheight,width=\modelwidth]
\nextgroupplot[height=\modelheight,width=\modelwidth,
               axis on top, separate axis lines,
               xmin=\xminloc, xmax=\xmaxloc, 
               ymin=\zminloc, ymax=\zmaxloc, 
               y dir=reverse,        
               xtick={5,10},
               ytick={5,10},     
               xticklabel pos=right, xlabel near ticks,
               x label style={xshift=-0.9cm, yshift=-0.50cm}, 
               y label style={xshift= 0.7cm, yshift=-0.50cm},
               colormap/jet,colorbar,
               colorbar style={title={{\scriptsize{$c_0$ (\si{\km\per\second})}}},
               title style={yshift=-2mm, xshift=0mm},
               width=2mm,xshift=-2mm,
               ytick={1.48,1.51,1.54},
               },point meta min=\cmin,point meta max=\cmax,
               label style     ={font=\scriptsize},
               tick label style={font=\scriptsize},
               legend style    ={font=\scriptsize\selectfont},
               ]
    \addplot [forget plot] graphics [xmin=\xmin,xmax=\xmax,ymin=\zmin,ymax=\zmax] {{\modelfileCp}.png};


\nextgroupplot[width=\modelwidth, height=\modelheight,
               axis on top, separate axis lines,
               xmin=\xminloc, xmax=\xmaxloc, 
               ymin=\zminloc, ymax=\zmaxloc, 
               y dir=reverse,        
               xtick={5,10},
               ytick={5,10}, 
               xticklabel pos=right, xlabel near ticks,
               ytick={},yticklabels={,,},ylabel={},
               colormap/jet,colorbar,
               colorbar style={title={{\scriptsize{$\rho$ (\num{e3}$\,$\si{\kg\per\meter\cubed})}}},
               title style={yshift=-2mm, xshift=0mm},
               width=2mm,xshift=-2mm,
               ytick={0.7,0.9,1.10},
               },point meta min=\cminr,point meta max=\cmaxr,
               label style     ={font=\scriptsize},
               tick label style={font=\scriptsize},
               legend style    ={font=\scriptsize\selectfont},
               ]
    \addplot [forget plot] graphics [xmin=\xmin,xmax=\xmax,ymin=\zmin,ymax=\zmax] {{\modelfileRho}.png};

\end{groupplot}
                           
\end{tikzpicture}}}  \hfill
  \renewcommand{\modelfileCp} {paramI_cp_data-kf_fwi-kv_600kHz}
  \renewcommand{\modelfileRho}{paramI_rho_data-kf_fwi-kv_600kHz}
  \subfloat[][Wave speed and density obtained after
              inverting with respect to $(I_0, \rho)$.]
              {\makebox[0.48\linewidth][c]{
\begin{tikzpicture}

  \pgfmathsetmacro{\xmin}{0}
  \pgfmathsetmacro{\xmax}{18}
  \pgfmathsetmacro{\zmin}{0}
  \pgfmathsetmacro{\zmax}{18}
  \pgfmathsetmacro{\cmin}{1.47}
  \pgfmathsetmacro{\cmax}{1.55}
  \pgfmathsetmacro{\cminr}{0.65}
  \pgfmathsetmacro{\cmaxr}{1.20}

\begin{groupplot}[group style={group size=2 by 1,
                  horizontal sep=11mm,
                  group name=mygroupplot},
                  enlargelimits=false
                  enlarge y limits=true,
                  enlarge x limits=false,
                  legend columns=1,
                  label style     ={font=\scriptsize},
                  tick label style={font=\scriptsize},
                  legend style    ={font=\scriptsize\selectfont},
                  enlargelimits=false,
                  height=\modelheight,width=\modelwidth]
\nextgroupplot[height=\modelheight,width=\modelwidth,
               axis on top, separate axis lines,
               xmin=\xminloc, xmax=\xmaxloc, 
               ymin=\zminloc, ymax=\zmaxloc, 
               y dir=reverse,        
               xtick={5,10},
               ytick={5,10},     
               xticklabel pos=right, xlabel near ticks,
               x label style={xshift=-0.9cm, yshift=-0.50cm}, 
               y label style={xshift= 0.7cm, yshift=-0.50cm},
               colormap/jet,colorbar,
               colorbar style={title={{\scriptsize{$c_0$ (\si{\km\per\second})}}},
               title style={yshift=-2mm, xshift=0mm},
               width=2mm,xshift=-2mm,
               ytick={1.48,1.51,1.54},
               },point meta min=\cmin,point meta max=\cmax,
               label style     ={font=\scriptsize},
               tick label style={font=\scriptsize},
               legend style    ={font=\scriptsize\selectfont},
               ]
    \addplot [forget plot] graphics [xmin=\xmin,xmax=\xmax,ymin=\zmin,ymax=\zmax] {{\modelfileCp}.png};


\nextgroupplot[width=\modelwidth, height=\modelheight,
               axis on top, separate axis lines,
               xmin=\xminloc, xmax=\xmaxloc, 
               ymin=\zminloc, ymax=\zmaxloc, 
               y dir=reverse,        
               xtick={5,10},
               ytick={5,10}, 
               xticklabel pos=right, xlabel near ticks,
               ytick={},yticklabels={,,},ylabel={},
               colormap/jet,colorbar,
               colorbar style={title={{\scriptsize{$\rho$ (\num{e3}$\,$\si{\kg\per\meter\cubed})}}},
               title style={yshift=-2mm, xshift=0mm},
               width=2mm,xshift=-2mm,
               ytick={0.7,0.9,1.10},
               },point meta min=\cminr,point meta max=\cmaxr,
               label style     ={font=\scriptsize},
               tick label style={font=\scriptsize},
               legend style    ={font=\scriptsize\selectfont},
               ]
    \addplot [forget plot] graphics [xmin=\xmin,xmax=\xmax,ymin=\zmin,ymax=\zmax] {{\modelfileRho}.png};

\end{groupplot}
                           
\end{tikzpicture}}} \\

  \renewcommand{\modelfileCp} {paramV_cp_data-kf_fwi-kv_600kHz}
  \renewcommand{\modelfileRho}{paramV_rho_data-kf_fwi-kv_600kHz}
  \subfloat[][Wave speed and density obtained after
              inverting with respect to $(c_0, \rho)$.]
              {\makebox[0.48\linewidth][c]{
\begin{tikzpicture}

  \pgfmathsetmacro{\xmin}{0}
  \pgfmathsetmacro{\xmax}{18}
  \pgfmathsetmacro{\zmin}{0}
  \pgfmathsetmacro{\zmax}{18}
  \pgfmathsetmacro{\cmin}{1.47}
  \pgfmathsetmacro{\cmax}{1.55}
  \pgfmathsetmacro{\cminr}{0.65}
  \pgfmathsetmacro{\cmaxr}{1.20}

\begin{groupplot}[group style={group size=2 by 1,
                  horizontal sep=11mm,
                  group name=mygroupplot},
                  enlargelimits=false
                  enlarge y limits=true,
                  enlarge x limits=false,
                  legend columns=1,
                  label style     ={font=\scriptsize},
                  tick label style={font=\scriptsize},
                  legend style    ={font=\scriptsize\selectfont},
                  enlargelimits=false,
                  height=\modelheight,width=\modelwidth]
\nextgroupplot[height=\modelheight,width=\modelwidth,
               axis on top, separate axis lines,
               xmin=\xminloc, xmax=\xmaxloc, 
               ymin=\zminloc, ymax=\zmaxloc, 
               y dir=reverse,        
               xtick={5,10},
               ytick={5,10},     
               xticklabel pos=right, xlabel near ticks,
               x label style={xshift=-0.9cm, yshift=-0.50cm}, 
               y label style={xshift= 0.7cm, yshift=-0.50cm},
               colormap/jet,colorbar,
               colorbar style={title={{\scriptsize{$c_0$ (\si{\km\per\second})}}},
               title style={yshift=-2mm, xshift=0mm},
               width=2mm,xshift=-2mm,
               ytick={1.48,1.51,1.54},
               },point meta min=\cmin,point meta max=\cmax,
               label style     ={font=\scriptsize},
               tick label style={font=\scriptsize},
               legend style    ={font=\scriptsize\selectfont},
               ]
    \addplot [forget plot] graphics [xmin=\xmin,xmax=\xmax,ymin=\zmin,ymax=\zmax] {{\modelfileCp}.png};


\nextgroupplot[width=\modelwidth, height=\modelheight,
               axis on top, separate axis lines,
               xmin=\xminloc, xmax=\xmaxloc, 
               ymin=\zminloc, ymax=\zmaxloc, 
               y dir=reverse,        
               xtick={5,10},
               ytick={5,10}, 
               xticklabel pos=right, xlabel near ticks,
               ytick={},yticklabels={,,},ylabel={},
               colormap/jet,colorbar,
               colorbar style={title={{\scriptsize{$\rho$ (\num{e3}$\,$\si{\kg\per\meter\cubed})}}},
               title style={yshift=-2mm, xshift=0mm},
               width=2mm,xshift=-2mm,
               ytick={0.7,0.9,1.10},
               },point meta min=\cminr,point meta max=\cmaxr,
               label style     ={font=\scriptsize},
               tick label style={font=\scriptsize},
               legend style    ={font=\scriptsize\selectfont},
               ]
    \addplot [forget plot] graphics [xmin=\xmin,xmax=\xmax,ymin=\zmin,ymax=\zmax] {{\modelfileRho}.png};

\end{groupplot}
                           
\end{tikzpicture}}} \hfill
  \renewcommand{\modelfileCp} {paramIV_cp_data-kf_fwi-kv_600kHz}
  \renewcommand{\modelfileRho}{paramIV_rho_data-kf_fwi-kv_600kHz}
  \subfloat[][Wave speed and density obtained after
              inverting with respect to $(I_0, c_0)$.]
              {\makebox[0.48\linewidth][c]{
\begin{tikzpicture}

  \pgfmathsetmacro{\xmin}{0}
  \pgfmathsetmacro{\xmax}{18}
  \pgfmathsetmacro{\zmin}{0}
  \pgfmathsetmacro{\zmax}{18}
  \pgfmathsetmacro{\cmin}{1.47}
  \pgfmathsetmacro{\cmax}{1.55}
  \pgfmathsetmacro{\cminr}{0.65}
  \pgfmathsetmacro{\cmaxr}{1.20}

\begin{groupplot}[group style={group size=2 by 1,
                  horizontal sep=11mm,
                  group name=mygroupplot},
                  enlargelimits=false
                  enlarge y limits=true,
                  enlarge x limits=false,
                  legend columns=1,
                  label style     ={font=\scriptsize},
                  tick label style={font=\scriptsize},
                  legend style    ={font=\scriptsize\selectfont},
                  enlargelimits=false,
                  height=\modelheight,width=\modelwidth]
\nextgroupplot[height=\modelheight,width=\modelwidth,
               axis on top, separate axis lines,
               xmin=\xminloc, xmax=\xmaxloc, 
               ymin=\zminloc, ymax=\zmaxloc, 
               y dir=reverse,        
               xtick={5,10},
               ytick={5,10},     
               xticklabel pos=right, xlabel near ticks,
               x label style={xshift=-0.9cm, yshift=-0.50cm}, 
               y label style={xshift= 0.7cm, yshift=-0.50cm},
               colormap/jet,colorbar,
               colorbar style={title={{\scriptsize{$c_0$ (\si{\km\per\second})}}},
               title style={yshift=-2mm, xshift=0mm},
               width=2mm,xshift=-2mm,
               ytick={1.48,1.51,1.54},
               },point meta min=\cmin,point meta max=\cmax,
               label style     ={font=\scriptsize},
               tick label style={font=\scriptsize},
               legend style    ={font=\scriptsize\selectfont},
               ]
    \addplot [forget plot] graphics [xmin=\xmin,xmax=\xmax,ymin=\zmin,ymax=\zmax] {{\modelfileCp}.png};


\nextgroupplot[width=\modelwidth, height=\modelheight,
               axis on top, separate axis lines,
               xmin=\xminloc, xmax=\xmaxloc, 
               ymin=\zminloc, ymax=\zmaxloc, 
               y dir=reverse,        
               xtick={5,10},
               ytick={5,10}, 
               xticklabel pos=right, xlabel near ticks,
               ytick={},yticklabels={,,},ylabel={},
               colormap/jet,colorbar,
               colorbar style={title={{\scriptsize{$\rho$ (\num{e3}$\,$\si{\kg\per\meter\cubed})}}},
               title style={yshift=-2mm, xshift=0mm},
               width=2mm,xshift=-2mm,
               ytick={0.7,0.9,1.10},
               },point meta min=\cminr,point meta max=\cmaxr,
               label style     ={font=\scriptsize},
               tick label style={font=\scriptsize},
               legend style    ={font=\scriptsize\selectfont},
               ]
    \addplot [forget plot] graphics [xmin=\xmin,xmax=\xmax,ymin=\zmin,ymax=\zmax] {{\modelfileRho}.png};

\end{groupplot}
                           
\end{tikzpicture}}}
  \caption{Comparison of the reconstruction depending on 
           the parametrization of the inverse problem. 
           The reconstruction follows \cref{algo:FWI} 
           with sequential frequency progression 
           $\omegaR/(2\pi)=\{\num{200}, \, \num{300}, \, \num{400}, \, \num{600}\}$ \si{\kilo\Hz} and
           fixed $\omegaI=0$ with \num{30} iterations per frequency, 
           starting from homogeneous backgrounds (\cref{table:interval}).
           The domain units are given in \si{\centi\meter}.
           }
  \label{fig:fwi2d:absorbing:rho}
\end{figure}
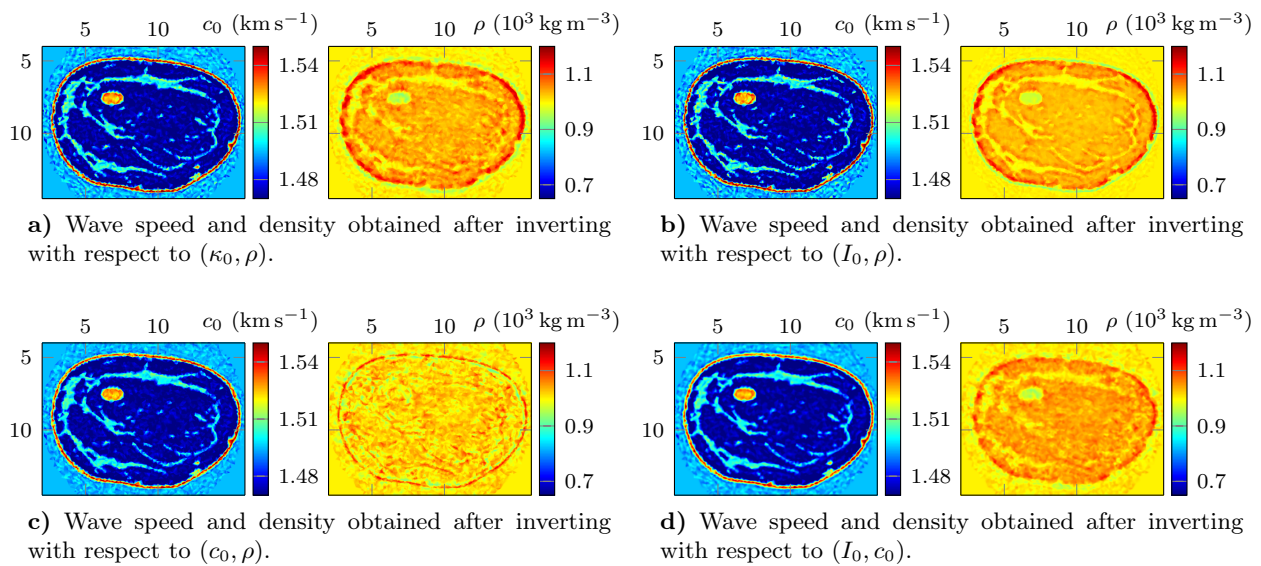

The assembled wave speed is similarly accurate with 
all choices of parametrizations, with the different 
tissue features correctly retrieved and the ellipse-shaped 
defect resulting in a contrast of wave speed is visually 
clearly identified.
However, we observe major differences in the 
reconstruction of the density: 
The parametrization $(\kappa_0, \rho)$,
$(I_0, c_0)$ and $(I_0, \rho)$ are able 
to identify the main layers and $(I_0, \rho)$ 
provides the best resolution, with smoother contrasts.
Nonetheless, the values of the density are incorrect
in all cases, and remain near to the initial value.
Here, the parametrization $(c_0, \rho)$ is the worst 
regarding the density reconstruction, with a noisy 
background where only the skin delimitation barely 
appears.

\begin{remark}
  We note that the difficulty of recovering the density 
  is similar in seismic imaging, \cite{Virieux2009,Faucher2017}.
  This is explained as density perturbations mostly modify
  the amplitude of the signals and, as indicated above, the
  FWI focuses on the differences in the signals phase-shift. 
  As an alternative, one could carry the FWI reconstruction 
  of one main parameter, such as the bulk modulus, wave speed 
  or impedance, and then infer the density via analytic formula 
  inherited from physics, such as Gardner's relation, \cite{Gardner1974}.
  Another alternative is to select a different misfit function.
\end{remark}

\subsection{Reconstruction with wall boundary conditions and complex frequencies}
\label{subsection:2d:wall-bc}

In the previous experiments, we have seen that 
the wave speed is accurately reconstructed, 
even in the case of an attenuation model error.
We have assumed that waves can escape the medium 
without incurring reflections, by using in the numerical simulations 
absorbing boundary conditions on all sides. 
We now instead consider wall boundary conditions
on the domain boundary $\Gamma$ with \cref{eq:bc:wall}, 
hence leading to multiple wave reflections from each side.
The background coefficients also 
have low attenuation (see \cref{table:interval}), 
further supporting these multiple reflections.
We shall see that the modification of the boundary condition
drastically alters the performance of the reconstruction, 
and that complex frequencies can be used to overcome the 
difficulties.
For the sake of conciseness, we only consider the 
data-set using the Kolsky--Futterman attenuation model, 
while the reconstruction procedure is carried out with 
the Kelvin--Voigt model.

\subsubsection{Reconstruction with ordinary frequencies ($\omegaR > 0$, $\omegaI=0$)} 

We keep the setup of the previous experiment where we 
only use Fourier frequencies with the progression 
$\omegaR/(2\pi) = \{\num{200}, \, \num{300}, \, \num{400}, \num{600} \}$ \si{\kilo\Hz},
with a fixed $\omegaI = 0$.
We invert with respect to the bulk modulus only, 
and picture the reconstruction in \cref{fig:fwi2d:wall:compare-Q}, 
where we compare different choices of fixed quality
factor for the reconstruction, that is, different
choices of initial values for the attenuation parameters.

\setlength{\modelwidth} {5.10cm}
\setlength{\modelheight}{4.20cm} 
\graphicspath{{figures/tumor_18x18/fwi_wall_no-noise/}}
\begin{figure}[ht!]\centering
  \renewcommand{\modelfile} {Q800_cp_data-kf_fwi-kv_600kHz_gauss2}
  \subfloat[][Reconstruction with $Q=800$ (value of $Q$ in the background).]
             {\makebox[.30\linewidth]{\begin{tikzpicture}

  \pgfmathsetmacro{\xmin}{0}
  \pgfmathsetmacro{\xmax}{18}
  \pgfmathsetmacro{\zmin}{0}
  \pgfmathsetmacro{\zmax}{18}
  \pgfmathsetmacro{\cmin}{1.47}
  \pgfmathsetmacro{\cmax}{1.55}
  
\begin{axis}[%
width=\modelwidth, height=\modelheight,
axis on top, separate axis lines,
xmin=\xminloc, xmax=\xmaxloc, xlabel={$x$ (\si{\cm})},
ymin=\zminloc, ymax=\zmaxloc, ylabel={$z$ (\si{\cm})}, y dir=reverse,
xtick={4,10,12},
ytick={4,10,14},
xticklabel pos=right, xlabel near ticks,
x label style={xshift=-0.5cm, yshift=-0.50cm}, 
y label style={xshift= 0.5cm, yshift=-0.50cm},
colormap/jet,colorbar,
colorbar style={title={{\scriptsize{$c_0$ (\si{\km\per\second})}}},
title style={yshift=-2mm, xshift=0mm},
width=2mm,xshift=-2mm,
ytick={1.48,1.51,1.54},
},point meta min=\cmin,point meta max=\cmax,
  label style={font=\scriptsize},
  tick label style={font=\scriptsize},
  legend style={font=\scriptsize\selectfont},
]
\addplot [forget plot] graphics [xmin=\xmin,xmax=\xmax,ymin=\zmin,ymax=\zmax] {{\modelfile}.png};
\end{axis}
\end{tikzpicture}%
}
              \label{fig:fwi2d:wall:compare-Q_a}} \hfill
  \renewcommand{\modelfile} {Q200_cp_data-kf_fwi-kv_600kHz_gauss2}
  \subfloat[][Reconstruction with $Q=200$.]
             {\makebox[.30\linewidth]{\begin{tikzpicture}

  \pgfmathsetmacro{\xmin}{0}
  \pgfmathsetmacro{\xmax}{18}
  \pgfmathsetmacro{\zmin}{0}
  \pgfmathsetmacro{\zmax}{18}
  \pgfmathsetmacro{\cmin}{1.47}
  \pgfmathsetmacro{\cmax}{1.55}
  
\begin{axis}[%
width=\modelwidth, height=\modelheight,
axis on top, separate axis lines,
xmin=\xminloc, xmax=\xmaxloc, xlabel={$x$ (\si{\cm})},
ymin=\zminloc, ymax=\zmaxloc, ylabel={$z$ (\si{\cm})}, y dir=reverse,
xtick={4,10,12},
ytick={4,10,14},
xticklabel pos=right, xlabel near ticks,
x label style={xshift=-0.5cm, yshift=-0.50cm}, 
y label style={xshift= 0.5cm, yshift=-0.50cm},
colormap/jet,colorbar,
colorbar style={title={{\scriptsize{$c_0$ (\si{\km\per\second})}}},
title style={yshift=-2mm, xshift=0mm},
width=2mm,xshift=-2mm,
ytick={1.48,1.51,1.54},
},point meta min=\cmin,point meta max=\cmax,
  label style={font=\scriptsize},
  tick label style={font=\scriptsize},
  legend style={font=\scriptsize\selectfont},
]
\addplot [forget plot] graphics [xmin=\xmin,xmax=\xmax,ymin=\zmin,ymax=\zmax] {{\modelfile}.png};
\end{axis}
\end{tikzpicture}%
}
              \label{fig:fwi2d:wall:compare-Q_b}} \hfill
  \renewcommand{\modelfile} {Q100_cp_data-kf_fwi-kv_600kHz_gauss2}
  \subfloat[][Reconstruction with $Q=100$.]
             {\makebox[.30\linewidth]{\begin{tikzpicture}

  \pgfmathsetmacro{\xmin}{0}
  \pgfmathsetmacro{\xmax}{18}
  \pgfmathsetmacro{\zmin}{0}
  \pgfmathsetmacro{\zmax}{18}
  \pgfmathsetmacro{\cmin}{1.47}
  \pgfmathsetmacro{\cmax}{1.55}
  
\begin{axis}[%
width=\modelwidth, height=\modelheight,
axis on top, separate axis lines,
xmin=\xminloc, xmax=\xmaxloc, xlabel={$x$ (\si{\cm})},
ymin=\zminloc, ymax=\zmaxloc, ylabel={$z$ (\si{\cm})}, y dir=reverse,
xtick={4,10,12},
ytick={4,10,14},
xticklabel pos=right, xlabel near ticks,
x label style={xshift=-0.5cm, yshift=-0.50cm}, 
y label style={xshift= 0.5cm, yshift=-0.50cm},
colormap/jet,colorbar,
colorbar style={title={{\scriptsize{$c_0$ (\si{\km\per\second})}}},
title style={yshift=-2mm, xshift=0mm},
width=2mm,xshift=-2mm,
ytick={1.48,1.51,1.54},
},point meta min=\cmin,point meta max=\cmax,
  label style={font=\scriptsize},
  tick label style={font=\scriptsize},
  legend style={font=\scriptsize\selectfont},
]
\addplot [forget plot] graphics [xmin=\xmin,xmax=\xmax,ymin=\zmin,ymax=\zmax] {{\modelfile}.png};
\end{axis}
\end{tikzpicture}%
}
              \label{fig:fwi2d:wall:compare-Q_c}} 
  \caption{Reconstruction of the wave speed 
           $c_0 = \sqrt{\kappa_0/\rho}$
           in the case of
           wall boundary conditions with different choices
           of (fixed) attenuation parameters.
           The reconstruction follows \cref{algo:FWI} with 
           sequential frequency progression 
           $\omegaR/(2\pi)=\{200, \,300, \, 400, \, 600\}$ \si{\kilo\Hz} and
           fixed $\omegaI=0$ with \num{30} iterations per frequency.
           The reconstruction starts from homogeneous wave speed and 
           density 
           (\cref{table:interval}),
           and we compare different choices of constant attenuation 
           parameters to observe the effect of the initial 
           guess for the attenuation model onto the reconstruction. 
           }
  \label{fig:fwi2d:wall:compare-Q}
\end{figure}
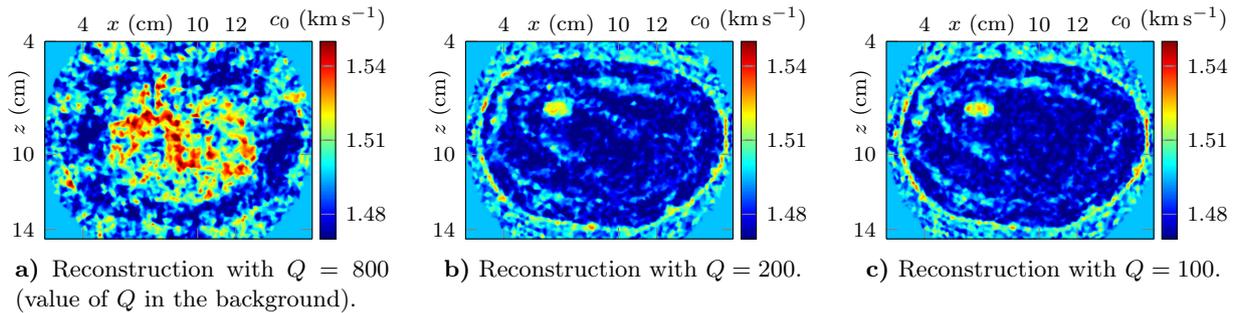

We see that the consideration of wall boundary on the sides
of the medium completely changes the behaviour and prevent
any accurate reconstruction when using the 
constant value of background attenuation
as initial guess, \cref{fig:fwi2d:wall:compare-Q_a}, 
while it is the ``natural'' choice of background value.
It appears that the multiple reflections from the boundaries 
totally annihilate the possibility of reconstruction, with 
information from the sample merged within the wave reflections.
Interestingly, when the initial guess for 
the quality factor is much lower than expected 
(i.e., the initial guess is too attenuating compared
to the real medium), 
the FWI reconstruction is slightly better, with
the tissue features appearing, see 
\cref{fig:fwi2d:wall:compare-Q_b,fig:fwi2d:wall:compare-Q_c}.
However, the reconstruction remains with a poor resolution
compared to the previous experiment assuming free-space propagation.

This experiment highlights the difficulty of having strong
multiple reflections from the medium boundary, drastically
reducing the performance of the reconstruction procedure. 
It would be necessary to treat those multiples in the measurement 
data, \cite{Yilmaz2001}, to clear
the reflections coming from the domain boundaries,
but this task is not trivial in practice.
As an alternative, one could 
perform the experiments in a medium with 
strong attenuation, to reduce or prevent
the multiple reflections.
In the next section, we show how to resolve the issue of 
the multiple boundary reflections, by encoding an artificial 
damping via the use of complex frequencies.

\subsubsection{Reconstruction using complex frequencies ($\omegaR > 0$, $\omegaI \geq 0$)}

We now investigate the use of complex frequencies for 
the case with wall boundary conditions, where we keep the 
same set of Fourier frequency, enriched with a Laplace
component. In this case, it corresponds to applying 
a complex Fourier transform to the original time-domain 
observed data, as illustrated in \cref{fig:time-domain-acquisition-cx-freq}.

For the reconstruction, we use the following sets:
\begin{itemize}\setlength{\itemsep}{-1pt}
  \item Fourier frequency $\omegaR/(2\pi)$: $\{200,\,300,\,400,\,600\}$ \si{\kilo\Hz},
  \item Laplace damping   $\omegaI$:        $\{\num{20e3},\,\num{15e3},\,\num{10e3}\}$ \si{\per\second}.
\end{itemize}
Following the analysis of the inverse problem for complex frequencies
carried out in \cite{Faucher2020basins}, the progression of content 
first varies the Laplace damping from high to low for a fixed 
Fourier frequency, which then varies from low high. 
Namely, the inverted complex frequencies follow the order:
$(200\si{\kilo\Hz}, \num{20e3})$, $(200\si{\kilo\Hz}, \num{15e3})$, $(200\si{\kilo\Hz}, \num{10e3})$, 
$(300\si{\kilo\Hz}, \num{20e3})$, $(300\si{\kilo\Hz}, \num{15e3})$, etc.
The reconstructed wave speed $c_0$ using the complex frequency 
set is pictured in \cref{fig:fwi2d:wall:CX}.

\setlength{\modelwidth} {7.00cm}
\setlength{\modelheight}{6.00cm} 
\graphicspath{{figures/tumor_18x18/fwi_wall_no-noise/}}
\begin{figure}[ht!]\centering
  \renewcommand{\modelfile} {CXfreq_cp_data-kf_fwi-kv_600kHz-15e3_gauss2}
  \begin{tikzpicture}

  \pgfmathsetmacro{\xmin}{0}
  \pgfmathsetmacro{\xmax}{18}
  \pgfmathsetmacro{\zmin}{0}
  \pgfmathsetmacro{\zmax}{18}
  \pgfmathsetmacro{\cmin}{1.47}
  \pgfmathsetmacro{\cmax}{1.55}
  
\begin{axis}[%
width=\modelwidth, height=\modelheight,
axis on top, separate axis lines,
xmin=\xminloc, xmax=\xmaxloc, xlabel={$x$ (\si{\cm})},
ymin=\zminloc, ymax=\zmaxloc, ylabel={$z$ (\si{\cm})}, y dir=reverse,
xtick={4,10,12},
ytick={4,10,14},
xticklabel pos=right, xlabel near ticks,
x label style={xshift=-0.5cm, yshift=-0.50cm}, 
y label style={xshift= 0.5cm, yshift=-0.50cm},
colormap/jet,colorbar,
colorbar style={title={{\scriptsize{$c_0$ (\si{\km\per\second})}}},
title style={yshift=-2mm, xshift=0mm},
width=2mm,xshift=-2mm,
ytick={1.48,1.51,1.54},
},point meta min=\cmin,point meta max=\cmax,
  label style={font=\scriptsize},
  tick label style={font=\scriptsize},
  legend style={font=\scriptsize\selectfont},
]
\addplot [forget plot] graphics [xmin=\xmin,xmax=\xmax,ymin=\zmin,ymax=\zmax] {{\modelfile}.png};
\end{axis}
\end{tikzpicture}%
  \caption{Reconstruction of medium with wall boundary condition using complex frequencies.
           The starting models correspond to 
           the constant background values of
           wave speed, density and quality factor
           of \cref{table:interval}.
           (i.e., using fixed attenuation parameters chosen such that
           $Q=800$ at \num{300} \si{\kilo\Hz}). 
           The corresponding case without using the complex 
           frequencies is pictured \cref{fig:fwi2d:wall:compare-Q_a}.}
  \label{fig:fwi2d:wall:CX}
\end{figure}
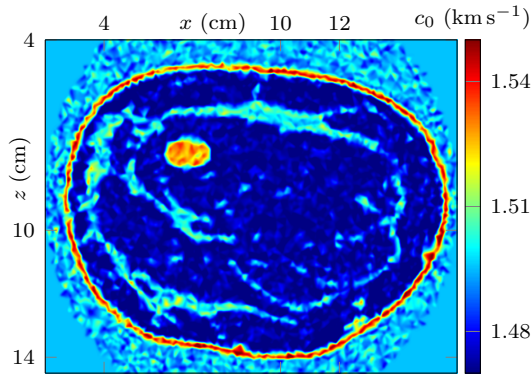
The reconstruction using complex frequencies is accurate,
with all the layers and the contrasting tumor appearing 
with the same high resolution as in the case of absorbing
boundary conditions.
Therefore the artificial damping has allowed to overcome
the multiple boundary reflections.
Note that this procedure can be applied for any data measured
in the time-domain, and only necessitates to apply the complex
Fourier transform as illustrated in \cref{fig:time-domain-acquisition-cx-freq}.
The performance of the use of the Laplace-Fourier transform 
with several complex frequencies can further be explained as
the transform introduces artificial damping with the complex frequency, 
hence enhancing the first arrivals (\cite{Kamei2013}) by reducing the 
multiple reflections.

\section{Numerical experiments in three dimensions}
\label{section:fwi-3d}

\subsection{Experimental setup}

We set up a three-dimensional experiment with
a breast model extracted from the OA-Breast 
Phantom data-set, see \cref{footnote:OA-breast},
and further include a contrasting 3D object 
embedded in the medium. 
In \cref{fig:3d:true} are displayed the wave 
speed, density and quality factor, which are
encompassed in a rectangular domain of size
\num{10} $\times$ \num{10} $\times$ \num{16} 
\si{\centi\meter\cubed}. 
The values of the parameters in the different 
layers of the sample follow the values of the 
previous experiment, and are prescribed in 
\cref{table:interval}.

\graphicspath{{figures/tumor_3d/main/}}
\setlength{\modelwidth} {7cm}
\setlength{\modelheight}{7cm}
\setlength{\jumpvert}   {3cm}
\begin{figure}[ht!]\centering
  \subfloat[][3D breast sample.]
           {{\raisebox{1em}{
             \includegraphics[scale=0.4]{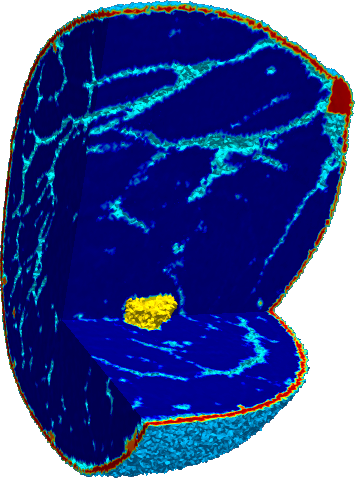}}}} \hspace*{2em}
  \renewcommand{\modelfile}{cp-true_scale1440-1610}
  \subfloat[][Wave speed model.]
           {
\begin{tikzpicture}
\pgfmathsetmacro{\xmin} {0.}
\pgfmathsetmacro{\xmax} {10.0}
\pgfmathsetmacro{\ymin} {0.}
\pgfmathsetmacro{\ymax} {16.0}
\pgfmathsetmacro{\zmin} {0.}
\pgfmathsetmacro{\zmax} {16.0}

\pgfmathsetmacro{\ycut} { 8.00}
\pgfmathsetmacro{\xcut} { 7.54}
\pgfmathsetmacro{\zcut} {10.80}
\pgfmathsetmacro{\vmin} {1.440}
\pgfmathsetmacro{\vmax} {1.610}

\matrix[column sep=-1.80em, row sep=-1.6em] {
\begin{axis}[yshift=0.7*\jumpvert,
  anchor=center,
  tick label style={font=\small},
  grid=both,minor tick num=1,
  xlabel={\scriptsize{$x$ (\si{\cm})}},
  ylabel={\scriptsize{$y$ (\si{\cm})}},
  zlabel={\scriptsize{$z$ (\si{\cm})}}, 
  ztick pos=left,
  3d box,width=1.4\modelwidth, 
  xmin=\xmin,ymin=\ymin,zmin=\zmin,xmax=\xmax,ymax=\ymax,zmax=\zmax,
  every axis x label/.style={at={(0.05, 0.10)},anchor=north},
  every axis y label/.style={at={(1.00, 0.15)},anchor=north}, 
  z label style={xshift= 0.5cm, yshift=-0.50cm},
  xtick={0,10},
  line width=.25pt,
  label style={font=\scriptsize},
  tick label style={font=\scriptsize},
  legend style={font=\scriptsize\selectfont},  
  ]
  \addplot3[fill=white] graphics[points={
            (\xcut,\ymin,\zmin)     => ( 16,215-52)
            (\xcut,\ymax,\zmin)     => (152,215-0)
            (\xcut,\ymin,\zmax)     => ( 16,215-208)
            (\xmin,\ycut,\zmax)     => (132,215-215)}
            ]{{\modelfile_3d}.png};
  
  \addplot3[dashed,black,line width=1pt] 
      coordinates {(\xmin,\ycut,\zmin) (\xmin,\ycut,\zmax)
                   (\xmax,\ycut,\zmax) (\xmax,\ycut,\zmin) (\xmin,\ycut,\zmin)};
  \addplot3[dashed,black,line width=1pt] 
      coordinates {(\xmin,\ymin,\zcut) (\xmin,\ymax,\zcut) 
                   (\xmax,\ymax,\zcut) (\xmax,\ymin,\zcut) (\xmin,\ymin,\zcut)};

  \addplot3[dashed,black,line width=1pt] 
      coordinates {(\xcut,\ymin,\zmin) (\xcut,\ymin,\zmax)
                   (\xcut,\ymax,\zmax) (\xcut,\ymax,\zmin) (\xcut,\ymin,\zmin)};

\end{axis}

&

\begin{axis}[yshift=0.1em,anchor=center,  
             width=0.51\modelwidth, height=0.64\modelheight,
             xmin=\xmin, xmax=\xmax,ymin=\zmin, ymax=\zmax,hide axis,
             tick label style={font=\small},y dir=reverse]
         \addplot [forget plot] graphics 
         [xmin=\xmin,xmax=\xmax,ymin=\zmin,ymax=\zmax] {{\modelfile-1}.png};
         \addplot[dashed,black,line width=1pt] 
       coordinates {(\xmin+.01,\zmin+.01) (\xmin+.01,\zmax-.01) 
                    (\xmax-.01,\zmax-.01) (\xmax-.01,\zmin+.01) (\xmin+.01,\zmin+.01)};
\end{axis}

\begin{axis}[yshift=\jumpvert,
             anchor=center,  width =0.64\modelwidth, 
                             height=0.64\modelheight,
             xmin=\xmin, xmax=\xmax,ymin=\zmin, ymax=\zmax,hide axis,
             tick label style={font=\small},y dir=reverse]
         \addplot [forget plot] graphics 
         [xmin=\xmin,xmax=\xmax,ymin=\zmin,ymax=\zmax] {{\modelfile-3}.png};
         \addplot[dashed,black,line width=1pt] 
       coordinates {(\xmin+.01,\zmin+.01) (\xmin+.01,\zmax-.01) 
                    (\xmax-.01,\zmax-.01) (\xmax-.01,\zmin+.01) (\xmin+.01,\zmin+.01)};
\end{axis}

\\

\begin{axis}[anchor=center, width=.7\modelwidth, height=.51\modelwidth,
             xmin=\xmin, xmax=\xmax,ymin=\ymin, ymax=\ymax,hide axis,
             tick label style={font=\small}, y dir=reverse]
            \addplot [forget plot] graphics 
                     [xmin=\xmin,xmax=\xmax,ymin=\ymin,ymax=\ymax] {{\modelfile-2}.png};
            
            \addplot[dashed,black,line width=1pt] 
                  coordinates {(\xmin+.01,\ymin+.01) (\xmin+.01,\ymax-.01) 
                               (\xmax-.01,\ymax-.01) (\xmax-.01,\ymin+.01) (\xmin+.01,\ymin+.01)};
\end{axis}

&

\begin{axis}[yshift=-1.80em,
     anchor=south,yshift=.5cm,hide axis,height=.5cm,colorbar/width=.2cm,
     colormap/jet,colorbar horizontal,
     colorbar style={separate axis lines,
                     title style={yshift=-3mm, xshift=0mm},
                     tick label style={font=\small},
                     title={\small{wave speed (\si{\km\per\second})}}},
                     point meta min=\vmin, point meta max=\vmax,
                     height=.5\modelwidth]
\end{axis}
\\};

\draw[black,-stealth,->,dashed,line width=1pt] (-0.2, 3.3) to[out=0,in=-180] ( 0.6, 3.3);
\draw[black,-stealth,->,dashed,line width=1pt] (-0.6, 0.0) to[out=0,in=-180] ( 1.0, 0.0);
\draw[black,-stealth,->,dashed,line width=1pt] (-1.75,-0.15) to[out=-90,in=90] (-1.75,-1.6);

\end{tikzpicture}}
  \renewcommand{\modelfile}{rho-true_scale700-1150}
  \subfloat[][Density model.]
           {
\begin{tikzpicture}
\pgfmathsetmacro{\xmin} {0.}
\pgfmathsetmacro{\xmax} {10.0}
\pgfmathsetmacro{\ymin} {0.}
\pgfmathsetmacro{\ymax} {16.0}
\pgfmathsetmacro{\zmin} {0.}
\pgfmathsetmacro{\zmax} {16.0}

\pgfmathsetmacro{\ycut} { 8.00}
\pgfmathsetmacro{\xcut} { 7.54}
\pgfmathsetmacro{\zcut} {10.80}
\pgfmathsetmacro{\vmin} {0.7}
\pgfmathsetmacro{\vmax} {1.15}

\matrix[column sep=-1.8em, row sep=-1.6em] {
\begin{axis}[yshift=0.7*\jumpvert,
  anchor=center,
  tick label style={font=\small},
  grid=both,minor tick num=1,
  xlabel={\scriptsize{$x$ (\si{\cm})}},
  ylabel={\scriptsize{$y$ (\si{\cm})}},
  zlabel={\scriptsize{$z$ (\si{\cm})}}, 
  ztick pos=left,
  3d box,width=1.4\modelwidth, 
  xmin=\xmin,ymin=\ymin,zmin=\zmin,xmax=\xmax,ymax=\ymax,zmax=\zmax,
  every axis x label/.style={at={(0.05, 0.10)},anchor=north},
  every axis y label/.style={at={(1.00, 0.15)},anchor=north}, 
  z label style={xshift= 0.5cm, yshift=-0.50cm},
  xtick={0,10},
  line width=.25pt,
  label style={font=\scriptsize},
  tick label style={font=\scriptsize},
  legend style={font=\scriptsize\selectfont},  
  ]
  \addplot3[fill=white] graphics[points={
            (\xcut,\ymin,\zmin)     => ( 16,215-52)
            (\xcut,\ymax,\zmin)     => (152,215-0)
            (\xcut,\ymin,\zmax)     => ( 16,215-208)
            (\xmin,\ycut,\zmax)     => (132,215-215)}
            ]{{\modelfile_3d}.png};
  
  \addplot3[dashed,black,line width=1pt] 
      coordinates {(\xmin,\ycut,\zmin) (\xmin,\ycut,\zmax)
                   (\xmax,\ycut,\zmax) (\xmax,\ycut,\zmin) (\xmin,\ycut,\zmin)};
  \addplot3[dashed,black,line width=1pt] 
      coordinates {(\xmin,\ymin,\zcut) (\xmin,\ymax,\zcut) 
                   (\xmax,\ymax,\zcut) (\xmax,\ymin,\zcut) (\xmin,\ymin,\zcut)};

  \addplot3[dashed,black,line width=1pt] 
      coordinates {(\xcut,\ymin,\zmin) (\xcut,\ymin,\zmax)
                   (\xcut,\ymax,\zmax) (\xcut,\ymax,\zmin) (\xcut,\ymin,\zmin)};

\end{axis}

&

\begin{axis}[yshift=0.1em,anchor=center,  
             width=0.51\modelwidth, height=0.64\modelheight,
             xmin=\xmin, xmax=\xmax,ymin=\zmin, ymax=\zmax,hide axis,
             tick label style={font=\small},y dir=reverse]
         \addplot [forget plot] graphics 
         [xmin=\xmin,xmax=\xmax,ymin=\zmin,ymax=\zmax] {{\modelfile-1}.png};
         \addplot[dashed,black,line width=1pt] 
       coordinates {(\xmin+.01,\zmin+.01) (\xmin+.01,\zmax-.01) 
                    (\xmax-.01,\zmax-.01) (\xmax-.01,\zmin+.01) (\xmin+.01,\zmin+.01)};
\end{axis}

\begin{axis}[yshift=\jumpvert,
             anchor=center,  width =0.64\modelwidth, 
                             height=0.64\modelheight,
             xmin=\xmin, xmax=\xmax,ymin=\zmin, ymax=\zmax,hide axis,
             tick label style={font=\small},y dir=reverse]
         \addplot [forget plot] graphics 
         [xmin=\xmin,xmax=\xmax,ymin=\zmin,ymax=\zmax] {{\modelfile-3}.png};
         \addplot[dashed,black,line width=1pt] 
       coordinates {(\xmin+.01,\zmin+.01) (\xmin+.01,\zmax-.01) 
                    (\xmax-.01,\zmax-.01) (\xmax-.01,\zmin+.01) (\xmin+.01,\zmin+.01)};
\end{axis}

\\

\begin{axis}[anchor=center, width=.7\modelwidth, height=.51\modelwidth,
             xmin=\xmin, xmax=\xmax,ymin=\ymin, ymax=\ymax,hide axis,
             tick label style={font=\small}, y dir=reverse]
            \addplot [forget plot] graphics 
                     [xmin=\xmin,xmax=\xmax,ymin=\ymin,ymax=\ymax] {{\modelfile-2}.png};
            
            \addplot[dashed,black,line width=1pt] 
                  coordinates {(\xmin+.01,\ymin+.01) (\xmin+.01,\ymax-.01) 
                               (\xmax-.01,\ymax-.01) (\xmax-.01,\ymin+.01) (\xmin+.01,\ymin+.01)};
\end{axis}

&

\begin{axis}[yshift=-1.80em,
     anchor=south,yshift=.5cm,hide axis,height=.5cm,colorbar/width=.2cm,
     colormap/jet,colorbar horizontal,
     colorbar style={separate axis lines,
                     title style={yshift=-3mm, xshift=0mm},
                     tick label style={font=\small},
                     title={\small{density (\num{e3}\si{\kg\per\meter\cubed})}}},
                     point meta min=\vmin, point meta max=\vmax,
                     height=.5\modelwidth]
\end{axis}
\\};

\draw[black,-stealth,->,dashed,line width=1pt] (-0.2, 3.3) to[out=0,in=-180] ( 0.6, 3.3);
\draw[black,-stealth,->,dashed,line width=1pt] (-0.6, 0.0) to[out=0,in=-180] ( 1.0, 0.0);
\draw[black,-stealth,->,dashed,line width=1pt] (-1.75,-0.15) to[out=-90,in=90] (-1.75,-1.6);

\end{tikzpicture}}
            \hspace*{-0.5em}
  \renewcommand{\modelfile}{Q-true_scale150-850}
  \subfloat[][Quality factor at \num{300} \si{\kilo\Hz}.]
           {
\begin{tikzpicture}
\pgfmathsetmacro{\xmin} {0.}
\pgfmathsetmacro{\xmax} {10.0}
\pgfmathsetmacro{\ymin} {0.}
\pgfmathsetmacro{\ymax} {16.0}
\pgfmathsetmacro{\zmin} {0.}
\pgfmathsetmacro{\zmax} {16.0}

\pgfmathsetmacro{\ycut} { 8.00}
\pgfmathsetmacro{\xcut} { 7.54}
\pgfmathsetmacro{\zcut} {10.80}
\pgfmathsetmacro{\vmin} {150}
\pgfmathsetmacro{\vmax} {850}

\matrix[column sep=-1.8em, row sep=-1.6em] {
\begin{axis}[yshift=0.7*\jumpvert,
  anchor=center,
  tick label style={font=\small},
  grid=both,minor tick num=1,
  xlabel={\scriptsize{$x$ (\si{\cm})}},
  ylabel={\scriptsize{$y$ (\si{\cm})}},
  zlabel={\scriptsize{$z$ (\si{\cm})}}, 
  ztick pos=left,
  3d box,width=1.4\modelwidth, 
  xmin=\xmin,ymin=\ymin,zmin=\zmin,xmax=\xmax,ymax=\ymax,zmax=\zmax,
  every axis x label/.style={at={(0.05, 0.10)},anchor=north},
  every axis y label/.style={at={(1.00, 0.15)},anchor=north}, 
  z label style={xshift= 0.5cm, yshift=-0.50cm},
  xtick={0,10},
  line width=.25pt,
  label style={font=\scriptsize},
  tick label style={font=\scriptsize},
  legend style={font=\scriptsize\selectfont},  
  ]
  \addplot3[fill=white] graphics[points={
            (\xcut,\ymin,\zmin)     => ( 16,215-52)
            (\xcut,\ymax,\zmin)     => (152,215-0)
            (\xcut,\ymin,\zmax)     => ( 16,215-208)
            (\xmin,\ycut,\zmax)     => (132,215-215)}
            ]{{\modelfile_3d}.png};
  
  \addplot3[dashed,black,line width=1pt] 
      coordinates {(\xmin,\ycut,\zmin) (\xmin,\ycut,\zmax)
                   (\xmax,\ycut,\zmax) (\xmax,\ycut,\zmin) (\xmin,\ycut,\zmin)};
  \addplot3[dashed,black,line width=1pt] 
      coordinates {(\xmin,\ymin,\zcut) (\xmin,\ymax,\zcut) 
                   (\xmax,\ymax,\zcut) (\xmax,\ymin,\zcut) (\xmin,\ymin,\zcut)};

  \addplot3[dashed,black,line width=1pt] 
      coordinates {(\xcut,\ymin,\zmin) (\xcut,\ymin,\zmax)
                   (\xcut,\ymax,\zmax) (\xcut,\ymax,\zmin) (\xcut,\ymin,\zmin)};

\end{axis}

&

\begin{axis}[yshift=0.1em,anchor=center,  
             width=0.51\modelwidth, height=0.64\modelheight,
             xmin=\xmin, xmax=\xmax,ymin=\zmin, ymax=\zmax,hide axis,
             tick label style={font=\small},y dir=reverse]
         \addplot [forget plot] graphics 
         [xmin=\xmin,xmax=\xmax,ymin=\zmin,ymax=\zmax] {{\modelfile-1}.png};
         \addplot[dashed,black,line width=1pt] 
       coordinates {(\xmin+.01,\zmin+.01) (\xmin+.01,\zmax-.01) 
                    (\xmax-.01,\zmax-.01) (\xmax-.01,\zmin+.01) (\xmin+.01,\zmin+.01)};
\end{axis}

\begin{axis}[yshift=\jumpvert,
             anchor=center,  width =0.64\modelwidth, 
                             height=0.64\modelheight,
             xmin=\xmin, xmax=\xmax,ymin=\zmin, ymax=\zmax,hide axis,
             tick label style={font=\small},y dir=reverse]
         \addplot [forget plot] graphics 
         [xmin=\xmin,xmax=\xmax,ymin=\zmin,ymax=\zmax] {{\modelfile-3}.png};
         \addplot[dashed,black,line width=1pt] 
       coordinates {(\xmin+.01,\zmin+.01) (\xmin+.01,\zmax-.01) 
                    (\xmax-.01,\zmax-.01) (\xmax-.01,\zmin+.01) (\xmin+.01,\zmin+.01)};
\end{axis}

\\

\begin{axis}[anchor=center, width=.7\modelwidth, height=.51\modelwidth,
             xmin=\xmin, xmax=\xmax,ymin=\ymin, ymax=\ymax,hide axis,
             tick label style={font=\small}, y dir=reverse]
            \addplot [forget plot] graphics 
                     [xmin=\xmin,xmax=\xmax,ymin=\ymin,ymax=\ymax] {{\modelfile-2}.png};
            
            \addplot[dashed,black,line width=1pt] 
                  coordinates {(\xmin+.01,\ymin+.01) (\xmin+.01,\ymax-.01) 
                               (\xmax-.01,\ymax-.01) (\xmax-.01,\ymin+.01) (\xmin+.01,\ymin+.01)};
\end{axis}

&

\begin{axis}[yshift=-1.80em,
     anchor=south,yshift=.5cm,hide axis,height=.5cm,colorbar/width=.2cm,
     colormap/jet,colorbar horizontal,
     colorbar style={separate axis lines,
                     title style={yshift=-3mm, xshift=0mm},
                     tick label style={font=\small},
                     title={\small{Quality factor}}},
                     point meta min=\vmin, point meta max=\vmax,
                     height=.5\modelwidth]
\end{axis}
\\};

\draw[black,-stealth,->,dashed,line width=1pt] (-0.2, 3.3) to[out=0,in=-180] ( 0.6, 3.3);
\draw[black,-stealth,->,dashed,line width=1pt] (-0.6, 0.0) to[out=0,in=-180] ( 1.0, 0.0);
\draw[black,-stealth,->,dashed,line width=1pt] (-1.75,-0.15) to[out=-90,in=90] (-1.75,-1.6);

\end{tikzpicture}}
  \caption{3D breast sample in a bounding domain of size 
           \num{10} $\times$ \num{10} $\times$ \num{16} 
           \si{\centi\meter\cubed}. The values of the 
           parameters depending on the type of tissues 
           are given in \cref{table:interval}.}
  \label{fig:3d:true}
\end{figure}
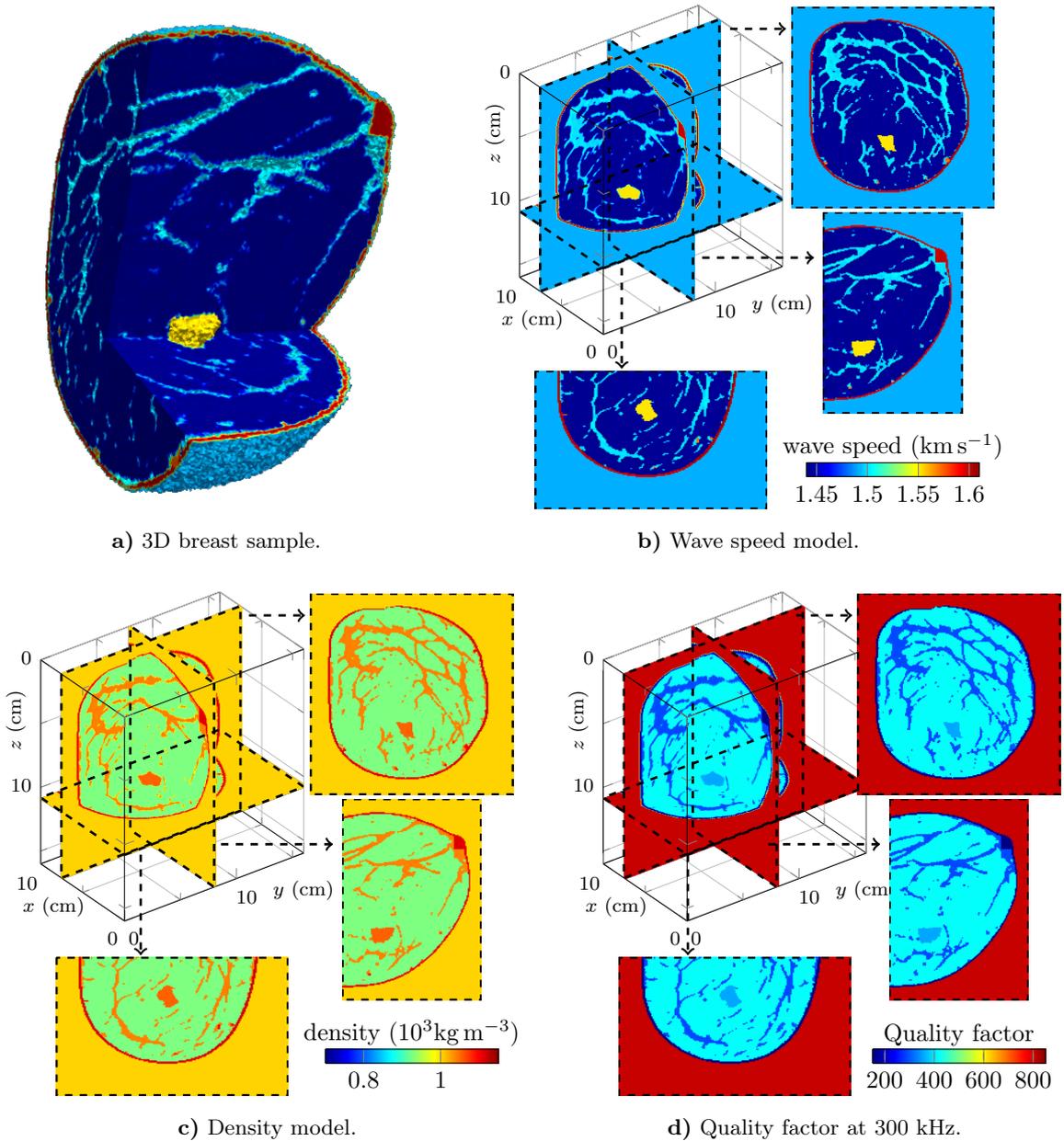

The data are generated from the sides of this 
bounding box, except for the plane where 
$x=\num{10}$ \si{\cm}, which correspond to the 
patient's chest in \cref{fig:3d:true}. That is, 
the data are generated from the five remaining sides.
We consider \num{224} point-sources and \num{1533} 
receivers equally partitioned between the different 
sides, with the receivers measuring the pressure 
field for each of the sources.
Note that the choice of using a rectangular acquisition 
setup is only motivated per simplicity, and we could also 
employ a spherical acquisition, only requiring to adapt the 
discretization mesh.

For the reconstruction, we start with constant 
background parameters with
$c_0=\num{1490}$ \si{\meter\per\second},
$\rho=\num{1000}$ \si{\kg\per\meter\cubed}
and constant attenuation parameters chosen
such that $Q(300\si{\kilo\Hz}) = \num{800}$.
Furthermore, we only invert for the bulk modulus
and keep the density and quality factor at their 
initial values during the entire iterations of FWI.
We also incorporate attenuation model error in 
the experiment: the data are generated using
the simplified Kolsky--Futterman model while the
inversion is carried out with wave propagation 
using the Kelvin--Voigt attenuation model.

\subsection{Reconstruction with absorbing boundary conditions (free-space)}

In the free-space problem, absorbing boundary conditions 
\cref{eq:bc:abc} are implemented on all sides of surrounding box.
In this case, following the observations of \cref{section:fwi-2d},
we expect that only the Fourier frequencies would be needed (i.e., 
taking $\omegaI=0$ to transform time-domain signals in 
\cref{fig:time-domain-acquisition-cx-freq}), and we select the 
following three frequencies $\omegaR/(2\pi)$ for the reconstruction: 
$\{\num{100}, \,\num{200}, \,\num{300}\}$ \si{\kilo\Hz}.
We perform \num{20} minimization iterations per frequency, 
hence a total of \num{60} iterations. 
The wave-speed (assembled from the reconstructed bulk modulus)
is shown in \cref{fig:3d:fwi-absorbing}.

\graphicspath{{figures/tumor_3d/fwi_absorbing/}}
\setlength{\modelwidth} {8cm}
\setlength{\modelheight}{8cm}
\setlength{\jumpvert}   {3.75cm}
\begin{figure}[ht!]\centering
  \renewcommand{\modelfile}{cp_data-kf_fwi-kv_300kHz_gauss05}  
\begin{tikzpicture}
\pgfmathsetmacro{\xmin} {0.}
\pgfmathsetmacro{\xmax} {10.0}
\pgfmathsetmacro{\ymin} {0.}
\pgfmathsetmacro{\ymax} {16.0}
\pgfmathsetmacro{\zmin} {0.}
\pgfmathsetmacro{\zmax} {16.0}

\pgfmathsetmacro{\ycut} { 8.00}
\pgfmathsetmacro{\xcut} { 7.54}
\pgfmathsetmacro{\zcut} {10.80}
\pgfmathsetmacro{\vmin} {1.440}
\pgfmathsetmacro{\vmax} {1.610}

\matrix[column sep=-0.80em, row sep=-2.75em] {
\begin{axis}[yshift=0.7*\jumpvert,
  anchor=center,
  tick label style={font=\small},
  grid=both,minor tick num=1,
  xlabel={\scriptsize{$x$ (\si{\cm})}},
  ylabel={\scriptsize{$y$ (\si{\cm})}},
  zlabel={\scriptsize{$z$ (\si{\cm})}}, 
  ztick pos=left,
  3d box,width=1.4\modelwidth, 
  xmin=\xmin,ymin=\ymin,zmin=\zmin,xmax=\xmax,ymax=\ymax,zmax=\zmax,
  every axis x label/.style={at={(0.05, 0.10)},anchor=north},
  every axis y label/.style={at={(1.00, 0.15)},anchor=north}, 
  z label style={xshift= 0.5cm, yshift=-0.50cm},
  xtick={0,10},
  line width=.25pt,
  label style={font=\scriptsize},
  tick label style={font=\scriptsize},
  legend style={font=\scriptsize\selectfont},  
  ]
  \addplot3[fill=white] graphics[points={
            (\xcut,\ymin,\zmin)     => ( 16,215-52)
            (\xcut,\ymax,\zmin)     => (152,215-0)
            (\xcut,\ymin,\zmax)     => ( 16,215-208)
            (\xmin,\ycut,\zmax)     => (132,215-215)}
            ]{{\modelfile_3d}.png};
  
  \addplot3[dashed,black,line width=1pt] 
      coordinates {(\xmin,\ycut,\zmin) (\xmin,\ycut,\zmax)
                   (\xmax,\ycut,\zmax) (\xmax,\ycut,\zmin) (\xmin,\ycut,\zmin)};
  \addplot3[dashed,black,line width=1pt] 
      coordinates {(\xmin,\ymin,\zcut) (\xmin,\ymax,\zcut) 
                   (\xmax,\ymax,\zcut) (\xmax,\ymin,\zcut) (\xmin,\ymin,\zcut)};

  \addplot3[dashed,black,line width=1pt] 
      coordinates {(\xcut,\ymin,\zmin) (\xcut,\ymin,\zmax)
                   (\xcut,\ymax,\zmax) (\xcut,\ymax,\zmin) (\xcut,\ymin,\zmin)};

\end{axis}

&

\begin{axis}[yshift=0.1em,anchor=center,  
             width=0.51\modelwidth, height=0.64\modelheight,
             xmin=\xmin, xmax=\xmax,ymin=\zmin, ymax=\zmax,hide axis,
             tick label style={font=\small},y dir=reverse]
         \addplot [forget plot] graphics 
         [xmin=\xmin,xmax=\xmax,ymin=\zmin,ymax=\zmax] {{\modelfile-1}.png};
         \addplot[dashed,black,line width=1pt] 
       coordinates {(\xmin+.01,\zmin+.01) (\xmin+.01,\zmax-.01) 
                    (\xmax-.01,\zmax-.01) (\xmax-.01,\zmin+.01) (\xmin+.01,\zmin+.01)};
\end{axis}

\begin{axis}[yshift=\jumpvert,
             anchor=center,  width =0.64\modelwidth, 
                             height=0.64\modelheight,
             xmin=\xmin, xmax=\xmax,ymin=\zmin, ymax=\zmax,hide axis,
             tick label style={font=\small},y dir=reverse]
         \addplot [forget plot] graphics 
         [xmin=\xmin,xmax=\xmax,ymin=\zmin,ymax=\zmax] {{\modelfile-3}.png};
         \addplot[dashed,black,line width=1pt] 
       coordinates {(\xmin+.01,\zmin+.01) (\xmin+.01,\zmax-.01) 
                    (\xmax-.01,\zmax-.01) (\xmax-.01,\zmin+.01) (\xmin+.01,\zmin+.01)};
\end{axis}

\\

\begin{axis}[anchor=center, width=.7\modelwidth, height=.51\modelwidth,
             xmin=\xmin, xmax=\xmax,ymin=\ymin, ymax=\ymax,hide axis,
             tick label style={font=\small}, y dir=reverse]
            \addplot [forget plot] graphics 
                     [xmin=\xmin,xmax=\xmax,ymin=\ymin,ymax=\ymax] {{\modelfile-2}.png};
            
            \addplot[dashed,black,line width=1pt] 
                  coordinates {(\xmin+.01,\ymin+.01) (\xmin+.01,\ymax-.01) 
                               (\xmax-.01,\ymax-.01) (\xmax-.01,\ymin+.01) (\xmin+.01,\ymin+.01)};
\end{axis}

&

\begin{axis}[yshift=-1.80em,
     anchor=south,yshift=.5cm,hide axis,height=.5cm,colorbar/width=.2cm,
     colormap/jet,colorbar horizontal,
     colorbar style={separate axis lines,
                     title style={yshift=-3mm, xshift=0mm},
                     tick label style={font=\small},
                     title={\small{wave speed (\si{\km\per\second})}}},
                     point meta min=\vmin, point meta max=\vmax,
                     height=.5\modelwidth]
\end{axis}
\\};

\draw[black,-stealth,->,dashed,line width=1pt] (-0.2, 4.0) to[out=0,in=-180] ( 0.8, 4.0);
\draw[black,-stealth,->,dashed,line width=1pt] (-0.6, 0.0) to[out=0,in=-180] ( 1.3, 0.0);
\draw[black,-stealth,->,dashed,line width=1pt] (-1.75,-0.15) to[out=-90,in=90] (-1.75,-1.6);

\end{tikzpicture} 
  \caption{Three-dimensional reconstruction of the breast 
           model \cref{fig:3d:true} assuming free-space 
           propagation. 
           The reconstruction starts from constant 
           background parameters. 
           The iterative minimization uses three frequencies from 
           \num{100} to \num{300} \si{\kilo\Hz}, with a total of
           \num{60} iterations.}
  \label{fig:3d:fwi-absorbing}
\end{figure}
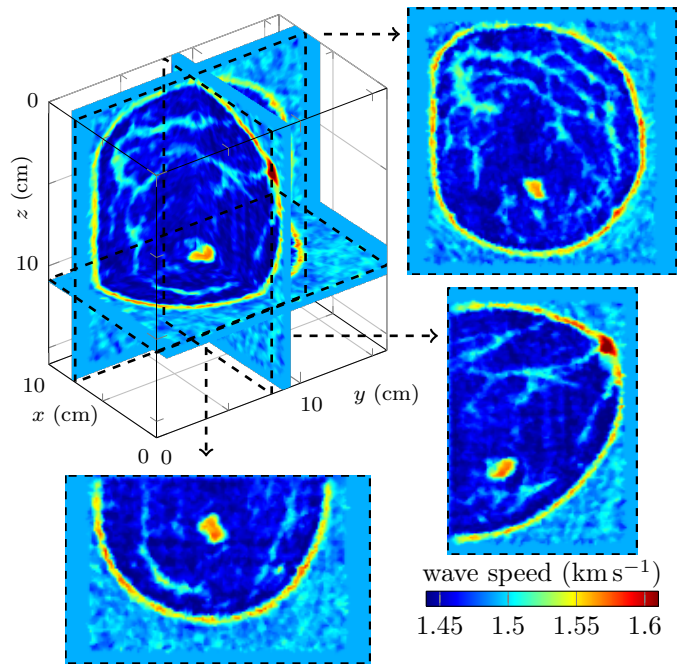

We see that the main features of the sample appear 
correctly on the reconstruction, with the appropriate
values of speed recovered.
The contrasting shape appears appropriately, in its correct 
location, as seen on the different cross-sections
of \cref{fig:3d:fwi-absorbing}.
We also observe some background oscillatory artifacts 
in the surroundings of the medium but the skin contour 
of the breast is well identified. 
Here, smaller details and tissue features are missing and may
require higher frequency-contents in the data (if available), 
to refine the accuracy of the reconstructed model.

\subsection{Reconstruction with wall boundary conditions and complex frequencies}

We consider wall boundary conditions \cref{eq:bc:wall} 
on the sides of the acquisition box. Here, we maintain an 
absorbing condition on the side of the chest (the plane 
at fixed $x=\num{10}$
\si{\cm} in \cref{fig:3d:true}) and use wall conditions on 
all of the other sides.
To alleviate the difficulties of handling the multiple 
reflections generated from the different sides, we use 
complex frequencies. 
The Fourier content remains the same as the one used above, 
with $\omegaR/(2\pi)=\{\num{100}, \,\num{200}, \,\num{300}\}$ 
\si{\kilo\Hz}, and we further incorporate two damping 
coefficients: $\omegaI=\{\num{e4}, \,\num{5e3}\}$ \si{\per\second}.
The sequential progression of frequency follows a 
low-to-high Fourier content and  high-to-low Laplace one, 
\cite{Faucher2020basins}, such that we use
$(\omegaR/(2\pi),   \, \omegaI) = \{
 (100\si{\kilo\Hz}, \, \num{  e4}), 
 (100\si{\kilo\Hz}, \, \num{ 5e3}),
 (200\si{\kilo\Hz}, \, \num{  e4}), \ldots \}$.

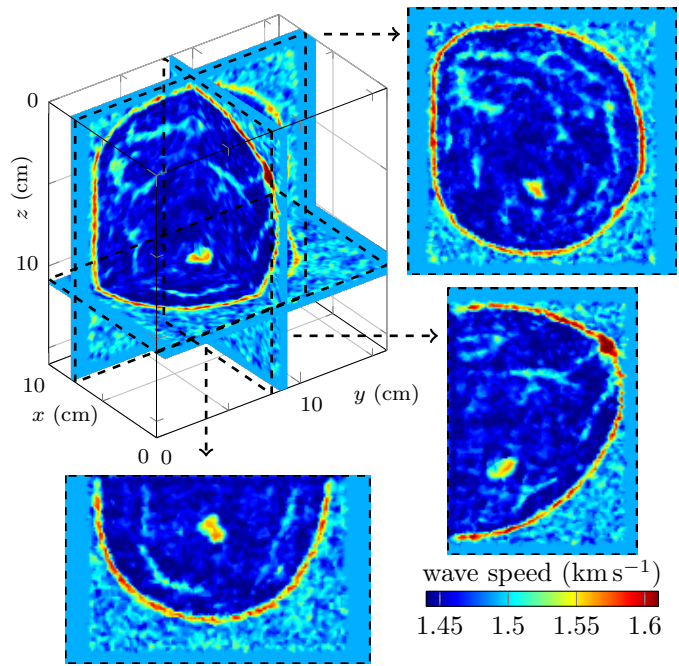
\begin{figure}[ht!]\centering
  \renewcommand{\modelfile}{cp_data-kf_fwi-kv-using2cx_300kHz-5e3_gauss0.5}
\begin{tikzpicture}
\pgfmathsetmacro{\xmin} {0.}
\pgfmathsetmacro{\xmax} {10.0}
\pgfmathsetmacro{\ymin} {0.}
\pgfmathsetmacro{\ymax} {16.0}
\pgfmathsetmacro{\zmin} {0.}
\pgfmathsetmacro{\zmax} {16.0}

\pgfmathsetmacro{\ycut} { 8.00}
\pgfmathsetmacro{\xcut} { 7.54}
\pgfmathsetmacro{\zcut} {10.80}
\pgfmathsetmacro{\vmin} {1.440}
\pgfmathsetmacro{\vmax} {1.610}

\matrix[column sep=-0.80em, row sep=-2.75em] {
\begin{axis}[yshift=0.7*\jumpvert,
  anchor=center,
  tick label style={font=\small},
  grid=both,minor tick num=1,
  xlabel={\scriptsize{$x$ (\si{\cm})}},
  ylabel={\scriptsize{$y$ (\si{\cm})}},
  zlabel={\scriptsize{$z$ (\si{\cm})}}, 
  ztick pos=left,
  3d box,width=1.4\modelwidth, 
  xmin=\xmin,ymin=\ymin,zmin=\zmin,xmax=\xmax,ymax=\ymax,zmax=\zmax,
  every axis x label/.style={at={(0.05, 0.10)},anchor=north},
  every axis y label/.style={at={(1.00, 0.15)},anchor=north}, 
  z label style={xshift= 0.5cm, yshift=-0.50cm},
  xtick={0,10},
  line width=.25pt,
  label style={font=\scriptsize},
  tick label style={font=\scriptsize},
  legend style={font=\scriptsize\selectfont},  
  ]
  \addplot3[fill=white] graphics[points={
            (\xcut,\ymin,\zmin)     => ( 16,215-52)
            (\xcut,\ymax,\zmin)     => (152,215-0)
            (\xcut,\ymin,\zmax)     => ( 16,215-208)
            (\xmin,\ycut,\zmax)     => (132,215-215)}
            ]{{\modelfile_3d}.png};
  
  \addplot3[dashed,black,line width=1pt] 
      coordinates {(\xmin,\ycut,\zmin) (\xmin,\ycut,\zmax)
                   (\xmax,\ycut,\zmax) (\xmax,\ycut,\zmin) (\xmin,\ycut,\zmin)};
  \addplot3[dashed,black,line width=1pt] 
      coordinates {(\xmin,\ymin,\zcut) (\xmin,\ymax,\zcut) 
                   (\xmax,\ymax,\zcut) (\xmax,\ymin,\zcut) (\xmin,\ymin,\zcut)};

  \addplot3[dashed,black,line width=1pt] 
      coordinates {(\xcut,\ymin,\zmin) (\xcut,\ymin,\zmax)
                   (\xcut,\ymax,\zmax) (\xcut,\ymax,\zmin) (\xcut,\ymin,\zmin)};

\end{axis}

&

\begin{axis}[yshift=0.1em,anchor=center,  
             width=0.51\modelwidth, height=0.64\modelheight,
             xmin=\xmin, xmax=\xmax,ymin=\zmin, ymax=\zmax,hide axis,
             tick label style={font=\small},y dir=reverse]
         \addplot [forget plot] graphics 
         [xmin=\xmin,xmax=\xmax,ymin=\zmin,ymax=\zmax] {{\modelfile-1}.png};
         \addplot[dashed,black,line width=1pt] 
       coordinates {(\xmin+.01,\zmin+.01) (\xmin+.01,\zmax-.01) 
                    (\xmax-.01,\zmax-.01) (\xmax-.01,\zmin+.01) (\xmin+.01,\zmin+.01)};
\end{axis}

\begin{axis}[yshift=\jumpvert,
             anchor=center,  width =0.64\modelwidth, 
                             height=0.64\modelheight,
             xmin=\xmin, xmax=\xmax,ymin=\zmin, ymax=\zmax,hide axis,
             tick label style={font=\small},y dir=reverse]
         \addplot [forget plot] graphics 
         [xmin=\xmin,xmax=\xmax,ymin=\zmin,ymax=\zmax] {{\modelfile-3}.png};
         \addplot[dashed,black,line width=1pt] 
       coordinates {(\xmin+.01,\zmin+.01) (\xmin+.01,\zmax-.01) 
                    (\xmax-.01,\zmax-.01) (\xmax-.01,\zmin+.01) (\xmin+.01,\zmin+.01)};
\end{axis}

\\

\begin{axis}[anchor=center, width=.7\modelwidth, height=.51\modelwidth,
             xmin=\xmin, xmax=\xmax,ymin=\ymin, ymax=\ymax,hide axis,
             tick label style={font=\small}, y dir=reverse]
            \addplot [forget plot] graphics 
                     [xmin=\xmin,xmax=\xmax,ymin=\ymin,ymax=\ymax] {{\modelfile-2}.png};
            
            \addplot[dashed,black,line width=1pt] 
                  coordinates {(\xmin+.01,\ymin+.01) (\xmin+.01,\ymax-.01) 
                               (\xmax-.01,\ymax-.01) (\xmax-.01,\ymin+.01) (\xmin+.01,\ymin+.01)};
\end{axis}

&

\begin{axis}[yshift=-1.80em,
     anchor=south,yshift=.5cm,hide axis,height=.5cm,colorbar/width=.2cm,
     colormap/jet,colorbar horizontal,
     colorbar style={separate axis lines,
                     title style={yshift=-3mm, xshift=0mm},
                     tick label style={font=\small},
                     title={\small{wave speed (\si{\km\per\second})}}},
                     point meta min=\vmin, point meta max=\vmax,
                     height=.5\modelwidth]
\end{axis}
\\};

\draw[black,-stealth,->,dashed,line width=1pt] (-0.2, 4.0) to[out=0,in=-180] ( 0.8, 4.0);
\draw[black,-stealth,->,dashed,line width=1pt] (-0.6, 0.0) to[out=0,in=-180] ( 1.3, 0.0);
\draw[black,-stealth,->,dashed,line width=1pt] (-1.75,-0.15) to[out=-90,in=90] (-1.75,-1.6);

\end{tikzpicture}
  \caption{Three-dimensional reconstruction of the breast 
           model \cref{fig:3d:true} assuming wall boundary
           on the sides of the acquisition box.
           The reconstruction starts from constant 
           background parameters and only the bulk modulus 
           is inverted while the density and quality factor 
           remain constant. 
           The iterative minimization uses six (complex)
           frequencies with $\omegaR/(2\pi)=\{
           \num{100}, \num{200}, \num{300}\}$ \si{\kilo\Hz}
           and $\omegaI=\{\num{e4}, \,\num{5e3}\}$ \si{\per\second},
           for a total of \num{120} iterations.}
  \label{fig:3d:fwi-wall}
\end{figure}

The reconstructed wave speed is pictured 
in \cref{fig:3d:fwi-wall}, where we see 
that we obtain a accuracy relatively similar 
to the case assuming the free-space 
propagation (\cref{fig:3d:fwi-absorbing}).
The main tissue features
are visible in the cross-sections and the skin
layer is recovered. In addition, the position, size and amplitude
of the contrast are obtained and clearly stand out on the reconstruction.
Similar to the free-space medium reconstruction, we observe some
background perturbations, which however does not prevent us 
from clearly identifying the main components of the sample.

\graphicspath{{figures/tumor_3d/fwi_wall/}}

\subsection{Computational cost}
\label{subsection:computational-cost}

The reconstruction relies on the iterative algorithm depicted
in \cref{algo:FWI}, where the computational cost comes from 
the resolution of the forward problems which has to be repeated 
with iterations. 
For the resolution of the time-harmonic problem, the discretization 
leads to relatively large matrix, and one further needs solve the 
resulting linear system for each of the source in the acquisition.
In our implementation, we made the following choices regarding the
computational implementation:
\begin{itemize}
  \item We use HDG discretization method, see \cref{subsection:hdg}, 
        which can help reduce the computational cost compared to
        other approaches, see \cite{Faucher2020adjoint}.
  \item Using HDG, the order of the approximation polynomial on each 
        cell of the discretization mesh is chosen independently depending
        on the local (to the cell) wavelength, to ensure the linear 
        system is as small as possible, cf. \cite{Faucher2020adjoint}.
  \item For the resolution of the linear systems, we use the direct solver 
        MUMPS, \cite{Amestoy2019}, that is particularly efficient for solving 
        systems with multiple right-hand sides, contrary to iterative solvers. 
        Namely, once the global matrix is factorized, the resolution for all of 
        the sources (i.e., \num{224} point-sources in our 3D experiments) 
        is fast. On the other hand, the matrix factorization is a computationally
        intensive step, especially regarding the memory consumption.
\end{itemize}

We use software \texttt{hawen} (see \cite{Hawen2021} 
and \cref{footnote:haven}) in all of the experiments. 
They are carried out on the cluster PlaFRIM\footnote{\url{https://www.plafrim.fr/}.}.
The computational costs are the following:
\begin{itemize}
  \item The 2D experiments of \cref{section:fwi-2d} are performed 
        on \num{36} physical cores, with \num{18} processors 
        and \num{2} threads per processors. 
        The computational time to obtain a reconstruction as in 
        \cref{fig:fwi2d:abosrbing:all} is about 10 \si{\min} 
        and corresponds to 120 inversion iterations.
        The memory required to factorize the matrix is about
        1\si{GiB} for all frequencies.
  \item The 3D experiments are performed on \num{540} physical 
        cores, with \num{90} processors 
        and \num{6} threads per processors. 
        The computational time to obtain a reconstruction as in 
        \cref{fig:3d:fwi-absorbing} is of 5 \si{\hour} 
        and corresponds to 60 inversion iterations.
        The memory required to factorize the matrix depends on
        the frequency (high frequency needs higher polynomial
        order hence generates larger systems), and varies from
        \num{200}\si{GiB} at \num{100}\si{\kilo\Hz} to 
        \num{800}\si{GiB} at \num{300}\si{\kilo\Hz}.
\end{itemize}

\section{Conclusion}

We have performed imaging of visco-acoustic 
media using an iterative reconstruction procedure.
Attenuation can be encoded in the propagation 
via different models, and we have the following results.
\begin{enumerate}
  \item We have implemented seven different models of attenuation 
             and compare the wave propagation. It highlights that, while
             they can coincide at a reference frequency, each model leads
             to different wave patterns.
  \item We have shown that our algorithm is robust with respect to 
attenuation model uncertainty, where we have used a relatively 
narrow band of frequency (here between \num{100} and \num{600} \si{\kilo\Hz}).
Namely, if one does not know the attenuation model that 
corresponds to the sample, it does not prevent from 
reconstructing its main properties. 
  \item While the wave sped and bulk modulus are accurately
reconstructed and allow to describe the key-features 
of the sample, the recovery of the density and quality
factor is more difficult, cf.~\cite{Virieux2009,Karaouglu2017} and the references therein. 
To recover them, one would need a parametrization 
that is more suitable to focus on amplitude variation,
with an appropriate cost function (i.e., instead of the
$L2$ norm), \cite{Karaouglu2017}.
  \item We have shown that the reconstruction assuming a propagation
in free-space performs well but considering wall boundary 
around the samples leads to multiple reflections that drastically
alter the reconstruction. As an alternative, we have proposed 
the use of complex frequencies which, by introducing an 
artificial damping, alleviate the difficulties and enable
the accurate discoveries of the sample features. 
It provides us with a candidate to further consider media 
containing objects of high contrast that generate strong 
reflections, such as bones and skull, which are known to be 
harder to image. 
This is the subject of ongoing studies.
\end{enumerate}
Our experiments use physical properties of tissues which, 
in this range of ultrasonic frequencies, are weakly attenuating (quality 
factor higher than 100). 
In the opposite situation where the medium is highly attenuating and 
prevent waves to propagate in the tissues with sufficient 
energy, the strategy we propose may have to be modified, this is part
of ongoing works.
In addition, the consideration of visco-elasticity is 
part of our ongoing research, where the difficulties 
come on the one hand from the incorporation of additional 
unknowns for inversion and on the other hand by the 
increased computational cost of solving a vector-wave
problem.

\section*{Acknowledgments}
 FF is funded by the Austrian Science 
 Fund (FWF) under the Lise Meitner 
 fellowship M 2791-N.
 For the experiments, we acknowledge 
 the use of the cluster PlaFRIM\footnote{supported by 
 Inria, CNRS (LABRI and IMB), Universit\'e de Bordeaux, 
 Bordeaux INP and Conseil R\'egional d'Aquitaine,
 see \url{https://www.plafrim.fr/}.}, and  of the 
 Vienna Scientific Cluster \texttt{VSC4}\footnote{\url{https://vsc.ac.at/}.}. 
 OS is supported by the Austrian
 Science Fund (FWF) with SFB F68, project F6807-N36 (Tomography with
 Uncertainties).
 The financial support by the Austrian Federal Ministry for Digital and
 Economic Affairs, the National Foundation for Research, Technology and
 Development and the Christian Doppler Research Association is gratefully
 acknowledged.

\appendix

\section{Time-harmonic formulation with complex frequencies}
\label{appendix:time-harmonic}
To illustrate the derivation of the frequency-domain equation \cref{eq:euler_main}
for complex frequencies, we start with the time-domain wave equation. 
For simplicity, we consider the case without attenuation to avoid 
possibly integro-differential equation. In this case, time-dependent acoustic
waves are given by the velocity and pressure fields, respectively 
$\widehat{\velocity}(\bx,t)$ and $\widehat{\pressure}(\bx,t)$, solutions to
\begin{subequations} \label{eq:euler-main-time}
\begin{empheq}[left={\empheqlbrace}]{align} 
   \rho(\bx) \,\partial_t \, \widehat{\velocity}(\bx,t) \,+\, \nabla \widehat{\pressure}(\bx,t) & \,=\, 0,  \\
   \dfrac{1}{\kappa_0(\bx)} \, \partial_t \widehat{\pressure}(\bx,t) \,+\, 
   \nabla \cdot \widehat{\velocity}(\bx,t) & \, = \, \widehat{g}(\bx,t).
\end{empheq} \end{subequations}
where $\widehat{g}$ is the time-dependent source, that is, the Ricker wavelet 
in the context of \cref{fig:time-domain-acquisition}.

We consider time-harmonic solutions of the following form:
\begin{subequations}\begin{align}
   \widehat{\pressure}(\bx,t) &\,=\, \pressure(\bx,\omega) \,e^{-\ii \,(\omegaR \,+\, \ii \omegaI) \,t} \,, \\
   \widehat{\velocity}(\bx,t) &\,=\, \velocity(\bx,\omega) \,e^{-\ii \,(\omegaR \,+\, \ii \omegaI) \,t} \,, \\
   \widehat{g}(\bx,t)         &\,=\, g(\bx,\omega) \,e^{-\ii \,(\omegaR \,+\, \ii \omegaI) \,t} \,,
\end{align}\end{subequations}
where we remind that by definition, $\omega := \omegaR \,+\, \ii \omegaI$.
Deriving with respect to time gives, 
\begin{subequations}\begin{align}
   \partial_t \widehat{\pressure}(\bx,t) &\,=\,-\ii\omega\, \pressure(\bx,\omega) \,e^{-\ii \,\omega \,t} \,, \\
   \partial_t \widehat{\velocity}(\bx,t) &\,=\,-\ii\omega\, \velocity(\bx,\omega) \,e^{-\ii \,\omega \,t} \, .
\end{align}\end{subequations}
Replacing in \cref{eq:euler-main-time}, we obtain,
\begin{subequations} 
\begin{empheq}[left={\empheqlbrace}]{align}
   -\ii\omega \, \rho(\bx) \,{\velocity}(\bx,\omega) \,e^{-\ii \,\omega \,t}
   \,+\, \nabla \big( {\pressure}(\bx,\omega)\,e^{-\ii \,\omega \,t} \big) & \,=\, 0,  \\
   \dfrac{-\ii\omega}{\kappa_0(\bx)} \, {\pressure}(\bx,\omega) \,e^{-\ii \,\omega \,t} \,+\, 
   \nabla \cdot \big( \velocity(\bx,\omega)\,e^{-\ii \,\omega \,t} \big) & \, = \, g(\bx,\omega)\,\,e^{-\ii \,\omega \,t}.
\end{empheq} \end{subequations}
Simplifying the terms $e^{-\ii \,\omega \,t}$ (that do not depend on space), 
we obtain the time-harmonic equation \cref{eq:euler_main} in the case without
attenuation.
We note that this transformation introduces two parameters ($\omegaR$ and $\omegaI$),
similarly to the continuous Gabor or the Wavelet transforms.

\scriptsize
\bibliographystyle{siamplain}
\bibliography{sections/bibliography}

\end{document}